\documentclass[12pt]{amsart}
\newenvironment{theorem}{\noindent \bf Theorem.\it}{ \rm }
\newenvironment{proposition}{\noindent \bf Proposition.\it}{ \rm }
\newenvironment{lemma}{\noindent \bf Lemma.\it}{ \rm }
\newenvironment{corollary}{\noindent \bf Corollary.\it}{ \rm }

\newcommand{{\ind}}{\operatorname{index}}
\newcommand{{\ad}}{\operatorname{ad}}
\newcommand{{\rk}}{\operatorname{rank}}
\newcommand{{\id}}{\operatorname{id}}
\newcommand{{\car}}{\operatorname{char}}
\newcommand{{\dist}}{\operatorname{dist}}

\newcommand{{\liea}}{\mathfrak{a}}
\newcommand{{\lieg}}{\mathfrak{g}}
\newcommand{{\liep}}{\mathfrak{p}}
\newcommand{{\lieh}}{\mathfrak{h}}
\newcommand{{\lien}}{\mathfrak{n}}
\newcommand{{\lier}}{\mathfrak{r}}
\newcommand{{\lieb}}{\mathfrak{b}}
\newcommand{{\liesl}}{\mathfrak{sl}}
\newcommand{{\Oint}}{\mathcal{O}_{\operatorname{int}}}
\newcommand{{\bbN}}{\mathbb{N}}
\newcommand{{\bbC}}{\mathbb{C}}
\newcommand{{\wt}}{\operatorname{wt}}
\newcommand{{\Stab}}{\operatorname{Stab}}
\newcommand{{\pitheta}}{\pi*\theta_s^{s'}*\pi'}
\newcommand{{\Pitheta}}{\mathbb{P}_{\lambda}*\theta_s^{s'}*\mathbb{P}_{\mu}}
\newcommand{{\Pithetahat}}{\widehat{\mathbb{P}}_{\lambda}*\theta_s^{s'}*\widehat{\mathbb{P}}_{\mu}}
\newcommand{{\Pithetahatdouble}}{\widehat{\mathbb{P}}_{2\lambda}*\theta_{1/2}^{1/2}*\widehat{\mathbb{P}}_{2\mu}}
\newcommand{{\Pinfty}}{\mathbb{P}_{\infty}}
\newcommand{{\Plambda}}{\mathbb{P}_{\lambda}}
\newcommand{{\Plambdahat}}{\widehat{\mathbb{P}}_{\lambda}}
\newcommand{{\Pmuhat}}{\widehat{\mathbb{P}}_{\mu}}
\begin{document}

\title[Path model for Generalized Kac-Moody algebras]
{A Littelmann path model for crystals of Generalized Kac-Moody algebras.}

\author{Anthony Joseph and Polyxeni Lamprou}
\thanks{Work supported in part by the European Community RTN network ``Liegrits'', Grant No. MRTN - CT -
2003 - 505078 and in part by Minerva Foundation, Germany, Grant No. 8466.}

\maketitle

\begin{abstract}
A Littelmann path model is constructed for crystals pertaining to a not necessarily symmetrizable
Borcherds-Cartan matrix. Here one must overcome several combinatorial problems coming from the
imaginary simple roots. The main results are an isomorphism theorem and a character formula of Borcherds-Kac-Weyl type
for the crystals. In the symmetrizable case, the isomorphism theorem implies that the crystals constructed by
this path model coincide with those of Jeong, Kang, Kashiwara and Shin obtained by taking $q\rightarrow 0$
limit in the quantized enveloping algebra.
\end{abstract}

\section{Introduction}
\subsection{}
The original proof of the Weyl character formula given in 1925 by Weyl following the work of Schur for $\mathfrak{gl}(n)$ underwent a number of simplifications with a
particularly notable one due to Bernstein, Gelfand and Gelfand \cite{BGG}. This proof was shown by Kac
\cite{Kac2} to extend to integrable modules for Kac-Moody algebras obtained from a symmetrizable Cartan
matrix. For affine Lie algebras the corresponding Weyl denominator formula spectacularly recovered and
generalized sum-product identities from number theory due to Fermat, Gauss and Jacobi.

More recently Borcherds \cite{B1} showed that the Kac-Moody theory extends with equally beautiful results
when imaginary simple roots are permitted.  In particular the Bernstein-Gelfand-Gelfand method gives a
character formula, somewhat more complicated than the Weyl-Kac formula for unitarizable highest weight
modules.

\subsection{}
In 1986, Drinfeld and Jimbo independently introduced quantized enveloping algebras involving a parameter $q$. A little later, Lusztig \cite{L} and Kashiwara \cite{K1} considered a $q\rightarrow 0$ limit of these algebras and the integrable modules over them. In the Kashiwara theory $q$ was interpreted as the temperature and the modules were deemed to ``crystallize'' into a simpler form.  Kashiwara, drawing in part an observation of Date, Jimbo and Miwa \cite[Introduction]{K1}, required that the tensor product be without sums or coefficients different from $0$ or $1$,
leading to a tight combinatorial structure. Naturally an integrable highest weight module gives rise to a normal highest weight crystal (which can be viewed as a rather special graph).  Since much structure is lost in the process the latter are not uniquely defined by their highest weights.  However using the tensor structure, one obtains a unique closed (under tensor product) family of normal highest weight crystals.  More recently Jeong, Kang and
Kashiwara \cite{JKK} have extended this theory to include simple imaginary roots (as in Borcherds) but still
with the assumption that the Cartan matrix is symmetrizable (which is needed for quantization).

\subsection{}
Shortly after Kashiwara introduced crystals, Littelmann \cite{L1, L2} found a purely combinatorial path model
for them based on the Cartan matrix which was no longer required to be symmetrizable. He constructed a closed
family of normal highest weight crystals and computed their characters. This was based on Lakshmibai-Seshadri
paths, themselves described by Bruhat sequences in the Weyl group together with an intregrality condition.

\subsection{}
In this paper we extend Littelmann's path model to include imaginary simple roots. Specifically this means proving Proposition \ref{unique} and Theorems \ref{iso}, \ref{embedding} and \ref{character}, which we describe qualitatively below. This involves a number of
combinatorial complications.  Instead of the Weyl group we use a monoid with generators defined by both the
real and the imaginary simple roots. Here the presence of non-invertible elements ultimately means that the
normal highest weight crystals are not strict subcrystals of the full crystal defined by all possible paths.
Besides they are normal only with respect to the real simple roots. This makes it more difficult to show that
``generalized'' Lakshmibai-Seshadri paths describe the required normal highest weight crystals, which is the content of Proposition \ref{unique}. It becomes
correspondingly more difficult to show that this set of crystals is closed with respect to tensor product.
However this being achieved (Theorem \ref{iso}) we recover in the symmetrizable case, the Kashiwara crystals by the crystal embedding theorem \ref{embedding} and the resulting uniqueness (Theorem \ref{Binfty}). Here we remark that the embedding theorem is valid also in the non-symmetrizable case extending thereby the work of Jeong, Kang and Kashiwara \cite{JKK}. Finally we prove (Theorem \ref{character}) a version of Littelmann's combinatorial character formula for the crystals in this family.
Unlike \cite{JKK}, this does not need the Cartan matrix to be symmetrizable, though in any case a very similar formula to that
of Borcherds is obtained.

\subsection{}
Unfortunately Littelmann's combinatorial formula does not recover the Weyl denominator formula (known to hold
in the non-symmetrizable totally real case by independent work of Kumar \cite{Ku} and Mathieu \cite{M}). The
question this entails and many others remain open.\\

\noindent{\bf Acknowledgements.} This work started when the second author was visiting the University of Cologne as a Liegrits predoc. She would like to take the opportunity to thank P. Littelmann for his hospitality and his guidance during her stay.

\section{Preliminaries}

\subsection{Generalized Kac-Moody algebras}
Unless otherwise specified all numerical values are assumed rational. In particular, all vector spaces are
over $\mathbb{Q}$. We denote by $\mathbb{N}$ the set of natural numbers and we set
$\mathbb{N}^+:=\mathbb{N}\setminus \{0\}$.

\subsubsection{} \label{BCmatrix}
Let $I$ be a countable index set. We call $A = (a_{ij})_{i, j\in\ I}$ a {\it Borcherds-Cartan matrix} if the
following are satisfied :
\begin{enumerate}
\item $a_{ii} = 2$ or $a_{ii}\in\ -\mathbb{N}$ for all $i$,
\item $a_{ij}\in\ -\mathbb{N}$, for all $i\neq j$,
\item $a_{ij} = 0$ if and only if $a_{ji} = 0$.
\end{enumerate}
We call an index $i$ real if $a_{ii} = 2$ and we denote by $I^{re}$ the set of real indices. Otherwise, we
call an index $i$ imaginary and we denote by $I^{im} = I\setminus I^{re}$ the set of imaginary indices.

If $I = I^{re}$ and is finite, then $A$ is a generalized Cartan matrix in the language of \cite[Section
1.1]{Kac}. The matrix $A$ is called {\it symmetrizable} if there exists a diagonal matrix $S =
\operatorname{diag}\{s_i\in\ \mathbb{N}^+\,|\,i\in\ I\}$ such that $SA$ is symmetric.

\subsubsection{}\label{notation}
Let $\lieg$ be the generalized Kac-Moody algebra associated to a Borcherds-Cartan matrix $A$, $\lieh$ a fixed
Cartan subalgebra of $\lieg$, $\Pi=\{\alpha_i\,|\,i\in\
I\}\subset \lieh^*$ the set of simple roots, $\Pi^{\vee}=\{\alpha_i^{\vee}\,|\,i\in\ I\}\subset \lieh$ the
set of simple coroots such that $\alpha_i^{\vee}(\alpha_j)=a_{ij}$ and $\Delta$ the root system of $\lieg$ (for more details see \cite{B1, B2}).

\subsubsection{}
Let $P = \{\lambda\in\ \lieh^*\,|\, \alpha_i^{\vee}(\lambda)\in\ \mathbb{Z},\,\mbox{for all}\,\, i\in\ I\}$
be the weight lattice of $\lieg$, $Q = \bigoplus\limits_{i\in\ I}\mathbb{Z}\alpha_i$ be the root lattice and
$Q^+=\bigoplus\limits_{i\in\ I}\mathbb{N}\alpha_i$. Of course $\Delta\subset Q\subset P$. Set $P^+
=\{\lambda\in\ P\,|\,\alpha_i^{\vee}(\lambda)\geq 0, \,\mbox{for all}\, i\in\ I\}$.

\subsubsection{}
Define a partial order in $Q$ by setting $\beta\succ \gamma$ if and only if $\beta -\gamma\in\ Q^+$. Let
$\Delta^+ = \{\beta\in\ \Delta\,|\, \beta\succ 0\}$ be the set of positive roots and $\Delta^-=-\Delta^+$ the
set of negative roots. One has that $\Delta = \Delta^+\sqcup \Delta^-$.

\subsubsection{}\label{monoid}
For all $i\in\ I$ let $r_i$ be the linear map $r_i :\lieh^*\rightarrow \lieh^*$ defined by $$r_i (x) = x
-\alpha_i^{\vee}(x)\alpha_i.$$ Note that $r_i$ is a reflection (and thus $r_i^2=\id$) if and only if $i\in\
I^{re}$. Otherwise, if $i\in\ I^{im}$, $r_i$ has infinite order. Set $T = \langle r_i\,|\,i\in\ I\rangle$ to be the monoid
generated by all the $r_i,\, i\in\ I$ and denote by $\id$ its neutral element. Let $W$ be the group generated
by the reflections $r_i,\, i\in\ I^{re}$ and call it the Weyl group of $\lieg$. Then of course $W$ lies in
$T$. For any $\tau\in\ T$ we may write $\tau = r_{i_1}r_{i_2}\cdots r_{i_{\ell}}$ where $i_j\in\ I$, for all
$j$ with $1\leq j\leq \ell$. We call this a reduced expression if $\ell$ takes its minimal value which we
define to be the reduced length $\ell(\tau)$ of $\tau$.

\subsubsection{}
For all $i\in\ I$, we define $r_i$ on $\lieh$ by :
$$r_i(h) = h - h(\alpha_i)\alpha_i^{\vee}.$$
One checks that if $i\in I^{re}$ then $(r_ih)(r_i\lambda)=h(\lambda)$ for all $h\in\ \lieh$ and all $\lambda\in\ \lieh^*$.

\subsubsection{}\label{Wgeometry}
Set $\mathcal{C}=\{\mu\in\ \lieh^*\,|\,\alpha_i^{\vee}(\mu)\geq 0, \, \mbox{for all}\, i\in\ I^{re}\}$, the
set of dominant elements of $\lieh^*$. (Notice that we consider only real indices). We call a weight in
$\mathcal{C}$ a dominant weight. One has that $P^+\subset \mathcal{C}$. Notice that $-\alpha_i\in\ P^+\subset
\mathcal{C}$ for all $i\in\ I^{im}$. By \cite[Proposition 3.12]{Kac}, for all $\lambda\in\ \mathcal{C}$ one
has $W\lambda\cap \mathcal{C}=\{\lambda\}$ and $w\lambda\in\ \lambda-\mathbb{N}\Delta^+_{re}$.
Choose $\rho\in\ \lieh^*$ such that $\alpha_i^{\vee}(\rho)=\frac{1}{2}a_{ii}$. Then $\rho\in\ \mathcal{C}$,
but in general $\rho\notin P^+$.

\subsubsection{}\label{stabilizer}
Let $\lambda\in\ \mathcal{C}$ and denote by $W_{\lambda}$ the stabilizer of $\lambda$ in $W$. Then
$W_{\lambda}$ is generated by the simple reflections which stabilize $\lambda$, that is
$W_{\lambda}=\langle r_i\in W\,|\,r_i\lambda=\lambda\rangle$. Even when $|I^{re}|=\infty$, the proof is as in \cite[Proposition
3.12]{Kac}.

\subsubsection{}\label{coroot}
Denote by $\Pi_{re}=\{\alpha_i\,|\,i\in\ I^{re}\}$ and by $\Pi_{im}=\{\alpha_i\, |\, i\in\ I^{im}\}$ the sets
of real and imaginary simple roots respectively. \\

\begin{lemma} Take $\alpha_i,\, \alpha_j\in\ \Pi$ and $w, \tilde w\in\ W$. If $w\alpha_i=\tilde w\alpha_j$,
then $w\alpha_i^{\vee}=\tilde w\alpha_j^{\vee}$.
\end{lemma}
\begin{proof}If $\alpha_i,\,\alpha_j\in\ \Pi_{re}$, the assertion obtains from \cite[Section 5.1]{Kac}. Suppose $\alpha_i\in\ \Pi_{im}$. Since $w\alpha_i\in\
\alpha_i+\mathbb{N}\Pi_{re}$, the hypothesis forces $\alpha_i=\alpha_j$. It then suffices to prove the
assertion for $\tilde w=\id$ and $w\in\ \operatorname{Stab}_W(\alpha_i)$. Since $-\alpha_i\in\ \mathcal{C}$,
by section \ref{stabilizer}, we can write $w=r_{i_1}\cdots r_{i_k}$, with $\alpha_{i_t}^{\vee}(\alpha_i)=0$,
for all $t$, with $1\leq t\leq k$. Then $\alpha_i^{\vee}(\alpha_{i_t})=0$, for all $t$, with $1\leq t\leq k$,
so $w\in\ \operatorname{Stab}_W(\alpha_i^{\vee})$, as required.
\end{proof}

\noindent {\bf Definition.} By the above lemma, we may define $\beta^{\vee}\in\ \lieh$, for all $\beta\in\
W\Pi$, through $\beta^{\vee}=w\alpha_i^{\vee}$, given $\beta=w\alpha_i$.

\subsubsection{}\label{delta}
Take $i\in\ I^{re}$. Through \cite[Lemma 3.8]{Kac} we obtain $r_i(\Delta^+\setminus \{\alpha_i\})\subset
\Delta^+\setminus \{\alpha_i\}$. In particular, $\Delta$ is $W$-stable. Call a root $\beta\in\ W\Pi$ real if
$\beta^{\vee}(\beta)=2$ and imaginary if $\beta^{\vee}(\beta)\leq 0$. Set $\Delta_{re} =W\Pi_{re}$ and
$\Delta_{im}=W\Pi_{im}$. Define also $\Delta_{re}^+=\Delta_{re}\cap \Delta^+$, $\Delta^-_{re}=-\Delta_{re}^+$
and notice that $\Delta_{im}\subset \Delta^+$.  In general $\Delta_{re}\sqcup (\Delta_{im}\sqcup
-\Delta_{im})\subset \Delta$ is a strict inclusion. However, its complement in $\Delta$ does not play any
role in our analysis.

\subsubsection{}\label{improperties}
One could roughly say that everything we know about the Weyl group and the real roots in the Kac-Moody case,
also holds for the generalized Kac-Moody algebras. The imaginary roots need some attention. The following
result will be repeatedly used in the sequel.\\

\begin{lemma} Take $i\in\ I^{im}$, then
\begin{enumerate}\label{star}
\item $\beta^{\vee}(\alpha_i)\leq 0$, for all $\beta\in\ \Delta^+_{re}\sqcup \Delta_{im}$,
\item $\alpha_i^{\vee}(\beta)\leq 0$, for all $\beta\in\ Q^+$.
\end{enumerate}
\end{lemma}
\begin{proof}
Indeed, for (1) take $\beta=w\alpha_j$; then $\beta^{\vee}(\alpha_i)=\alpha_j^{\vee}(w^{-1}\alpha_i)$ and
$w^{-1}\alpha_i\in\ \alpha_i+\mathbb{N}\Pi_{re}$, because $-\alpha_i\in\ \mathcal{C}$. Hence the assertion
for $\alpha_j\in\ \Pi_{im}$. For $\alpha_j\in\ \Pi_{re}$, one must show that $(w\alpha_j)^{\vee}\in\
\mathbb{N}\Pi_{re}^{\vee}$. This is stated in \cite[Section 5.1]{Kac}. Finally, (2) is an immediate
consequence of the properties of the matrix $A$.
\end{proof}

\subsection{Dominant elements in $T\lambda$}
In this section we give a characterization of the dominant weights in the $T$-orbit $T\lambda$ of a weight
$\lambda\in\ P^+$.

\subsubsection{}
\begin{lemma}\label{Tlambda} For all $\lambda\in\ P^+$ one has that $T\lambda\subset \lambda -Q^+$. In
particular, $\alpha_i^{\vee}(\mu)\geq 0$ for all $\mu\in\ T\lambda$ and all $i\in\ I^{im}$.
\end{lemma}
\begin{proof} We will prove by induction on $\ell(\tau)$ that
\[\tau \lambda\in\ W\lambda -\mathbb{N}\Delta_{im}.\]
Then since $W\lambda\subset \lambda-\mathbb{N}\Delta_{re}^+$ and
$\Delta_{im}\subset \Delta^+$ by section \ref{delta}, we will have that $T\lambda\subset
\lambda-\mathbb{N}\Delta^+=\lambda-Q^+$.

For $\tau = \id$ the statement is obvious. Let $\tau\lambda=w\lambda-\beta = \lambda-\gamma$, with $\beta\in\
\mathbb{N}\Delta_{im}$ and $\gamma\in\ \mathbb{N}\Delta^+$. Take $i\in\ I^{im}$, then since by lemma
\ref{star} (2), $\alpha_i^{\vee}(\lambda-\gamma)\geq 0$ one has
$$r_i\tau\lambda \in\ \tau\lambda-\mathbb{N}\alpha_i\subset
W\lambda-\mathbb{N}\Delta_{im}.$$ Take $i\in\ I^{re}$. By section \ref{delta} we have that
$r_i\Delta_{im}\subset \Delta_{im}$ and so $$r_i\tau\lambda\in\ r_iw\lambda-\mathbb{N}\Delta_{im}\subset
W\lambda-\mathbb{N}\Delta_{im}.$$ Hence the assertion.
\end{proof}

\subsubsection{}
\begin{lemma}\label{stab} The stabilizer of $\lambda\in\ P^+$ in $T$ is generated by the $r_i,\, i\in\ I$ which
stabilize $\lambda$, that is $\Stab_{T}(\lambda)=\langle r_i\,|\, \alpha_i^{\vee}(\lambda)=0\rangle$.
\end{lemma}
\begin{proof} Set $S := \langle r_i\,|\, \alpha_i(\lambda)=0\rangle$. Clearly, $S\subset \Stab_T(\lambda)$. Let $\tau\in\ \Stab_T(\lambda)$; we will show that $\tau\in\ S$. We argue by induction on $\ell(\tau)$.
If $\tau=r_i$, for $i\in\ I$, the assertion is clear. Let $\tau\in\ \Stab_{T}(\lambda)$ be such that
$\ell(\tau)>1$ and write $\tau=r_i\tau'$, with $\ell(\tau')<\ell(\tau)$. Then, by the previous lemma
$r_i\tau'\lambda=r_i(w\lambda-\beta)=r_i(\lambda-\gamma)$, with $\beta\in\ \mathbb{N}\Delta_{im}$ and
$\gamma\in\ \mathbb{N}\Delta^+$.

If $i\in\ I^{im}$, $\alpha_i^{\vee}(\tau'\lambda)\geq 0$, which forces $\tau'\lambda=\lambda$ and
$\alpha_i^{\vee}(\tau'\lambda)=0$. In particular, $\alpha_i^{\vee}(\lambda)=0$ and $\tau'\in\
\Stab_T(\lambda)$. Then $\tau'\in\ S$, by the induction hypothesis and $r_i\in\ S$, hence $\tau\in\ S$.

If $i\in\ I^{re}$, $\lambda = r_i\tau'\lambda=r_iw\lambda-r_i\beta$, hence $\beta=0$ and $r_iw\in\
W_{\lambda}\subset S$. Then $\tau'\lambda =w\lambda=\lambda$, so by the induction hypothesis $\tau'\in\ S$
and since $r_i\in\ S$, we get $\tau\in\ S$. Hence the assertion.
\end{proof}

\subsubsection{}
Let $\lambda\in\ P^+$ and recall section \ref{monoid}. One would like to know which elements in $T\lambda$
are dominant. Here we remark that by lemma \ref{Tlambda} one has that $T\lambda\cap P^+=T\lambda\cap
\mathcal{C}$. By section \ref{Wgeometry}, for all $w\in\ W$, with $w\notin W_{\lambda}$, $w\lambda$ is not
dominant. On the other hand, notice that for all dominant $\mu\in\ T\lambda$ and all $i\in\ I^{im}$, $r_i\mu$
is also dominant. Indeed, for all $j\in\ I$ we have that $\alpha_j^{\vee}(r_i\mu)
=\alpha_j^{\vee}(\mu)-\alpha_i^{\vee}(\mu)a_{ji}\geq 0$, since $\mu$ is dominant and $a_{ji}\leq 0$. In
particular, $r_{i_1}r_{i_2}\cdots r_{i_k}\lambda$ is dominant for all
$i_1, i_2, \dots, i_k\in\ I^{im}$.\\

\begin{lemma}\label{dominance} Let $\mu\in\ T\lambda\cap P^+$ and $i\in\ I^{im}$ and assume that $r_iw\mu\notin P^+$ for some $w\neq \id$ in $W$. Then $r_iw\mu = r_jr_iw'\mu$, for some $j\in\
I^{re}$ with $w':=r_jw$ and $\ell (w')=\ell (w)-1$. Consequently there exist $w_1,\, w_2\in\ W$ such that
$w=w_1w_2$ and $\ell(w)=\ell(w_1)+\ell(w_2)$, with $r_iw_1=w_1r_i$. Moreover, $\mu':=r_iw_2\mu$ is dominant
and $r_iw\mu = w_1\mu'$.
\end{lemma}
\begin{proof} By assumption, $r_iw\mu$ is not
dominant, hence there exists a $j\in\ I^{re}$ such that $\alpha_j^{\vee}(r_iw\mu)<0$. This gives
$$\alpha_j^{\vee}(w\mu)-\alpha^{\vee}_i(w\mu)a_{ji}<0,$$
hence
\begin{equation}\label{dominant}\alpha_j^{\vee}(w\mu)<\alpha^{\vee}_i(w\mu)a_{ji}.
\end{equation}
Now since $\mu\in T\lambda$, one has that $r_jw\mu\in\ T\lambda$ and so, by lemma \ref{Tlambda}, $\alpha_i^{\vee}(r_jw\mu)\geq 0$ which in turn gives
:
\begin{equation}\label{dominant2}
\alpha_i^{\vee}(w\mu)-\alpha^{\vee}_j(w\mu)a_{ij}\geq 0,
\end{equation}
Suppose that $a_{ji}$ (and so $a_{ij}$) is not equal to zero and hence $a_{ij},\, a_{ji}<0$. Then equations
(\ref{dominant}) and (\ref{dominant2}) give $$\alpha_i^{\vee}(w\mu)(1-a_{ij}a_{ji})>0.$$ But this is
impossible since $1-a_{ij}a_{ji}\leq 0$ and again by lemma \ref{Tlambda}, $\alpha_i^{\vee}(w\mu)\geq 0$. We
conclude that $a_{ij}=a_{ji}=0$, which implies that $r_i$ and $r_j$ commute and
$\alpha_j^{\vee}(r_iw\mu)=\alpha_j^{\vee}(w\mu)<0$. Since $\mu\in\ P^+$, the last inequality forces
$w=r_jw'$, for some $w'\in\ W$ with $\ell(w')= \ell(w)-1$. Finally, $r_iw = r_ir_jw'=r_jr_iw'$. By repeating
the procedure for $r_iw'\mu$, the last assertion follows.
\end{proof}

\subsubsection{}
\begin{lemma}\label{length} Let $\mu\in\ P^+$ and $\tau\in\ T$. If
$\alpha_j^{\vee}(\tau\mu)<0$, for some $j\in\ I^{re}$, then $\ell(r_j\tau)<\ell(\tau)$.
\end{lemma}
\begin{proof}
Let
\begin{equation}\label{reduced}
\tau = w_0r_{i_1}w_1\cdots w_{k-1}r_{i_k}w_k,
\end{equation}
with $w_t\in\ W,\, 0\leq t\leq k$ and $i_s\in\ I^{im},\, 1\leq s\leq k$ be a reduced expression of $\tau$. By
the previous lemma, we can write $r_{i_k}w_k\mu$ as $w'_k\mu'$, with $\mu' = r_{i_k}w''_k\mu\in\ P^+$,
$r_{i_k}w'_k=w'_kr_{i_k}$ and
\begin{equation}\label{lengthsadd}\ell(w_k)=\ell(w'_k)+\ell(w''_k).
\end{equation}
Thus we get a new expression for $\tau$ :
$$\tau = w_0r_{i_1}w_1\cdots
r_{i_{k-1}}w_{k-1}'r_{i_k}w_k'',$$ where $w_{k-1}'=w_{k-1}w_k'$. By (\ref{lengthsadd}) and since the
expression (\ref{reduced}) of $\tau$ is reduced we get $\ell(w_{k-1}')=\ell(w_{k-1})+\ell(w_k')$. Repeating
this procedure, we obtain $\tau\mu = w_0'\nu$ and $\nu = \tau'\mu\in\ P^+$, with
$\ell(\tau)=\ell(w_0')+\ell(\tau')$. Let $j\in\ I^{re}$. Then $\alpha_j^{\vee}(\tau\mu)<0$ implies that
$\alpha_j^{\vee}(w_0'\nu)<0$ and so $\ell(r_jw_0')<\ell(w_0')$ which in turn gives that
$\ell(r_j\tau)=\ell(r_jw_0'\tau')<\ell(w_0'\tau')=\ell(\tau)$.
\end{proof}

\subsubsection{}
For any $\mu\in\ T\lambda,\, \lambda\in\ P^+$, call $\tau\lambda$ a minimal representative of $\mu$ if $\mu =
\tau\lambda$ and for every $\tau'$ such that $\mu = \tau'\lambda$ one has that $\ell(\tau)\leq
\ell(\tau')$. 
We have the following result :\\

\begin{corollary}\label{formofdom} An element $\mu\neq \lambda$ in $T\lambda$ is dominant if and only if every minimal representative of $\mu$ is of the form $r_i\tau\lambda$, for $\tau\in\ T,\, i\in\
I^{im}$.
\end{corollary}
\begin{proof} Suppose that $\mu=\tau \lambda$ is a minimal representative of $\mu$. As in the proof of lemma \ref{length} one can write $\mu=w_0\nu$, with $\nu = r_i\tau'\lambda\in\ P^+,\, i\in\ I^{im}, \tau'\in\ T$ and $\tau=w_0r_i\tau'$, where lengths add. If $\mu$
is dominant, then $w_0\in\ W_{\nu}$. But then $\mu=w_0r_i\tau'\lambda=r_i\tau'\lambda$ which implies that
$w_0=\id$, hence every minimal representative of $\mu$ starts with some $r_i$ with $i\in\ I^{im}$. If $\mu$ is
not dominant, then by lemma \ref{dominance}, there exists a minimal representative of $\mu$ starting with $r_j$,
where $j\in\ I^{re}$. Hence the assertion.
\end{proof}

\subsubsection{}\label{domred}
We may express this consequence of lemma \ref{length} in the following fashion. Write
$\tau\in\ T$ as in equation (\ref{reduced}). Call $\tau$ a dominant reduced expression if $\tau$ is
reduced and successively the $\ell(w_k),\, \ell(w_{k-1}),\dots, \ell(w_0)$ take their minimal values. Set
$\tau'=r_{i_1}w_1\cdots r_{i_k}w_k$. Then for all $\mu\in P^+,\, \tau'\mu$ is dominant and $\tau\mu$ is dominant if and only if $w_0\in\ \Stab_W(\tau'\mu)$. In particular, if $w_0=\id$ in a dominant reduced expression of $\tau$, then $\tau\mu\in P^+$ for all $\mu\in P^+$.

\subsubsection{}
\begin{lemma}\label{twodominant} Let $\lambda,\, \mu\in\ P^+$. Then, for all $\tau\in\ T$ one has that $\beta^{\vee}(\tau\mu)>0$ implies $\beta^{\vee}(\tau\lambda)\geq 0$, for all $\beta\in\ W\Pi\cap
\Delta^+$.
\end{lemma}
\begin{proof} Since $\beta=w\alpha_i\in\ W\Pi$ one has $\beta^{\vee}(\tau\mu)=\alpha_i^{\vee}(w^{-1}\tau\mu)$ which reduces us to the case $\beta=\alpha_i\in\ \Pi$. For $\alpha_i\in\ \Pi_{im}$ one always has that $\alpha_i^{\vee}(\tau\lambda)\geq 0$, by lemma \ref{Tlambda}. Suppose that $\alpha_i\in\ \Pi_{re}$. Take $\tau\in\ T$ and let $\tau=w_0r_{i_1}w_1r_{i_2}\cdots r_{i_k}w_k$ be a dominant reduced expression of $\tau$; then $r_{i_1}w_1r_{i_2}\cdots r_{i_k}w_k$ is a dominant reduced expression for $\tau'$, where $\tau=w_0\tau'$. By Section \ref{domred}, one has that $\tau'\lambda, \,\tau'\mu\in\
P^+$.

Suppose that $\alpha_i^{\vee}(w_0\tau'\mu)>0$, then $(w_0^{-1}\alpha_i)^{\vee}(\tau'\mu)>0$, which implies
that $w_0^{-1}\alpha_i\in\ \Delta^+$, since $\tau'\mu\in\ P^+$. Then
$(w_0^{-1}\alpha_i)^{\vee}(\tau'\lambda)\geq 0$, since $\tau'\lambda\in\ P^+$ and so
$\alpha_i^{\vee}(\tau\lambda)\geq 0$.
\end{proof}

\noindent {\bf Remark.} We prove in Lemma \ref{twodom} a general fact about ``Bruhat sequences'' in $T$ which generalizes a well-known result for $W$.  For technical reasons this is postponed for the moment.

\section{Generalized crystals}

\subsection{The notion of a crystal}

\subsubsection{}\label{crystal}
{\bf Definition.} A generalized crystal $B$ is a set endowed with the maps $\wt: B\rightarrow P$,
$\varepsilon_i, \varphi_i: B\rightarrow \mathbb{Z}\cup \{-\infty\}$, $e_i, f_i:B\rightarrow B\cup\{0\}$
satisfying the rules :
\begin{enumerate}
\item For all $i\in\ I$ and all $b\in\ B$, $\varphi_i(b)=\varepsilon_i(b)+ \alpha_i^{\vee}(\wt b)$.
\item For all $i\in\ I$ if $b, e_ib\in\ B$, then $\wt (e_ib)= \wt b +\alpha_i$.
\item For all $i\in\ I$ if $b, e_ib\in\ B$, then $\varepsilon_i(e_ib)=\varepsilon_i(b)-1$ if $i\in\ I^{re}$
and $\varepsilon_i(e_ib) = \varepsilon_i(b)$ if $i\in\ I^{im}$.
\item For all $i\in\ I$ and all $b, b'\in\ B$
one has $b'=e_ib$ if and only if $f_ib'=b$.
\item If for $b\in\ B, i\in\ I, \varphi_i(b)=-\infty$, then $e_ib=f_ib=0$.
\end{enumerate}

\subsubsection{}\label{simple}
{\bf Remarks.}
\begin{enumerate}
\item The axioms imply the following further properties. First $\varphi_i(e_ib)=\varphi_i(b)+1$, if $i\in\ I^{re}$ and $\varphi_i(e_ib) = \varphi_i(b)+a_{ii}$, if $i\in\
I^{im}$. Second (a) $\wt f_ib=\wt b-\alpha_i$, (b) $\varepsilon_i(f_ib)=\varepsilon_i(b)+1$, if $i\in\
I^{re}$ and $\varepsilon_i(f_ib) = \varepsilon_i(b)$, if $i\in\ I^{im}$, (c) $\varphi_i(f_ib)=\varphi_i(b)-1$
if $i\in\ I^{re}$ and $\varphi_i(f_ib)=\varphi_i(b)-a_{ii}$ if $i\in\ I^{im}$.
\item The crystal graph of a crystal $B$ is the colored graph having vertices the elements of $B$ and arrows $b\stackrel{i}{\rightarrow}b'$ if $f_ib=b'$.
\item This definition is due to Jeong, Kang, Kashiwara and Shin \cite{JKKS}. We omit the term
``generalized'' in the sequel.
\end{enumerate}

\subsubsection{}\label{specialcrystals}
For any $\mu\in\ P$ set $B_{\mu} = \{b\in\ B\,|\, \wt b=\mu\}$. If $b\in B_\mu$, we say that $b$ is of weight $\mu$. If all $B_{\mu}$ are finite, define the
formal character of $B$ to be
$$\operatorname{char}\,B := \sum\limits_{b\in\ B}e^{\wt b} = \sum\limits_{\mu\in\
P}|B_{\mu}\,|\,e^{\mu}$$

Call a crystal $B$ {\it upper normal} if $\varepsilon_i(b)=\max\, \{n\in\ \mathbb{N}\,|\, e_i^nb\neq 0\}$ for
all $i\in\ I^{re}$, {\it lower normal} if $\varphi_i(b)=\max\, \{n\in\ \mathbb{N}\,|\, f_i^nb\neq 0\}$ for
all $i\in\ I^{re}$ and {\it normal} if it is both upper and lower normal.

Denote by $\mathcal{F}$ the monoid generated by the $f_i; i\in\ I$. A crystal $B$ is called a {\it highest
weight crystal} of highest weight $\lambda$ if there exists an element $b_{\lambda}\in\ B_\lambda$, such that $B=\mathcal{F}b_\lambda$. Notice that this implies that $e_ib_\lambda=0$ for all $i\in\ I$, but the
converse can fail. Despite the obvious analogy to highest weight modules, this condition is rather weak (see
also remark in section \ref{crystaltensorproduct}). Indeed, given a crystal $B$ and an element $b_\lambda\in\
B_{\lambda}$, we obtain a highest weight subcrystal $\mathcal{F}b_\lambda$ of $B$, simply by declaring $e_ib'=0$,
whenever $e_ib'\notin \mathcal{F}b_\lambda$.

\subsubsection{}\label{categoryB}
Let $\mathcal{B}$ be the set of crystals $B$ which for all $b\in\ B$ and all $i\in\ I^{im}$ satisfy :
\begin{enumerate}
\item $\alpha_i^{\vee}(\wt b)\geq 0$,
\item $\varepsilon_i(b) = 0$ and consequently $\varphi_i(b)=\alpha_i^{\vee}(\wt b)$,
\item $f_ib\neq 0$ if and only if $\varphi_i(b)>0$.
\end{enumerate}

\subsubsection{}
\begin{lemma} \label{conditionfore} Let $B\in\ \mathcal{B}$ and take $i\in\ I^{im},\,b\in\ B$. If
$\alpha_i^{\vee}(\wt b)\leq -a_{ii}$, then $e_ib=0$. In particular, $e_ib=0$ if $\alpha_i^{\vee}(\wt b)=0$.
\end{lemma}
\begin{proof} Suppose that $e_ib\neq 0$, then $f_i(e_ib)\neq 0$ by \ref{crystal} (4), and so $0<\varphi_i(e_ib) = \alpha_i^{\vee}(\wt e_ib)$. By \ref{crystal} (2), $\wt e_ib=\wt b+\alpha_i$ and
so $\alpha_i^{\vee}(\wt e_ib)> 0$ implies that $\alpha_i^{\vee}(\wt b)>-a_{ii}$.
\end{proof}
\noindent {\bf Remark.} The converse of the above lemma is false.

\subsubsection{}\label{morphism}
{\bf Definition.} A morphism $\psi$ of crystals $B_1, B_2$ is a map
$$\psi:B_1\longrightarrow B_2\cup\{0\}$$ such that:
\begin{enumerate}
\item $\wt (\psi(b))=\wt b$, $\varepsilon_i(\psi(b))= \varepsilon_i(b)$, $\varphi_i(\psi(b))=\varphi_i(b)$ for
all $i\in\ I$.
\item $\psi(f_ib)=f_i(\psi(b))$, if $f_ib\neq 0$.
\end{enumerate}

Notice that if $\psi$ is a crystal morphism, then also $\psi(e_ib)=e_i\psi(b)$ if $e_ib\neq 0$. One says that $B_1$ is a subcrystal of $B_2$ if $\psi$ is an embedding. An embedding is said to be {\it
strict}, if $f_i$ and $e_i$ commute with $\psi$ for all $i\in\ I$. If $\psi$ is a strict embedding, then $B_1$ is said
to be a {\it strict subcrystal} of $B_2$. The crystal graph of a subcrystal $B_1$ of $B_2$ is obtained by
removing the arrows between vertices of $B_1$ and vertices of $B_2\setminus B_1$ in the crystal graph of
$B_2$.

\subsection{Crystal tensor product}

\subsubsection{}\label{crystaltensorproduct}
{\bf Definition.} Let $B_1, B_2$ be two crystals. Their tensor product $B_1\otimes B_2$ is $B_1\times B_2$ as
a set, with crystal operations defined as follows. Set $b = b_1\otimes b_2$ with $b_1\in\ B_1,\, b_2\in\
B_2$. Then :\\
(1) $\wt b=\wt b_1 + \wt b_2$.\\
(2) $\varepsilon_i(b)=\max \{\varepsilon_i(b_1),\, \varepsilon_i(b_2)-\alpha_i^{\vee}(\wt b_1)\}$.\\
(3) $\varphi_i(b) = \max\{\varphi_i(b_1)+\alpha_i^{\vee}(\wt b_2),\, \varphi_i(b_2)\}$.\\
(4) For all $i\in\ I$,
\begin{equation*}\label{crystalfi}
f_ib =\left\{
\begin{array}{lc} f_ib_1\otimes b_2, \qquad \mbox{if}\qquad \varphi_i(b_1)>\varepsilon_i(b_2),\\
b_1\otimes f_ib_2, \qquad \mbox{if}\qquad \varphi_i(b_1)\leq \varepsilon_i(b_2).
\end{array}
\right.\tag{a}
\end{equation*}
(5) For all $i\in\ I^{re}$,
\begin{equation*}\label{crystalere} e_ib =\left\{
\begin{array}{lc} e_ib_1\otimes b_2, \qquad \mbox{if}\qquad \varphi_i(b_1)\geq \varepsilon_i(b_2),\\
b_1\otimes e_ib_2, \qquad \mbox{if}\qquad \varphi_i(b_1)< \varepsilon_i(b_2),
\end{array}
\right.\tag{b}
\end{equation*}
and for all $i\in\ I^{im}$ we set
\begin{equation*}\label{crystaleim} e_ib =\left\{
\begin{array}{lc} e_ib_1\otimes b_2, \qquad \mbox{if}\qquad \varphi_i(b_1)> \varepsilon_i(b_2)-a_{ii},\\
0,\qquad \mbox{if}\qquad \varepsilon_i(b_2)<\varphi_i(b_1)\leq \varepsilon_i(b_2)-a_{ii},\\
b_1\otimes e_ib_2, \qquad \mbox{if}\qquad \varphi_i(b_1)\leq \varepsilon_i(b_2).
\end{array}
\right.\tag{c}
\end{equation*}

It is straightforward to verify that $B_1\otimes B_2$ endowed with the above operations is indeed a crystal
\cite[Lemma 2.10]{JKKS}. Moreover, as in the Kac-Moody case, the tensor product of two normal crystals is a
normal crystal.\\

\noindent {\bf Remark.} If $\mathcal{F}b_{\lambda}$ and $\mathcal{F}b_{\mu}$ are highest weight crystals, it
is not obvious that $\mathcal{F}(b_{\lambda}\otimes b_{\mu})$ is a strict subcrystal of
$\mathcal{F}b_{\lambda}\otimes \mathcal{F}b_{\mu}$.

\subsubsection{}\label{tensorinB}
Let $B_1, \,B_2$ be crystals in $\mathcal{B}$, form their tensor product $B :=B_1\otimes B_2$ and let
$b:=b_1\otimes b_2\in\ B$. Take $i\in\ I^{im}$. The formulae \ref{crystaltensorproduct} (\ref{crystalfi}) and
(\ref{crystaleim}) simplify as follows :
\begin{equation*}\label{finB} f_ib
=\left\{
\begin{array}{lc} f_ib_1\otimes b_2, \qquad \mbox{if}\qquad \varphi_i(b_1)>0,\\
b_1\otimes f_ib_2, \qquad \mbox{if}\qquad \varphi_i(b_1)=0,
\end{array}
\right.\tag{a'}
\end{equation*}
and
\begin{equation*}\label{einB} e_ib =\left\{
\begin{array}{lc} e_ib_1\otimes b_2, \qquad \mbox{if}\qquad \varphi_i(b_1)> 0,\\
b_1\otimes e_ib_2, \qquad \mbox{if}\qquad \varphi_i(b_1)=0.
\end{array}
\right.\tag{c'}
\end{equation*}
Indeed, equation (\ref{finB}) above immediately obtains from \ref{crystaltensorproduct} (\ref{crystalfi})
since $\varphi_i(b_1)\geq 0=\varepsilon_i(b_2)$. For $e_ib$ notice that the only case where equation
\ref{crystaltensorproduct} (\ref{crystaleim}) and equation (\ref{einB}) above can differ is when
$0<\varphi_i(b_1)\leq -a_{ii}$. But then by lemma \ref{conditionfore} one has that $e_ib_1=0$ and so
$e_i(b_1\otimes b_2)=0$ by either (c) or (c').

We show that $B\in \mathcal{B}$ and thus the set $\mathcal{B}$ is closed under tensor products. Indeed notice that $\alpha_i^{\vee}(\wt
b)=\alpha_i^{\vee}(\wt b_1)+\alpha_i^{\vee}(\wt b_2)\geq 0$ and $\varepsilon_i(b) = \max
\{\varepsilon_i(b_1),\, \varepsilon_i(b_2)-\alpha_i^{\vee}(\wt b_1)\}=0$. Now if $\varphi_i(b)>0$, then
either $\varphi_i(b_1)>0$ and $f_ib=(f_ib_1)\otimes b_2\neq 0$ or $\varphi_i(b_1)=0,\, \varphi_i(b_2)>0$ and
$f_ib=b_1\otimes f_ib_2\neq 0$. On the other hand, if $\varphi_i(b)=0$, then
$\varphi_i(b_1)=\varphi_i(b_2)=0$ which implies that $f_ib_1=f_ib_2=0$ and so $f_ib=0$. We conclude that
$f_ib\neq 0$ if and only if $\varphi_i(b)>0$ as required.

\subsection{The crystal $B(\infty)$}

\subsubsection{}\label{elementary}
For any index $i\in\ I$ we define the {\it elementary crystal} $B_i$ \cite[Example 2.14]{JKKS} to be the set $B_i = \{b_i(-n)\,|\,n\in\ \mathbb{N}\}$ with crystal operations:\\
$$
\begin{array}{l}\wt b_i(-n)=-n\alpha_i,\\
e_ib_i(-n) = b_i(-n+1), \quad f_ib_i(-n)=b_i(-n-1)\\
e_jb_i(-n)=f_jb_i(-n)=0, \quad \mbox{if }\, i\neq j\\
\varepsilon_i(b_i(-n)) = n, \quad \varphi_i(b_i(-n))= -n, \quad \mbox{if }\,  i\in\ I^{re}\\
\varepsilon_i(b_i(-n)) = 0, \quad \varphi_i(b_i(-n))= -na_{ii}, \quad \mbox{if }\, i\in\ I^{im}\\
\varepsilon_j(b_i(-n))=\varphi_j(b_i(-n))=-\infty, \quad \mbox{if } \, i\neq j,
\end{array}
$$
where we have set $b_i(-n)=0$ for all $n<0$.

\subsubsection{}\label{Binfty}
\begin{theorem}  There exists a unique (up to isomorphism) crystal, denoted by $B(\infty)$, with the properties :
\begin{enumerate}
\item There exists an element $b_0\in\ B(\infty)$ of weight zero.
\item The set of weights of $B(\infty)$ lies in $- Q^+$.
\item For any element $b\in\ B(\infty)$ with $b\neq b_0$, there exists some $i\in\ I$ such that $e_ib\neq 0$.
\item For all $i\in\ I$ there exists a unique strict embedding $\Psi_i : B(\infty)\longrightarrow B(\infty)\otimes B_i$, sending $b_0$ to $b_0\otimes b_i(0)$.
\end{enumerate}
\end{theorem}
The description of $B(\infty)$ which results is given in section \ref{infty} below.

\subsubsection{} \label{BJcrystal}
Let $J=\{i_1, i_2, \dots\}$ where $i_j\in\ I$ is a countable sequence with the property that for all $i\in\
I$ and all $j\in\ \mathbb{N}^+$, there exists $k>j$ such that $i_k=i$. It is convenient to assume that
$i_j\neq i_{j+1}$ for all $j\in\ \mathbb{N}^+$. Set $B(k)=B_{i_k}\otimes \cdots \otimes B_{i_1}$ and for
$k\leq l$, let $\psi_{k, l}:B(k)\rightarrow B(l)$ be the map $b\mapsto b_{i_l}(0)\otimes\cdots\otimes
b_{i_k+1}(0)\otimes b$. Let $B_J(\infty)$ be the inductive limit of the family $\{B(k)\}_{k\geq 1}$. Then
$B_J(\infty)$ is the crystal in which an element $b$ takes the form
\begin{equation*} b=\dots\otimes b_{i_2}(-m_2)\otimes b_{i_1}(-m_1),
\end{equation*}
with $m_k\in\ \mathbb{N}$ and $m_k=0$ for $k\gg0$.

\subsubsection{}\label{infty}
Let $B$ be a crystal satisfying properties (1)-(4) of theorem \ref{Binfty}. Then $b_0$ is the unique element
of weight zero in $B$. Indeed, if $b\neq b_0$ and $\wt b=0$ then $e_ib\neq 0$ for some $i\in\ I$. But then
$\wt e_ib=\alpha_i\notin -Q^+$ contradicting property (2). It follows that $B=\mathcal{F}b_0$.

Iterating (4) we have a strict embedding :
$$B\hookrightarrow B\otimes B_{i_1}\hookrightarrow B\otimes B_{i_2}\otimes B_{i_1}\hookrightarrow \cdots  \hookrightarrow B\otimes B_{i_r}\otimes \cdots \otimes B_{i_2}\otimes B_{i_1},$$
for all $r>0$. There exists $N>0$ such that any element $b\in\ B$ takes the form
$$b_0\otimes b_{i_N}(-m_N)\otimes \cdots \otimes b_{i_1}(-m_1).$$ Associating $\cdots \otimes
b_{i_{N+1}}(0)\otimes b_{i_N}(-m_N)\otimes \cdots \otimes b_{i_1}(-m_1)$ to $b$ we obtain a strict embedding
$B\hookrightarrow B_J(\infty)$. Now $B_J(\infty)$ admits a unique element $b_{\infty}$ of weight zero given
by taking all the $m_k=0$ for all $k\in\ \mathbb{N}^+$. Then $B$ is the strict subcrystal of $B_J(\infty)$
generated by $b_{\infty}$. We conclude that a crystal satisfying (1)-(4) of theorem \ref{Binfty} is unique.

\subsubsection{} {\bf Remark.} For $A$ symmetrizable and $a_{ii}\in\ -2\mathbb{N}^+$ if $i\in\ I^{im}$, theorem \ref{Binfty} is due to Jeong,
Kang and Kashiwara \cite[Theorem 4.1]{JKKS}. Their proof is not combinatorial. We shall prove it
combinatorially and in general by constructing a path model. Moreover, the crystal structure of $B_J(\infty)$ is given explicitly
in \cite[Example 2.17]{JKKS}. We describe it in section \ref{Kashiwara}, where some further properties of
$B_J(\infty)$ are discussed .

\section{A path model for crystals defined by a Borcherds-Cartan matrix}\label{Generalized Lakshmibai-Seshadri}
According to our general conventions, all intervals are considered in $\mathbb{Q}$, that is we write $[a, b]$
for $\{c\in\ \mathbb{Q}\,|\,a\leq c\leq b \}$. Let $X$ be a topological space. A function $\pi :[0,
1]\rightarrow X$ is said to be continuous (or just a path) if it is the restriction of a continuous function
on the real interval. Actually, we shall mainly use piecewise linear functions.

Let $\mathbb{P}$ be the set of paths $\pi : [0,\, 1]\rightarrow \mathbb{Q}P$ such that $\pi(0) = 0$ and
$\pi(1)\in\ P$. We consider two paths $\pi, \,\pi'\in\ \mathbb{P}$ equivalent if $\pi = \pi'$ up to
parametrization, i.e. if there exists a non-decreasing continuous function $\phi : [0, 1]\rightarrow [0, 1]$
such that $\pi(\phi(t)) = \pi'(t)$ for all $t\in\ [0, 1]$.

\subsection{The operators $f_i,\, e_i$}

\subsubsection{}\label{definitionofh}
For all $\pi\in\ \mathbb{P}$ and all $i\in\ I$, set $h_i^{\pi}(t) := \alpha_i^{\vee}(\pi(t)),\, t\in\ [0, 1]$
and let $m_i^{\pi}$ be the minimal integral value of the function $h_i^{\pi}$, that is
$$m_i^{\pi} = \min \{h_i^{\pi}(t)\cap \mathbb{Z}| t\in\ [0, 1]\}.$$
(Notice that since $\pi (0)= 0$, one has that $h_i^{\pi}(0)=0$, hence the function $h_i^{\pi}$ attains
integral values.) The action of $f_i,\, e_i$ for $i\in\ I$ is defined in the following sections.

\subsubsection{}\label{definitionoff}
Let $f_+^i(\pi)\in\ [0, 1]$ be maximal such that $h_i^{\pi}(f^i_+(\pi)) = m^{\pi}_i$. If $f_+^i(\pi)= 1$, set $f_i\pi = 0$. Otherwise, since $\pi (1)\in\ P$ and so $h_i^{\pi}(1)\in\ \mathbb{Z}$, it follows that the function $h_i^\pi$ attains the value $m_i^\pi$ in the interval $[f_+^i(\pi), 1]$. Let $f_-^i(\pi)\in\
[f_+^i(\pi), 1]$ be minimal such that $h_i^{\pi}(f^i_-(\pi)) = m_i^{\pi} + 1$. Then define $(f_i\pi)(t)$ to
be the path
:\\
$$(f_i\pi) (t) = \left\{
\begin{array}{ll} \pi(t), &t\in\ [0, f_+^i(\pi)],\\
\pi(f_+^i(\pi)) + r_i (\pi(t)-\pi(f_+^i(\pi))), &t\in\ [f_+^i(\pi), f_-^i(\pi)],\\
\pi(t)-\alpha_i, &t\in\ [f_-^i(\pi), 1].
\end{array} \right.$$

\subsubsection{}\label{definitionofere}
Take $i\in\ I^{re}$, and let $e_+^i(\pi)\in\ [0, 1]$ be minimal such that $h^{\pi}_i(e_+^i(\pi)) = m^{\pi}_i$. If $e_+^i(\pi)=0$, set $e_i\pi = 0$. If $e_+^i(\pi)>0$ and since $h_i^\pi(0)=0$, the function $h_i^\pi$ attains the value $m_i^\pi+1$ in the interval $[0, e_+^i(\pi)]$. Let $e_-^i(\pi)\in\ [0, e_+^i(\pi)]$ be maximal such that $h_i^{\pi}(e_-^i(\pi)) = m^{\pi}_i+1$. The path $e_i\pi$ is then defined by :\\
$$(e_i\pi) (t) = \left\{
\begin{array}{ll} \pi(t), &t\in\ [0, e_-^i(\pi)],\\
\pi(e_-^i(\pi)) + r_i (\pi(t)-\pi(e_-^i(\pi))), &t\in\ [e_-^i(\pi), e_+^i(\pi)],\\
\pi(t) + \alpha_i, &t\in\ [e_+^i(\pi), 1].
\end{array} \right.$$

\subsubsection{}\label{definitionofeim}
Take $i\in\ I^{im}$. Define $r_i^{-1}:\lieh^*\rightarrow \lieh^*$ to be the map :
$$r_i^{-1}(x) = x+\frac{1}{1-a_{ii}}\alpha_i^{\vee}(x)\alpha_i.$$
One checks that $r_ir_i^{-1}=r_i^{-1}r_i=\id$. Recall the number $f_+^i(\pi)$ defined in section
\ref{definitionoff} and set $e_-^i(\pi):=f_+^i(\pi)$. If $e_-^i(\pi)=1$ or $h_i^{\pi}(t)<m_i^{\pi}+1-a_{ii}$
for all $t\in\ [e_-^i(\pi), 1]$, set $e_i\pi=0$. Otherwise let $e_+^i(\pi)\in\ [e_-^i(\pi), 1]$ be minimal such that
$h_i^{\pi}(e_+^i)=m_i^{\pi}+1-a_{ii}$. If $h_i^{\pi}(t)\leq m_i^{\pi}-a_{ii}$ for some $t\in\ [e_+^i(\pi), 1]$, set
$e_i\pi= 0$. Otherwise, define $e_i\pi$ to be the path :\\
$$(e_i\pi)(t) = \left\{
\begin{array}{ll} \pi(t), &t\in\ [0, e_-^i(\pi)],\\
\pi(e_-^i(\pi))+ r_i^{-1} (\pi(t)-\pi(e_-^i(\pi))), &t\in\ [e_-^i(\pi), e_+^i(\pi)],\\
\pi(t) + \alpha_i, &t\in\ [e_+^i(\pi), 1].
\end{array} \right.$$

\subsubsection{}\label{rem} {\bf{Remarks.}}
\begin{enumerate}
\item The definition of $f_i,\, e_i$ for $i\in\ I^{re}$ is as in \cite[Section 2]{L1}. Notice that in \cite{L1} the condition under which Littelmann sets $f_i\pi = 0$  is that $h_i^{\pi}(1)-h_i^{\pi}(f_+^i(\pi))< 1$. This is equivalent to equality in $h_i^{\pi}(1)\geq m_i^{\pi}$ and so to $f^i_+(\pi)=1$ if we consider only paths with endpoint in
$P$.
\item Littelmann gave a different and more involved definition for the root operators in \cite{L2}. These operators were compatible with the ``stretching of paths'', an essential tool in the proof of his Isomorphism Theorem. However, for the paths appearing in this paper (namely, integral and monotone paths, see sections \ref{integral} and \ref{monotone}) the two definitions coincide. In order to prove an analogue of Littelmann's Isomorphism theorem for generalized Kac-Moody algebras, one would have to consider more general paths and hence use the definitions of \cite{L2}. We note here that the Isomorphism theorem is equivalent to the tensor product decomposition, proven in \cite{JKK} for crystal bases, if paths admit ``stretching''.
\item It is easy to verify that for $e_i, \, f_i,\, i\in\ I$ defined above, $f_i\pi=\pi'$ if and only if
$e_i\pi'=\pi$. For $i\in\ I^{re}$, this is done in \cite{L2}.
\item If $f_i\pi\neq 0$, one has that $(f_i\pi)(1) = \pi(1) -\alpha_i$. Similarly, if $e_i\pi\neq 0$ then $(e_i\pi)(1) = \pi(1) +\alpha_i$.
\end{enumerate}

\subsubsection{}
\begin{lemma}\label{computationofhfi} Take $i\in\ I$ and let $\pi\in\ \mathbb{P}$ be such that $f_i\pi\neq 0$.
\begin{enumerate}
\item If $i\in\ I^{im}$, then $m_i^{f_i\pi} = m_i^{\pi}$ and $f^i_+(f_i\pi) = f_+^i(\pi)$, whereas $f_-^i(f_i\pi)\leq f_-^i(\pi)$, with equality if and only if $a_{ii} =
0$. In particular, $f_i^k\pi\neq 0$, for all $k\geq 0$.
\item If $i\in\ I^{re}$, then $m_i^{f_i\pi} = m_i^{\pi}-1$ and $f^i_+(f_i\pi) = f_-^i(\pi)$. In particular, since $h_i^{f_i\pi}(1)=h_i^{\pi}(1)-2$, there exists $k\in\ \mathbb{N}$ such that $f_i^k\pi=0$.
\end{enumerate}
\end{lemma}
\begin{proof} Consider (1) and let $i\in\ I^{im}$. By definition, $h_i^{\pi}(t)\cap \mathbb{Z}\geq h_i^{\pi}(f_+^i(\pi))=m_i^{\pi}$ and this inequality is strict for $t>f_+^i(\pi)$. We will compute the function $h_i^{f_i\pi}(t)$. Recall definition \ref{definitionoff}. One has that for $t\in\ [0, f_+^i(\pi)]$
\begin{equation}\label{computationhfi1}h_i^{f_i\pi}(t)=h_i^{\pi}(t).
\end{equation}
Now, for $t\in\ [f_+^i(\pi), f_-^i(\pi)]$,
\begin{equation}\label{computationhfi2}h_i^{f_i\pi}(t) = h_i^{\pi}(t)-a_{ii}(h_i^{\pi}(t)-h_i^{\pi}(f_+^i(\pi)))\geq h_i^{\pi}(t).
\end{equation}
Finally, for $t\in\ [f_-^i(\pi), 1]$,
\begin{equation}\label{computationhfi3}h_i^{f_i\pi}(t)=h_i^{\pi}(t)-a_{ii}\geq h_i^{\pi}(t).
\end{equation}
By equations (\ref{computationhfi1}), (\ref{computationhfi2}), (\ref{computationhfi3}), we conclude that
$h_i^{f_i\pi}(t)\geq h_i^{\pi}(t)$ for all $t\in\ [0, 1]$ and so $h_i^{f_i\pi}(t)\cap \mathbb{Z}\geq
m_i^{\pi}$. Since also $h_i^{f_i\pi}(f_+^i(\pi)) = h_i^{\pi}(f_+^i(\pi))= m_i^{\pi}$, we conclude that
$m_i^{f_i\pi}=m_i^{\pi}$. Also $h_i^{f_i\pi}(t)\cap \mathbb{Z}> m_i^{\pi}$ for $t>f_+^i(\pi)$ and thus
$f_+^i(f_i\pi)=f_+^i(\pi)$. By (\ref{computationhfi2}) we obtain that $f_-^i(f_i\pi)\leq f_-^i(\pi)$ with
equality if and only if $a_{ii}=0$. Finally, since $f_i\pi=0$ if and only if $f_+^i(\pi)=1$ and
$f_+^i(f_i\pi)=f_+^i(\pi)$, one has that $f_i\pi\neq 0$ implies that $f_i^2\pi\neq 0$ and inductively,
$f_i^k\pi\neq 0$ for all $k\geq 0$. Hence (1).

Statement (2) which we have included for comparison is implicit in \cite[Proposition 1.5]{L1}. It may be
similarly verified by substituting $a_{ii}=2$ in the first parts of equations (\ref{computationhfi2}) and
(\ref{computationhfi3}).
\end{proof}

\subsection{The Crystal structure of $\mathbb{P}$}

\subsubsection{}\label{crystalstructure}
We will endow $\mathbb{P}$ with a normal crystal structure. We define the operators $f_i, \,e_i, \,i\in\ I$
as in sections \ref{definitionoff}, \ref{definitionofere}, \ref{definitionofeim}. We set $\wt \pi =\pi(1)$.
For $i\in\ I^{re}$ we set $\varepsilon_i(\pi) =-m_i^{\pi}$. For $i\in\ I^{im}$, we set
$\varepsilon_i(\pi)=0$. Then $\varphi_i$ can be recovered by the formula $\varphi_i(\pi) =
\varepsilon_i(\pi)+\alpha_i^{\vee}(\wt \pi)$. From section \ref{rem} and lemma \ref{computationofhfi} one checks the following:\\

\begin{lemma}\label{pathcrystal} The set of paths $\mathbb{P}$ together with the maps $e_i, \,f_i,\, \varepsilon_i,\,\varphi_i, \,\wt$ for all $i\in\ I$ defined above, is a normal crystal.
\end{lemma}

\subsubsection{Concatenation of paths}
We define the tensor product of $\pi_1,\, \pi_2\in\ \mathbb{P}$ to be the concatenation of the two paths :\\
$$(\pi_1\otimes \pi_2)(t) =
\left\{
\begin{array}{lr} \pi_1(t/s), \qquad t\in\ [0, s],\\
\pi_1(1)+\pi_2(\frac{t-s}{1-s}), \qquad t\in\ [s, 1],
\end{array}
\right.$$ for any rational number $s\in\ [0, 1]$.\\

\begin{lemma}
The crystal operations on $\mathbb{P}\otimes \mathbb{P}\subset \mathbb{P}$ satisfy the crystal tensor product
rules defined in section \ref{crystaltensorproduct}.
\end{lemma}
\begin{proof} This is straightforward; a point to remark is that $(\pi_1\otimes \pi_2)(s)\in\ P$, otherwise
the insertion of $e_i$ or $f_i$ will simultaneously change both $\pi_1$ and $\pi_2$.
\end{proof}

\section{Generalized Lakshmibai-Seshadri paths}

\subsection{Distance of two weights in $T\lambda$}
Notation is as in sections \ref{BCmatrix}-\ref{delta}.

\subsubsection{}
Let $\lambda\in\ P^+$ and let $\mu,\, \nu\in\ T\lambda$ be two weights in the $T$ orbit of $\lambda$. We
write $\mu>\nu$ if there exists a sequence of weights $\mu := \lambda_0, \lambda_1, \dots, \lambda_{s-1},
\lambda_s:=\nu$ and positive roots $\beta_1, \dots, \beta_s\in\ W\Pi\cap \Delta^+ = \Delta_{re}^+\sqcup
\Delta_{im}$ such that $\lambda_{i-1}=r_{\beta_i}\lambda_i$ and $\beta_i^{\vee}(\lambda_i)>0$, for all $i$,
with $1\leq i\leq s$. Note that $\mu=r_{\beta}\nu$, with $\beta\in W\Pi\cap \Delta^+$, one has $\mu>\nu$ if
and only if $\beta^{\vee}(\nu)>0$.

We call the {\it distance} of $\mu$ and $\nu$ and write $\dist\, (\mu, \nu)$ the maximal length of such
sequences. If $\mu = r_{\beta}\nu>\nu$ and $\dist\,(\mu, \nu)=1$ we write $\nu \stackrel{\beta}{\leftarrow}
\mu$. Since $\beta$ is uniquely determined by the pair $(\mu, \nu)$, we can omit it and write $\nu\leftarrow
\mu$.

\subsubsection{}
{\bf Remarks.}\label{Bruhat}
\begin{enumerate}
\item If $\beta_i\in\ \Delta^+_{re}$ for all $i,\, 1\leq i\leq s$ then $\id\leftarrow r_{\beta_s}\leftarrow \cdots \leftarrow r_{\beta_1}r_{\beta_2}\cdots r_{\beta_s}$ is a Bruhat sequence in $W/W_{\lambda}$, where recall that $W_{\lambda}$ stands for the stabilizer of $\lambda$ in $W$.
\item Let $\mu =r_i\nu$ and $\alpha_i^{\vee}(\nu)>0$. Then $\dist\, (\mu, \nu)=1$. Indeed, note that
$$\nu:=\lambda_s\stackrel{\beta_s}{\leftarrow}\lambda_{s-1}\cdots\stackrel{\beta_2}{\leftarrow}\lambda_1\stackrel{\beta_1}{\leftarrow}\lambda_0=:\mu$$
with $\beta_i\in\ \Delta^+$ means that $\nu = \mu+\sum\limits_{t=1}^sn_t\beta_t$, with $n_t\in\
\mathbb{N}^+$. Note that $\mu>\nu$ implies that $\mu \prec \nu$. The converse fails.
\end{enumerate}

\subsubsection{} The following is exactly as in \cite[Lemma 4.1]{L2}.\\

\begin{lemma}\label{distances} Let $\alpha_i\in\ \Pi_{re}$ be a simple real root and let $\mu\geq \nu$ be two weights in $T\lambda$ with $\lambda\in\
P^+$. Then :
\begin{enumerate}
\item If $\alpha_i^{\vee}(\mu)<0$ and $\alpha_i^{\vee}(\nu)\geq 0$, then $r_i\mu\geq \nu$ and $\dist \,(r_i\mu, \nu)<\dist\,(\mu, \nu)$.
\item If $\alpha_i^{\vee}(\mu)\leq 0$ and $\alpha_i^{\vee}(\nu)> 0$, then $\mu\geq r_i\nu$ and $\dist \,(\mu, r_i\nu)<\dist\,(\mu, \nu)$.
\item If $\alpha_i^{\vee}(\mu)\alpha_i^{\vee}(\nu)> 0$, then $r_i\mu\geq r_i\nu$ and $\dist\,(r_i\mu, r_i\nu)= \dist \,(\mu, \nu)$.
\end{enumerate}
\end{lemma}

\subsubsection{}
\begin{lemma}\label{positivityofhim} Let $\mu\geq \nu\in\ T\lambda$ be such that $\nu : = \nu_s\stackrel{\beta_s}{\leftarrow} \nu_{s-1}\stackrel{\beta_{s-1}}{\leftarrow}\dots \stackrel{\beta_2}{\leftarrow}\nu_1\stackrel{\beta_1}{\leftarrow}\nu_0=:\mu$ where $\beta_j\in\ W\Pi\cap \Delta^+$ with $1\leq j\leq s$ and let $i\in\ I^{im}$. Then
$\alpha_i^{\vee}(\mu)\geq \alpha_i^{\vee}(\nu)$ and $\alpha_i^{\vee}(\mu) = \alpha_i^{\vee}(\nu)$ if and only
if $r_i$ commutes with $r_{\beta_j}$ for all $j$ with $1\leq j\leq s$.
\end{lemma}
\begin{proof} Since $\mu = \nu -\sum\limits_{j=1}^s\beta_j^{\vee}(\nu_j)\beta_j\in\ \nu-\sum\limits_{j=1}^s\mathbb{N}^+\beta_j$ and $\alpha_i^{\vee}(\beta_j)\leq
0$ for all $j$, with $1\leq j\leq s$, we conclude that $\alpha_i^{\vee}(\mu)\geq \alpha_i^{\vee}(\nu)$ and
equality holds if and only if $\alpha_i^{\vee}(\beta_j)= 0$ for all $j$, with $1\leq j\leq s$. The latter is
equivalent to $r_ir_{\beta_j}=r_{\beta_j}r_i$ for all $j$, with $1\leq j\leq s$.
\end{proof}

\subsubsection{}
\begin{lemma}\label{distanceim} Let $\alpha_i\in\ \Pi_{im}$ be a simple imaginary root and let $\mu\geq \nu$ be two weights in $T\lambda$ with $\lambda\in\
P^+$. If  $\alpha_i^{\vee}(\mu)=\alpha_i^{\vee}(\nu)\geq 0$, then $r_i\mu\geq r_i\nu$ and $\dist\,(r_i\mu,
r_i\nu)=\dist\, (\mu, \nu)$.
\end{lemma}
\begin{proof} Set $\dist\, (\mu, \nu)=s\geq 1$, then
$$\nu : = \nu_s\stackrel{\beta_s}{\leftarrow} \nu_{s-1}\stackrel{\beta_{s-1}}{\leftarrow}\dots \stackrel{\beta_2}{\leftarrow}\nu_1\stackrel{\beta_1}{\leftarrow}\nu_0=:\mu,$$
for some $\beta_j\in\ W\Pi\cap \Delta^+$, where $1\leq j\leq s$. Since by assumption
$\alpha_i^{\vee}(\mu)=\alpha_i^{\vee}(\nu)$, lemma \ref{positivityofhim} gives that $r_i$ commutes with
$r_{\beta_j}$ for all $j$, with $1\leq j\leq s$. Notice that this means that $r_i\mu =r_{\beta_1}\cdots
r_{\beta_s}r_i\nu$ and $\beta_j^{\vee}(r_{\beta_{j+1}}\cdots r_{\beta_s}r_i\nu) =
\beta_j^{\vee}(r_{\beta_{j+1}}\cdots r_{\beta_s}\nu)>0$. Hence $r_i\mu\geq r_i\nu$ and $\dist \,(r_i\mu,
r_i\nu)\geq s$. Suppose that $\dist\,(r_i\mu, r_i\nu)>s$. This means that there exist positive roots
$\gamma_j$, with $1\leq j\leq t$ and $t>s$, such that
$$r_i\nu : = \nu'_t\stackrel{\gamma_t}{\leftarrow}
\nu'_{t-1}\stackrel{\gamma_{t-1}}{\leftarrow}\dots
\stackrel{\gamma_2}{\leftarrow}\nu'_1\stackrel{\gamma_1}{\leftarrow}\nu'_0=: r_i\mu.
$$
But our hypothesis $\alpha_i^{\vee}(\mu)=\alpha_i^{\vee}(\nu)$ also implies that
$\alpha_i^{\vee}(r_i\mu)=\alpha_i^{\vee}(r_i\nu)$. By lemma \ref{positivityofhim}, $r_i$ commutes with
$r_{\gamma_j}$ for all $j$, with $1\leq j\leq t$. This gives us $r_i\mu=r_ir_{\gamma_1}\cdots
r_{\gamma_t}\nu$ and so $\mu=r_{\gamma_1}\cdots r_{\gamma_t}\nu$, therefore $\dist\,(\mu,\nu)\geq t>s$, which
is a contradiction.
\end{proof}

\subsection{Generalized Lakshmibai-Seshadri paths}
\subsubsection{}\label{achain}
Let $a$ with $0<a\leq 1$ be a rational number and let $\mu > \nu$ be two weights in $T\lambda$. An $a$-chain
for the pair $(\mu, \nu)$ is a sequence of weights in $T\lambda$ :
\begin{equation*}\nu : = \nu_s\stackrel{\beta_s}{\leftarrow} \nu_{s-1}\stackrel{\beta_{s-1}}{\leftarrow}\dots
\stackrel{\beta_2}{\leftarrow}\nu_1\stackrel{\beta_1}{\leftarrow}\nu_0=:\mu,
\end{equation*}
such that for all $i$ with $1\leq i\leq s$ :\\
(a) $a\beta^{\vee}_i(\nu_i)\in\ \mathbb{N}^+$, if $\beta_i\in\
\Delta_{re}^+$.\\
(b) $a\beta_i^{\vee}(\nu_i)=1$, if $\beta_i\in\ \Delta_{im}$. \\
Observe that if $a=1$, then condition (a) is automatically satisfied.\\
For the above $a$-chain one has
\begin{equation}\label{goodweight} a(\mu-\nu)=\sum\limits_{i=0}^{s-1}a(\nu_i-\nu_{i+1}) =
-\sum\limits_{i=0}^{s-1}a\beta_i^{\vee}(\nu_i)\beta_i\in\ -Q^+.
\end{equation}
\subsubsection{}\label{defgls}
Suppose we have :
\begin{enumerate}
\item $\boldsymbol{\lambda} = (\lambda_1 > \lambda_2 >\dots > \lambda_s)$, a sequence of elements in $T\lambda$,
\item ${\bf{a}} = (a_0=0 <a_1<\dots <a_s=1)$, a sequence of rational numbers,
\end{enumerate}
and set $\pi := (\boldsymbol{\lambda}, {\bf{a}})$ to be the path :
\begin{equation}\label{gls}\pi(t) = \sum\limits_{i=1}^{j-1}(a_i-a_{i-1})\lambda_i + (t-a_{j-1})\lambda_j,\, a_{j-1}\leq t\leq a_j.
\end{equation}
A Generalized Lakshmibai-Seshadri path $\pi = ({\boldsymbol{\lambda}}, {\bf{a}})$ of shape $\lambda$ is the
path given in (\ref{gls}) such that : \\
(a) there exists an $a_i$-chain for $(\lambda_i, \lambda_{i+1})$ for all $i$ with $1\leq i\leq s$,\\
(b) if $\lambda_s\neq \lambda$ there exists a $1$-chain for $(\lambda_s, \lambda)$. \\
We sometimes write
\begin{equation*}\pi = (\boldsymbol{\lambda},\, {\bf{a}}) = (\lambda_1, \lambda_2,\dots, \lambda_s\,;\,a_0 = 0, a_1,\dots, a_{s-1}, a_s = 1).
\end{equation*}
For short, we write GLS path for Generalized Lakshmibai-Seshadri path.
\subsubsection{}\label{remarks}
{\bf Remarks.}\\
(1) Equation (\ref{gls}) of $\pi$ can be also written as follows. Let $t \in\ [a_{j-1}, a_j]$, then :\\
\begin{equation}\label{otherexpression} \pi(t) = \sum\limits_{i=1}^{j-1}(a_i-a_{i-1})\lambda_i +
(t-a_{j-1})\lambda_j = \sum\limits_{i=1}^{j-1}a_i(\lambda_i-\lambda_{i+1}) + t\lambda_j.
\end{equation}
(2)\label{wt} By equation (\ref{otherexpression}) we have that \\
\begin{equation*} \wt \pi = \pi(1)=\sum\limits_{i=1}^{s}(a_i-a_{i-1})\lambda_i = \sum\limits_{i=1}^{s-1}a_i(\lambda_i-\lambda_{i+1})+\lambda_s.
\end{equation*}
Now by equation (\ref{goodweight}), we have $a_i(\lambda_i-\lambda_{i+1})\in\ -Q^+$ for all $i$ with $1\leq
i\leq s-1$ and $\lambda_s\in\ \lambda-Q^+$. In particular, $\wt \pi$ is an integral weight in the
intersection of
$\lambda-Q^+$ and the convex hull of $T\lambda$.\\

\subsubsection{Example}
Let $A = (-k)$ with $k\geq 0$ be an $1\times 1$ matrix, $\lieg$ the associated generalized Kac-Moody algebra,
$\alpha$ the unique simple (imaginary) root, $r : = r_{\alpha}$. Let $\lambda$ be a dominant weight in the
weight lattice of $\lieg$ such that $\alpha^{\vee}(\lambda)=m>0$. One checks that the only GLS paths of shape $\lambda$ are : \\
\begin{center}
$\begin{array}{lr} \pi_0 =(\lambda ; 0, 1),\\
\pi_1 = (r\lambda, \lambda; 0, \frac{1}{m}, 1),\\
\pi_2= (r^2\lambda, r\lambda, \lambda;0, \frac{1}{m(1+k)}, \frac{1}{m}, 1),\\
\pi_3= (r^3\lambda, r^2\lambda, r\lambda, \lambda;0, \frac{1}{m(1+k)^2}, \frac{1}{m(1+k)}, \frac{1}{m},
1),\\
........\\
\pi_s= (r^s\lambda, r^{s-1}\lambda,\dots, r\lambda, \lambda ; 0, \frac{1}{m(1+k)^{s-1}},
\frac{1}{m(1+k)^{s-2}},\dots, \frac{1}{m(1+k)}, \frac{1}{m}, 1),\\
........
\end{array}$
\end{center}
Recall section \ref{definitionoff} and set $f:=f_{\alpha}$. One further checks that $\pi_i=f^i\pi_{\lambda}$.
Notice that the linear path $(r\lambda)t = (r\lambda ; 0, 1)$ is not always a GLS path unlike the Kac-Moody
case. One sees that $(r\lambda)t$ is a GLS path if and only if $m=1$. Again one sees that $(r^s\lambda)t$ is
a GLS path for all $s\in\ \mathbb{N}$, if and only if $m=1$ and $k=0$.

\subsubsection{} \label{aim}
For all $\lambda\in\ P^+$ we denote by $\mathbb{P}_{\lambda}$ the set of all GLS paths of shape $\lambda$. It
is proven in \cite{L1} that when $I^{im}=\emptyset$ the set $\mathbb{P}_{\lambda}$ is stable under the action
of the root operators $f_i,\, e_i,\, i\in\ I$ defined in sections \ref{definitionoff}, \ref{definitionofere}
and $\mathbb{P}_{\lambda} = \mathcal{F}\pi_{\lambda}$, where $\pi_{\lambda}$ is the linear path
$\pi_{\lambda}(t) = \lambda t = (\lambda;0,1)$.

Recall sections \ref{specialcrystals},  \ref{morphism} and \ref{crystalstructure}; the above imply that $\Plambda$ is a highest weight
crystal and a strict subcrystal of $\mathbb{P}$. Furthermore, it is proven in \cite{J} that $\mathbb{P}_{\lambda}$ is isomorphic
(as a crystal) to the crystal associated with the crystal basis of the (unique) highest weight module
$V(\lambda)$ of highest weight $\lambda$ over the quantized enveloping algebra of a Kac-Moody algebra
$\lieg$. Finally, by \cite[Section 9]{L2} $\car V(\lambda) = \car \mathbb{P}_{\lambda}$.

Our aim is to prove analogous results in the generalized Kac-Moody case. However, this is not
straightforward. Already $\Plambda$ will not be a strict subcrystal of $\mathbb{P}$. This results in a number
of complications, in particular to show that it is a highest weight crystal and with respect to the joining
of paths (section \ref{joining}). The latter is needed to prove that the $\Plambda,\, \lambda\in\ P^+$ form a
closed family of highest weight crystals and as a consequence that this family is unique (section 8). The
proof of the character formula (section 9) poses some particular challenges and is significantly more
complicated.

\subsection{Some integrality properties of the Generalized Lakshmibai-Seshadri paths}
In order to study the action of the operators $e_i,\, f_i,\,i\in\ I$ on the set of GLS paths $\Plambda$ we
need certain preliminary results which we give in this section.

\subsubsection{}
Recall sections \ref{definitionofh}-\ref{definitionofeim}. \\

\begin{lemma}\label{fimgls} Suppose that $\pi = (\lambda_1, \lambda_2,\dots, \lambda_s \,;\,0, a_1, \dots, a_{s-1}, 1)$ is a Generalized Lakshmibai-Seshadri path of shape $\lambda\in\ P^+$ and let $i\in\ I^{im}$. Then the function $h_i^{\pi}$ is increasing, $m_i^{\pi} = 0$ and one of the following is true :
\begin{enumerate}
\item $f_+^i(\pi) = 0,\,\,f_i\pi\neq 0$ and $h_i^{\pi}$ is strictly increasing in a neighbourhood of $0$,
\item $f_+^i(\pi)=1,\,\, f_i\pi=0$ and $h_i^{\pi}=0$.
\end{enumerate}
Moreover, $e_i\pi=0$ if and only if $\alpha_i^{\vee}(\wt \pi)<1-a_{ii}$.
\end{lemma}
\begin{proof} Take $i\in\ I^{im}$; by lemma \ref{Tlambda} and since the $\lambda_j$ are in $T\lambda$, one has that $\alpha_i^{\vee}(\lambda_j)\geq 0$, for all $j$ with $1\leq j\leq
s$. Since $\lambda_j>\lambda_{j+1}$, by lemma \ref{positivityofhim} we obtain $\alpha_i^{\vee}(\lambda_1)\geq
\alpha_i^{\vee}(\lambda_2)\geq \cdots \geq \alpha_i^{\vee}(\lambda_s)\geq 0$. Substitution in
(\ref{otherexpression}) shows that $h_i^{\pi}$ is increasing in $[0, 1]$. Since $h_i^{\pi}(0)=0$, we have
that $m_i^{\pi} = 0$.

Moreover, either $\alpha_i^{\vee}(\lambda_1)>0$, in which case $h_i^{\pi}$ increases strictly in $[0, a_1]$,
or $\alpha_i^{\vee}(\lambda_1)=0$ and so $h_i^{\pi}(t)=0$ for all $t\in\ [0, 1]$. In the first case
$f_+^i(\pi) = 0$ and $f_i\pi\neq 0$. In the second case $f_+^i(\pi)=1$ and so, by definition, $f_i\pi=0$.

Since $h_i^{\pi}$ is increasing and $m_i^{\pi}=0$, one obtains $e_i\pi=0$ if and only if $h_i^{\pi}(1) =
\alpha_i^{\vee}(\wt \pi)<1-a_{ii}$.
\end{proof}

\subsubsection{}
\begin{lemma}\label{simpleroot} Let $\nu : = \nu_s\stackrel{\beta_s}{\leftarrow}\nu_{s-1}\stackrel{\beta_{s-1}}{\leftarrow}\cdots
\nu_1\stackrel{\beta_1}{\leftarrow}\nu_0=:\mu$ and take $\alpha_i\in\ \Pi_{re}$. If $r_i\mu<\mu$ and
$r_i\nu\geq \nu$ or $r_i\mu\leq \mu$ and $r_i\nu> \nu$, then $\alpha_i=\beta_{\ell}$ for some $\ell$, with
$1\leq \ell\leq s$.
\end{lemma}
\begin{proof} Assume that $r_i\mu<\mu$ and $r_i\nu\geq \nu$ and recall that $\nu_t =r_{\beta_{t+1}}\nu_{t+1}$ for all $t$, with $0\leq t\leq s-1$. By the hypothesis, there
exists $\ell$ with $1\leq \ell \leq s$ such that $r_i\nu_t< \nu_t$ for all $t$, with  $1\le t\le \ell-1$ and
$r_i\nu_\ell\ge \nu_\ell$, so then $\alpha_i^{\vee}(\nu_\ell)\ge 0$. By lemma \ref{distances} (1)
with $\nu=\nu_{\ell}$ and $\mu=\nu_{\ell-1}$, one has that $\dist\, (r_i\nu_{\ell-1},
\nu_{\ell})<\dist\,(\nu_{\ell-1}, \nu_{\ell})=1$. This implies that $\nu_{\ell-1}=r_i\nu_{\ell}$ and
$\alpha_i = \beta_{\ell}$. The second case follows similarly using lemma \ref{distances} (2).
\end{proof}

\subsubsection{}
The following lemma is similar to \cite[Lemma 4.3]{L2}. We will give the proof in order to outline the fact
that the real operators behave exactly as in the purely real case.\\

\begin{lemma}\label{leftchain} Let $i\in\ I^{re}$ and let $(\mu, \nu)$ be a pair of weights in $T\lambda$, with $\mu>\nu$.
\begin{enumerate}
\item If $r_i\mu<\mu$ and $r_i\nu\geq \nu$, then $r_i\mu\geq \nu$ and if there exists an $a$-chain for
$(\mu, \nu)$, there exist one for the pair $(r_i\mu, \nu)$.
\item If $r_i\mu\leq \mu$ and $r_i\nu> \nu$, then $\mu\geq r_i\nu$ and if there exists an $a$-chain for
$(\mu, \nu)$, there exist one for the pair $(\mu, r_i\nu)$.
\end{enumerate}
\end{lemma}
\begin{proof} Let $$\nu : = \nu_s\stackrel{\beta_s}{\leftarrow}\nu_{s-1}\stackrel{\beta_{s-1}}{\leftarrow}\cdots \nu_1\stackrel{\beta_1}{\leftarrow}\nu_0=:\mu$$
be an $a$-chain for $(\mu, \nu)$ and suppose that $r_i\mu<\mu$ and $r_i\nu\geq \nu$. By lemma
\ref{simpleroot} there exists $\ell$, with $1\leq \ell\leq s$ such that $\alpha_i = \beta_{\ell}$. We will
prove that :
$$\nu : = \nu_s\stackrel{\beta_s}{\leftarrow}\nu_{s-1}\stackrel{\beta_{s-1}}{\leftarrow}\cdots\stackrel{\beta_{\ell+1}}{\leftarrow} \nu_{\ell}=r_i\nu_{\ell-1}\stackrel{\beta_{\ell-1}'}{\leftarrow}r_i\nu_{\ell-2}\cdots r_i\nu_1\stackrel{\beta_1'}{\leftarrow}r_i\nu_0=:r_i\mu,$$
where $\beta'_t=r_i\beta_t$ for all $t$, with $1\leq t\leq \ell-1$, is an $a$-chain.

First of all, since $r_{\beta_t'}=r_ir_{\beta_t}r_i$ we have that $r_{\beta'_t}r_i\nu_t =
r_ir_{\beta_t}\nu_t=r_i\nu_{t-1}$. Again by lemma \ref{simpleroot}, $r_i\nu_t<\nu_t$ for all $t$ with
$1\leq t\leq \ell-1$, so lemma \ref{distances}  (3) gives $\dist\,(r_i\nu_{t-1}, r_i\nu_t)=\dist\,(\nu_{t-1},
\nu_t)=1$ for all $t$, with $1<t\leq \ell-1$. Finally
$a\beta_t'^{\vee}(r_i\nu_t)=a(r_i\beta_t)^{\vee}(r_i\nu_t)=a\beta_t^{\vee}(\nu_t)$ and $\beta_t\in\
\Delta_{im}$ if and only if $\beta_t'\in\ \Delta_{im}$ which imply that the number
$a\beta_t'^{\vee}(r_i\nu_t)$ is an integer and is equal to $1$ if $\beta_t'$ is imaginary. The proof of the
second statement is similar.
\end{proof}

\subsubsection{}
\begin{corollary}\label{middlechain} Let $i\in\ I^{re}$ and let $(\mu, \nu)$ be a pair of weights in $T\lambda$ such that $\mu>\nu$.
\begin{enumerate}
\item If $r_i\mu< \mu$ and $r_i\nu< \nu$, then $r_i\mu\geq  r_i\nu$ and if there exists an $a$-chain
for $(\mu, \nu)$, there exist one for the pair $(r_i\mu, r_i\nu)$.
\item If $r_i\mu>\mu$ and $r_i\nu> \nu$, then $r_i\mu\geq r_i\nu$ and if there exists an $a$-chain
for $(\mu, \nu)$, there exist one for the pair $(r_i\mu, r_i\nu)$.
\end{enumerate}
\end{corollary}
\begin{proof} This follows by the proof of lemma \ref{leftchain}. For example, to prove (1) take $\nu = \nu_{\ell-1}$ in the previous lemma.
\end{proof}

\subsubsection{}
\begin{lemma}\label{middlechainim} Suppose $\mu,\,\nu\in\ T\lambda$ with $\mu>\nu$ and $i\in\ I^{im}$. If $\alpha_i^{\vee}(\mu) =
\alpha_i^{\vee}(\nu)$, then $r_i\mu>r_i\nu$ and if there exists an $a$-chain for the pair $(\mu, \nu)$, then
there exists one for $(r_i\mu, r_i\nu)$.
\end{lemma}
\begin{proof} Let
$$\nu : = \nu_s\stackrel{\beta_s}{\leftarrow}\nu_{s-1}\stackrel{\beta_{s-1}}{\leftarrow}\cdots \nu_1\stackrel{\beta_1}{\leftarrow}\nu_0=:\mu$$
be an $a$-chain for $(\mu, \nu)$ and $i\in\ I^{im}$. Suppose that $\alpha_i^{\vee}(\mu) =
\alpha_i^{\vee}(\nu)\geq 0$. If equality holds, then $r_i\mu=\mu,\, r_i\nu=\nu$ and so there is nothing to
prove. Thus we can assume $\alpha_i^{\vee}(\mu) = \alpha_i^{\vee}(\nu)> 0$. By lemma \ref{positivityofhim},
$r_i$ commutes with $r_{\beta_j}$ for all $j$ with $1\leq j\leq s$. Hence
$$r_i\nu = r_i\nu_s\stackrel{\beta_s}{\leftarrow}r_i\nu_{s-1}\stackrel{\beta_{s-1}}{\leftarrow}\cdots
r_i\nu_1\stackrel{\beta_1}{\leftarrow}r_i\nu_0=r_i\mu$$ is an $a$-chain for $(r_i\mu, r_i\nu)$. Indeed,
$r_i\nu_{j-1} = r_{\beta_j}r_i\nu_j$ and $a\beta_j^{\vee}(r_i\nu_j) = a\beta_j^{\vee}(\nu_j)$ for all $j$
with $1\leq j\leq s$. Finally, $\dist \,(r_i\nu_j, r_i\nu_{j+1}) = \dist\,(\nu_j, \nu_{j+1})=1$ by lemma
\ref{distanceim}.
\end{proof}

\subsubsection{}
Call a path $\pi$ {\it integral} if for all $i\in\ I$, the minimal value of the function $h_i^{\pi}$ is an
integer, i.e. $\min\,\{h_i^{\pi}(t)\, |\, t\in\ [0, 1]\}\in\ \mathbb{Z}$ for all $i\in\ I$. In other words,
$\pi$ is integral, if $m_i^{\pi} := \min\,\{h_i^{\pi}(t)\cap \mathbb{Z}\,|\,t\in\ [0, 1]\} = \min\,\{h_i^{\pi}(t)\,|\,t\in\
[0,1]\}$. We will prove that a GLS path is integral. For this we need the following preliminary result :\\

\begin{lemma}\label{weakchain} Suppose that $\nu : = \nu_s\stackrel{\beta_s}{\leftarrow}\nu_{s-1}\stackrel{\beta_{s-1}}{\leftarrow}\cdots \nu_1\stackrel{\beta_1}{\leftarrow}\mu := \nu_0$ is an $a$-chain. For all $t,\, 1\leq t\leq s$ one has that $a\beta^{\vee}_t(\mu)\in\ \mathbb{Z}$.
\end{lemma}
\begin{proof} We will prove the assertion by induction on $\dist\,(\mu, \nu)$. First, if $\dist\,(\mu, \nu)=1$ and hence $\mu = r_{\beta}\nu$, the assertion follows by the definition of an $a$-chain. In the general case we have that $\mu = r_{\beta_1}r_{\beta_2}\dots r_{\beta_s}\nu$. That $a\beta_1^{\vee}(\mu)\in\ \mathbb{Z}$ follows again by definition. It remains to prove that $a\beta_i^{\vee}(\mu)\in\ \mathbb{Z}$ for $i\neq 1$. Indeed, one has :
$$a\beta_i^{\vee}(r_{\beta_1}r_{\beta_2}\dots r_{\beta_s}\nu) = a\beta_i^{\vee}(r_{\beta_2}\dots r_{\beta_s}\nu)-a\beta_1^{\vee}(r_{\beta_2}\dots r_{\beta_s}\nu)\beta_i^{\vee}(\beta_1).$$
Notice that $a\beta_i^{\vee}(r_{\beta_2}\dots r_{\beta_s}\nu)\in\ \mathbb{Z}$ by the induction hypothesis to
a strictly shorter $a$-chain in which $\nu_1$ replaces $\nu_0$ and since $a\beta_1^{\vee}(r_{\beta_2}\dots
r_{\beta_s}\nu)\in\ \mathbb{Z}$ by the definition of an $a$-chain. Finally, of course
$\beta_i^{\vee}(\beta_1)\in\ \mathbb{Z}$ and hence the result.
\end{proof}

\subsubsection{}
Let $\pi$ be a GLS path. A key fact is that if $h_i^{\pi}$ attains a local minimum at $t_0\in\ [0, 1]$, then
$h_i^{\pi}(t_0)\in\ \mathbb{Z}$. Since $\pi(0)=0$ and $\pi(1)\in\ P$, it is enough to consider $t_0\in\ ]0,
1[$ and then by lemma \ref{fimgls} to take $i\in\ I^{re}$. We say that $h_i^{\pi}$ attains a left local
minimum at $t_0\in\ ]0, 1[$, if there exists $\varepsilon>0$ such that $h_i^{\pi}$ is strictly decreasing in
$[t_0-\varepsilon, t_0]$ and increasing in $[t_0, t_0+\varepsilon]$. Presenting $\pi$ as in equation
(\ref{gls}) it is obvious that $t_0=a_j$ for some $j$, with $0< j< s$. Moreover,
$\alpha_i^{\vee}(\lambda_j)<0$ and $\alpha_i^{\vee}(\lambda_{j+1})\geq 0$ (equivalently,
$r_i\lambda_j<\lambda_j$ and $r_i\lambda_{j+1}\geq
\lambda_{j+1}$).\\

\begin{lemma}\label{integral} Let $\pi$ be a Generalized Lakshmibai-Seshadri path of shape $\lambda$ and $t_0$ a left local
minimum of $h_i^{\pi}$, for $i\in\ I^{re}$. Then $h_i^{\pi}(t_0)\in\ \mathbb{Z}$.
\end{lemma}
\begin{proof} By the first part of lemma \ref{simpleroot} and the hypothesis one obtains $\alpha_i=\beta_t$
for some $\beta_t$ in the $a_j$-chain for the pair $(\lambda_j, \lambda_{j+1})$. Then
$a_j\alpha_i^{\vee}(\lambda_j)=a_j\beta_t^{\vee}(\lambda_j)\in\ \mathbb{Z}$, by lemma \ref{weakchain}.
Conclude by substituting into equations (\ref{goodweight}) and (\ref{otherexpression}).
\end{proof}

\noindent {\bf Remark.} Define a right local minimum by shifting ``strict'' to the right in the definition of
a left local minimum. Using the second part of lemma \ref{simpleroot} we obtain as in lemma \ref{integral}
above that if $h_i^{\pi}$, with $i\in\ I^{re}$, attains a right local minimum at $t_0\in\ ]0, 1[$, then
$h_i^{\pi}(t_0)\in\ \mathbb{Z}$. In the sequel, a local minimum means a right or left local minimum.

\subsubsection{}
The lemma below immediately follows by definition of a GLS path and equation (\ref{otherexpression}) :\\

\begin{lemma}\label{integralvalues} Let $\pi = (\boldsymbol{\lambda}, {\bf a})$ be a Generalized Lakshmibai-Seshadri path and let $t_0\in\ [0, 1]$ with $a_k<t_0\leq a_{k+1}$ and $i\in\ I$. Then $h_i^{\pi}(t_0)\in\ \mathbb{Z}$ if and only if $t_0\alpha_i^{\vee}(\lambda_{k+1})\in\ \mathbb{Z}$. In particular, if $f_+^i(\pi)=a_j$ and $a_p<f_-^i(\pi)\leq a_{p+1}$ then $a_j\alpha_i^{\vee}(\lambda_{j+1}),\, f_-^i\alpha_i^{\vee}(\lambda_{p+1})\in\ \mathbb{Z}$.
\end{lemma}

\subsubsection{}
Call a path $\pi(t)$ {\it monotone} if for all $i\in\ I$ such that $f_i\pi\neq 0$, the function $h^{\pi}_i$ is strictly increasing in $[f_+^i(\pi), f_-^i(\pi)]$ and for all $t>f_-^i(\pi)$ one has that $h_i^{\pi}(t)\geq
m_i^{\pi}+1$.\\

\begin{lemma}\label{monotone} A Generalized Lakshmibai-Seshadri path is monotone.
\end{lemma}
\begin{proof} By definition of $f_+^i(\pi),\, f_-^i(\pi)$ and by continuity, the function $h_i^{\pi}$ does not attain integral values in the interval $]f_+^i(\pi), f_-^i(\pi)[$.
If $h_i^{\pi}$ is not increasing in $[f_+^i(\pi), f_-^i(\pi)]$ and since
$h_i^{\pi}(f_+^i(\pi))<h_i^{\pi}(f_-^i(\pi))$ it follows that $h_i^{\pi}$ attains a local minimum at some
$t_0\in\ ]f_+^i(\pi), f_-^i(\pi)[$. But then by lemma \ref{integral}, $h_i^{\pi}(t_0)\in\ \mathbb{Z}$ which
contradicts our first observation.

The second part is similar. We can assume $f_-^i(\pi)<1$ and then $h_i^{\pi}(1)\geq
m_i^{\pi}+1=h_i^{\pi}(f_-^i(\pi))$, by definition of $f_+^i(\pi)$. If $h_i^{\pi}(t)<m_i^{\pi}+1$ for some
$t\in\ ]f_-^i(\pi), 1[$, we obtain a local minimum $t_0$ in this interval, with $h_i^{\pi}(t_0)<
m_i^{\pi}+1$, forcing $h_i^{\pi}(t_0)= m_i^{\pi}$ by integrality, and so contradicting the definition of
$f_+^i(\pi)$.
\end{proof}

\section{The crystal structure of the set of Generalized Lakshmibai-Seshadri paths}

\subsection{The action of the $f_i,\, i\in\ I$}
In this section we will show that $\mathbb{P}_{\lambda}$ is stable under the action of the root operators
$f_i,\, i\in\ I$.

\subsubsection{}\label{recallnotation}
Let $\pi = (\boldsymbol{\lambda}, {\bf a}) = (\lambda_1, \lambda_2, \dots, \lambda_s\,;\,0, a_1, \dots,
a_{s-1}, a_s = 1)$ be a GLS path of shape $\lambda$, and recall sections \ref{definitionofh},
\ref{definitionoff}. In this section it is convenient for brevity to just suppose $\lambda_1\geq
\lambda_2\geq \cdots\geq \lambda_s$ and $0=a_0\leq a_1\leq \cdots\leq a_{s-1}\leq a_s=1$. We recover
strictness by just dropping some terms in the expression for $\pi$.

For simplicity in the the rest of this section we set $f_-^i : = f_-^i(\pi)$ and $f_+^i ; = f_+^i(\pi)$.

\subsubsection{}

\begin{proposition}\label{actionfre} Let $i\in\ I^{re}$ and $\pi$ as above and such that $f_i\pi\ne 0$. Then $f_+^i=a_t$, for some $t$, with $1\le t< s$ and $a_{p-1}<f_-^i\le a_p$ for some $p$, with $t+1\le p\le s$. The path $f_i\pi$ is equal to
$$f_i\pi = (\lambda_1, \dots, \lambda_t, r_i\lambda_{t+1},\dots, r_i\lambda_p, \lambda_p,\dots \lambda_s ; a_0, a_1, \dots, a_{p-1}, f_-^i, a_p,\dots, a_s
),$$
and is a Generalized Lakshmibai-Seshadri path of shape $\lambda$. In particular, the set of Generalized
Lakshmibai-Seshadri paths of shape $\lambda$ is stable under the action of $f_i,\, i\in\ I^{re}$.
\end{proposition}
\begin{proof}
The proof of this proposition is exactly as in \cite[Proposition 4.7]{L2}, but we still give the proof for
completeness. The fact that $f_+^i=a_t$ for some $t$, with $1\le t<s$ follows by the integrality of $\pi$ proven in lemma \ref{integral}-then it is clear that the resulting path is of the above form. The only thing one has to check
is that conditions (a), (b) of section \ref{defgls} still hold for this new path. By monotonicity of $\pi$
(see lemma \ref{monotone}) one has that $r_i\lambda_k>\lambda_k$ for all $k$, with $t<k\leq p$ and so by
corollary \ref{middlechain} (2) there exists an $a_k$-chain for $(r_i\lambda_k, r_i\lambda_{k+1})$ for all $k$,
with $t<k<p$. On the other hand, $r_i\lambda_t\leq \lambda_t$ and lemma \ref{leftchain} implies that there
exists an $a_t$-chain for $(\lambda_t, r_i\lambda_{t+1})$.  Finally, since $h_i^{\pi}(f_-^i)\in\ \mathbb{Z}$,
there exists an $f_-^i$-chain for $(r_i\lambda_p, \lambda_p)$.
\end{proof}

\subsubsection{}\label{imaginary}

\begin{proposition}\label{actionfim} Let $i\in\ I^{im}$, $\pi$ as in Section \ref{recallnotation} and such that $f_i\pi\ne 0$. Then $f_+^i=0$ and $a_{p-1}<f_-^i\le a_p$, for some $p$ with $1\le p\le s$ and the path $f_i\pi$ is equal to :
$$f_i\pi=(r_i\lambda_1, \dots, r_i\lambda_p, \lambda_p,\dots, \lambda_s ; a_0, a_1, \dots, a_{p-1}, f_-^i, a_p,\dots,
a_s),$$ with $\alpha_i^{\vee}(\lambda_j)=\alpha_i^{\vee}(\lambda_{j+1})$ for all $j$, with $1\leq j\leq p-1$
and is a Generalized Lakshmibai-Seshadri path. In particular, the set of Generalized Lakshmibai-Seshadri
paths of shape $\lambda$ is stable under the action of $f_i,\, i\in\ I^{im}$.
\end{proposition}
\begin{proof} By lemma \ref{fimgls}, since $f_i\pi\neq 0$, then $f_+^i=a_0=0$. Let $p$, with $1\leq p\leq s$ be such that $a_{p-1}<f_-^i\leq a_p$. By lemma \ref{fimgls} and the definition of $f_-^i$, the function $h_i^{\pi}$ is strictly increasing in the interval $[0, f_-^i]$. Thus $r_i\lambda_p>\lambda_p\geq \lambda_{p+1}$ and so the resulting path will be of the above form. We need to prove that $f_i\pi$ is a GLS path.
For this purpose it is sufficient to show the following :
\begin{enumerate}
\item For all $j$ with $1\leq j\leq p-1$ there exists an $a_j$-chain for the pair $(r_i\lambda_j, r_i\lambda_{j+1})$.
\item There exists an $f_-^i$-chain for the pair $(r_i\lambda_p, \lambda_p)$.
\end{enumerate}
Recall equation (\ref{otherexpression}); one has that
\[\pi(f_-^i)=\sum\limits_{k=1}^{p-1}a_k(\lambda_k-\lambda_{k+1})+f_-^i\lambda_p,\]
so that
\[1=h_i^{\pi}(f_-^i) = \sum\limits_{k=1}^{p-1}a_k\alpha_i^{\vee}(\lambda_k-\lambda_{k+1})+f_-^i\alpha_i^{\vee}(\lambda_p).\]
Since $\lambda_k>\lambda_{k+1}$ for all $k$ with $1\leq k\leq p-1$, by lemma \ref{positivityofhim} we have
that $\alpha_i^{\vee}(\lambda_k-\lambda_{k+1})\geq 0$. On the other hand,
$f_-^i\alpha_i^{\vee}(\lambda_p)\in\ \mathbb{Z}$ by lemma \ref{integralvalues} and is strictly positive by
monotonicity. This forces $\alpha_i^{\vee}(\lambda_k)=\alpha_i^{\vee}(\lambda_{k+1})$ for all $k$ with $1\leq
k\leq p-1$. Hence by lemma \ref{middlechainim} and since there exists an $a_k$-chain for $(\lambda_k,
\lambda_{k+1})$, there exists one for $(r_i\lambda_k, r_i\lambda_{k+1})$ for all $k$ with $1\leq k\leq p-1$
and (1) follows.

Finally, lemma \ref{integralvalues} gives that $f_-^i\alpha_i^{\vee}(\lambda_p)\in\ \mathbb{Z}$ and so
there exists an $f_-^i$-chain for $(r_i\lambda_p, \lambda_p)$ and (2) follows.
\end{proof}

\subsection{The action of the $e_i, i\in\ I$}\label{e}
In this section we study the action of the root operators $e_i,\, i\in\ I$ on the set $\Plambda$ of GLS paths
of shape $\lambda$.

\subsubsection{}
Let $\pi$ be as in section \ref{recallnotation}, $i\in\ I^{re}$ and recall section \ref{definitionofere}. By
lemma \ref{integral}, we have that $e_+^i:=e_+^i(\pi)=a_k$ for some $k$ with $1\leq k\leq s$. Then
$e_i\pi\neq 0$ if and only if $e_+^i>0$. Let $e_i\pi\neq 0$, then we can assume that $e_-^i:=e_-^i(\pi)$ is
such that $a_{q-1}\leq e_-^i<a_q$, for some $q$ with $0\leq q\leq k$. The proof of the proposition below is
similar to the proof of proposition \ref{actionfre} and so we omit it.\\

\begin{proposition}\label{ere} Assume that $i\in\ I^{re}$ and that $e_i\pi\neq 0$. Then $e_+^i> 0$ and hence the number $e_-^i\in\ [0, e_+^i]$ is defined. Let $e_+^i =a_k$ and $a_{q-1}\leq e_-^i<a_q$ for $1\leq q<k$. Then the path $e_i\pi$ is equal to :\\
$$e_i\pi=(\lambda_1, \dots, \lambda_q, r_i\lambda_q,\dots, r_i\lambda_k, \lambda_{k+1},\dots \lambda_s ; a_0, \dots, a_{q-1}, e^i_-, a_q,\dots,
a_s),$$
and is a Generalized Lakshmibai-Seshadri path of shape $\lambda$. In particular, the set of
Generalized Lakshmibai-Seshadri paths of shape $\lambda$ is stable under the action of $e_i,\, i\in\ I^{re}$.
\end{proposition}

\subsubsection{}\label{eim}
Let $\pi\in\ \Plambda$ and $i\in\ I^{im}$. It can happen that $e_i\pi\in\ \mathbb{P}\setminus
\mathbb{P}_{\lambda}$. Indeed, let $\lambda$ be such that $\alpha_i^{\vee}(\lambda)\geq 1-a_{ii}$. Recall
that for all $\pi\in\ \Plambda$ one has that $\wt \pi\in\ \lambda-Q^+$ and so $\alpha_i^{\vee}(\wt \pi)\geq
\alpha_i^{\vee}(\lambda)\geq 1-a_{ii}$. Then by lemma \ref{fimgls} we obtain that $e_i\pi\neq 0$ for all
$\pi\in\ \Plambda$. In particular, $e_i\pi_{\lambda}\neq 0$. But $\wt e_i\pi_{\lambda}=\lambda+\alpha_i$ and
so $e_i\pi\notin \Plambda$.

\subsection{Crystal Structure of $\mathbb{P}_{\lambda}$}

\subsubsection{}\label{imposeeim}
For all $i\in\ I$ we set $e_i\pi=0$ if and only if $e_i\pi\notin \mathbb{P}_{\lambda}$. For real indices,
$e_i\pi\notin \Plambda$ is equivalent to $e_i\pi = 0$ in $\mathbb{P}$. Notice that this means that the ``only
if'' of the last statement of lemma \ref{fimgls} will henceforth be violated. Recall the notation of section
\ref{aim}. Our aim is to show that $\Plambda = \mathcal{F}\pi_{\lambda}$.

\subsubsection{}\label{strict}
Recall that $\mathbb{P}$ has a crystal structure with crystal operations $\wt, e_i, f_i, \varepsilon_i,
\varphi_i$ for all $i\in\ I$ defined in section \ref{crystalstructure}. Consider the embedding $\psi :
\Plambda \hookrightarrow\mathbb{P}$. Then $\Plambda$ is a subcrystal of $\mathbb{P}$. However, it is not a
strict subcrystal of $\mathbb{P}$. Indeed, by propositions \ref{actionfre}, \ref{actionfim} and \ref{ere},
the map $\psi$ commutes with all the crystal operations except the $e_i,\, i\in\ I^{im}$, though we still
have $e_i\psi(\pi)=\psi (e_i\pi)$ if $e_i\pi\neq 0$, by definition \ref{crystal} (4) and because $f_i$
commutes with $\psi$.

\subsubsection{}\label{iff}
Take $i\in\ I^{im}$. Recall that we have set $\varepsilon_i(\pi)=0$ (section \ref{crystalstructure}). By
remark \ref{remarks} (2) one has that $\wt \pi\in\ \lambda-Q^+$ and so $\alpha_i^{\vee}(\wt \pi)\geq 0$. Finally,
$$f_i\pi=0\Leftrightarrow f_+^i(\pi)=1\Leftrightarrow \alpha_i^{\vee}(\wt \pi)=0\Leftrightarrow
\varphi_i(\pi)=0.$$ We conclude that $\Plambda\in\ \mathcal{B}$ (see section \ref{categoryB}).

\subsubsection{}
We will show that $\Plambda$ is a highest weight crystal (proposition \ref{unique}). For this we need the
following preliminary lemma.

Given a reduced decomposition $w=r_{i_1}r_{i_2}\cdots r_{i_t}$ of $w\in\ W$, set $\operatorname{Supp}(w) =
\{\alpha_{i_k},\,|\,1\leq k\leq t\}$. As is well-known it is independent of the choice of reduced
decomposition. \\

\begin{lemma}\label{chainfordom} Let $\mu, \nu\in\ T\lambda$ with $\lambda, \mu\in\ P^+$ and suppose that $\nu\stackrel{\beta}{\leftarrow}\mu$. Then $\beta$ is a simple imaginary root.
\end{lemma}
\begin{proof} By hypothesis $\mu=r_{\beta}\nu$ for some $\beta\in\ W\Pi\cap \Delta^+$ with $\beta^{\vee}(\nu)>0$. Then $\mu$ being dominant implies that
$\beta\in\ W\Pi_{im}$.

Let $\beta=w\alpha_i$ and $i\in\ I^{im}$ and suppose that $w\notin \operatorname{Stab}_W(\alpha_i)$. Then
$\mu=r_{\beta}\nu=wr_iw^{-1}\nu$. By corollary \ref{formofdom}, and since $\mu$ is dominant, every minimal representative of $\mu$ starts with $r_j,\, j\in\ I^{im}$. In particular, $w\in\
\operatorname{Stab}_W(r_iw^{-1}\nu)$. Recalling that $r_iw^{-1}\nu$ is dominant, this by \cite[Proposition
3.12]{Kac} implies that for every root $\alpha_j$ in $\operatorname{Supp}(w)$ one has that
$\alpha_j^{\vee}(r_iw^{-1}\nu)=0$. Then
$$\alpha_j^{\vee}(r_iw^{-1}\nu)=\alpha_j^{\vee}(w^{-1}\nu)-\alpha_i^{\vee}(w^{-1}\nu)a_{ji}=0,$$
and since $\alpha_i^{\vee}(w^{-1}\nu) = \beta^{\vee}(\nu)>0$ we must have that $\alpha_j^{\vee}(w^{-1}\nu)\leq 0$ for all $\alpha_j\in\ \operatorname{Supp}(w)$.

Let $w=r_kw_1$, with $\ell(w)=\ell(w_1)+1$. If $\alpha_k^{\vee}(\nu)=0$ then $w^{-1}\nu=w_1^{-1}\nu$ and
since $r_k\mu=\mu$ we can choose $\beta = w_1\alpha_i$.

In the above manner, we are reduced to the case where $\alpha_k^{\vee}(\nu)\neq 0$. Note that
$w_1^{-1}\alpha_k\in\ \Delta^+$ and that we can write $(w_1^{-1}\alpha_k)^{\vee}=\sum\limits_{\alpha_j\in\
\operatorname{Supp}(w)}n_j\alpha_j^{\vee}$ with $n_j\geq 0$ for all $j$. Then :

\begin{equation*} \alpha_k^{\vee}(\nu)=(w^{-1}\alpha_k)^{\vee}(w^{-1}\nu) =-(w_1^{-1}\alpha_k)^{\vee}(w^{-1}\nu)
=-\sum\limits_{\alpha_j\in\ \operatorname{Supp}(w)}n_j\alpha_j^{\vee}(w^{-1}\nu)\geq 0,
\end{equation*}
which by assumption on $r_k$ forces $\alpha_k^{\vee}(\nu)>0$ and so $r_k\nu>\nu$. Set $\beta_1=w_1\alpha_i$,
then $\beta_1^{\vee}(r_k\nu)=\alpha_i^{\vee}(w^{-1}\nu)>0$ and consequently $r_{\beta_1}r_k\nu =
w_1r_iw^{-1}\nu=\mu>r_k\nu$. We conclude that $\mu>r_k\nu>\nu$, which implies that $\dist (\mu, \nu)\geq 2$
which contradicts our hypothesis. Hence $w=\id, \mod \operatorname{Stab}_W(\alpha_i)$ and $\beta=\alpha_i\in\
\Pi_{im}$.
\end{proof}

\subsubsection{}
Let $\mathcal{E}$ be the monoid generated by the $e_i,\, i\in\ I$ and set $\Plambda^{\mathcal{E}}=\{\pi\in\
\Plambda\,|\, e_i\pi=0, \, \mbox{for all}\, i\in\ I\}$.\\

\begin{proposition}\label{unique} Let $\pi\in\ \mathbb{P}_{\lambda}$. Then $e_i\pi=0$ for all $i\in\ I$ if and only
if $\pi = \pi_{\lambda}$, that is $\Plambda^{\mathcal{E}}=\{\pi_{\lambda}\}$. Moreover,
$\mathbb{P}_{\lambda}=\mathcal{F}\pi_{\lambda}$.
\end{proposition}
\begin{proof} It is clear that $e_i\pi_{\lambda}=0$ for all $i\in\ I$, since for all $\pi\in\ \Plambda$, one has that $\wt \pi\prec \lambda$. Let
$$\pi = (\lambda_1, \lambda_2, \dots, \lambda_s:0, a_1, \dots, a_{s-1}, a_s=1)$$ be a path in
$\mathbb{P}_{\lambda}$ and notice that $\pi = \pi_{\lambda}$ if and only if $\lambda_1 = \lambda$.  Suppose
that $e_i\pi\notin \Plambda$, for all $i\in\ I$. Since for $i\in\ I^{re}$, the $e_i$ preserve the set
$\Plambda$ our assumption implies that $e_i\pi=0$ for all $i\in\ I^{re}$. This means that
$\alpha_i^{\vee}(\lambda_1)\geq 0$, for all $i\in\ I^{re}$, that is $\lambda_1$ is dominant (and different
from $\lambda$). On the other hand, by definition of a GLS path there exists an $a_1$-chain
$$\lambda_2:=\nu_s \stackrel{\beta_s}{\leftarrow}\cdots
\nu_1\stackrel{\beta_1}{\leftarrow}\nu_0=:\lambda_1.$$ (If $a_1=1$ we set $\lambda_2=\lambda$.) Then by lemma
\ref{chainfordom} we must have that $\beta_1=\alpha_i$ for some $i\in\ I^{im}$. Hence
$a_1\alpha_i^{\vee}(\nu_1)=1$ and applying proposition \ref{actionfim} we have that
$$\pi'=(\nu_1, \lambda_2, \dots, \lambda_s;0, a_1, \dots, a_s=1)\in \Plambda$$
is such that $f_i\pi'=\pi$ and so $e_i\pi\in\ \mathbb{P}_{\lambda}$.
We conclude that the only path in $\Plambda$ killed by all the $e_i,\,i\in\ I$ is $\pi_{\lambda}$.

We will prove now that $\Plambda=\mathcal{F}\pi_{\lambda}$. Since $\pi_{\lambda} = (\lambda ; 0, 1)\in\
\Plambda$ and $\Plambda$ is stable under the action of the $f_i,\, i\in\ I$ by propositions \ref{actionfre}
and \ref{actionfim}, one obtains $\mathcal{F}\pi_{\lambda}\subset \Plambda$. For the reverse inclusion it is
enough by definition \ref{crystal} (4) to show that $\pi_{\lambda}\in\ \mathcal{E}\pi$, for all $\pi\in\
\Plambda$. Since $\wt \pi\in\ \lambda-Q^+$ and $\wt (e_i\pi)=\wt \pi+\alpha_i$, we obtain $\mathcal{E}\pi\cap
\Plambda^{\mathcal{E}}\neq \emptyset$ and so the assertion follows from the first part.
\end{proof}

\section{Closed Families of Highest Weight Crystals}
Call a family $\{B(\lambda)|\lambda\in\ P^+\}$ of highest weight crystals {\it closed under tensor products}
or simply {\it closed} if for all $\lambda,\, \mu\in\ P^+$ the element $b_{\lambda}\otimes b_{\mu}$ of
$B(\lambda)\otimes B(\mu)$ generates a crystal isomorphic to $B(\lambda +\mu)$. Our aim now is to prove that
the family $\{\Plambda\,|\,\lambda\in\ P^+\}$ is closed.

Let $\lambda, \mu\in\ P^+$ and set $\nu : =\lambda+\mu\in\ P^+$. We need to show that the crystals generated
by $\pi_{\lambda}\otimes \pi_{\mu}$ and $\pi_{\nu}$ are isomorphic. As in \cite{L2}, the proof involves
deforming the path $\pi_{\lambda}\otimes \pi_{\mu}$ to $\pi_{\nu}$ without changing the crystal graph it
generates. To do this we need to introduce some operations on $\mathbb{P}$. The fact that the crystals
$\Plambda$ and the crystal generated by $\pi_{\lambda}\otimes \pi_{\mu}$ are not strict subcrystals of
$\mathbb{P}$ causes some significant extra difficulty.

\subsection{Deformations of paths}

\subsubsection{The join of two paths}\label{fusion}
Let $s\leq s'$ be two rational numbers in $[0, 1]$, $\theta$ the trivial path defined by
$\theta(t)=0$ for all $t\in\ [0,1]$ and let $\pi, \pi'\in\ \mathbb{P}$. Define $\pi*\theta_s^{s'}*\pi'$ to be
the path:
$$(\pi*\theta_s^{s'}*\pi')(t) =
\left\{
\begin{array}{lr}\pi(t), t\in\ [0, s],\\
\pi(s), t\in\ [s, s'],\\
\pi (s) + \pi'(t)-\pi'(s'), t\in\ [s', 1].
\end{array}
\right.$$ It is the concatenation of the truncated paths $\pi^s(t) :[0, s]\rightarrow \mathbb{Q}P,\,
\pi'_{s'}(t):[s', 1]\rightarrow \mathbb{Q}P$ and the trivial path $\theta$. Clearly, if $s = s'$ and $\pi =
\pi'$, then $\pitheta = \pi$. The reason for introducing this operation is explained in the section below.

\subsubsection{}\label{pix}
Take $\lambda,\, \mu\in\ P^+$. We recall that by our conventions $[0, 1]\subset \mathbb{Q}$. Let $x\in\ [0,
1]$ and set $\pi^x = (1-x)\pi_{\lambda}\otimes \pi_{\mu}+x\pi_{\nu}$. Then $\pi^x\in\ \mathbb{P}$ with $\wt
\pi^x = \nu$ for all $x\in\ [0, 1]$ and $\pi^0=\pi_{\lambda}\otimes \pi_{\mu},\,\pi^1=\pi_{\nu}$. One can
write $\pi^x = \pi_{\delta}\otimes \pi_{\delta'}$, where $\delta = (1-x)\lambda+\frac{1}{2}x\nu$ and $\delta'
= (1-x)\mu+\frac{1}{2}x\nu$. Of course $\delta+\delta' = \nu$, but $\delta, \delta'$ are not in general in the
weight lattice and thus $\pi_{\delta}, \pi_{\delta'}$ are not in $\mathbb{P}$. However, one can find a
positive integer $r$, such that $r\delta,\, r\delta'\in\ P$. Then $\pi^x =
\pi_{r\delta}*\theta_{1/r}^{1-1/r}*\pi_{r\delta'}$ up to parametrization.

In section \ref{pathdist} we give sufficient conditions for any two paths $\pi, \,\pi'$ to generate
isomorphic crystals. Then, in sections \ref{joining} and \ref{distofpix}, we show that the set of paths
$\{\pi^x, \,x\in\ \mathbb{Q}\}$ satisfies these conditions, and in particular that
$\mathcal{F}(\pi_{\lambda}\otimes \pi_{\mu})$ is a highest weight crystal isomorphic to
$\mathcal{F}\pi_{\nu}=\mathbb{P}_{\nu}$.

\subsection{Distance of paths}\label{pathdist}
\subsubsection{} Let $\mathcal{A}$ denote the monoid generated by the $e_i,\,f_i\in\ I$ and let $J\subset I$ be a finite subset of $I$. Denote by $\mathcal{A}_J,\,\mathcal{F}_J$ the monoids generated by the $e_i, f_i,\,i\in\ J $ and $f_i,\, i\in\ J$ respectively. Clearly, $\mathcal{A}_J\subset \mathcal{A}$ and $\mathcal{F}_J\subset \mathcal{F}$. Set $c_J = \max \{|a_{ij}|,\, i, j\in\ J\}$. For all $\pi, \,\pi'\in\ \mathbb{P}$, define their $J$-{\it distance} $d_J(\pi, \pi')$ to be :
$$d_J(\pi, \pi') = \max \{|\alpha_i^{\vee}(\pi(t)- \pi'(t))|,\, t\in\ [0,1], i\in\ J\}.$$

\subsubsection{}
The following lemma is the initial step in establishing the isomorphism theorem.\\

\begin{lemma}\label{distance} Let $\pi,\,\pi'$ be integral and monotone paths such that $d_J(\pi,\,\pi')<\epsilon<1$. Then for all $i\in\ J$ one has :\\
(1) $m^{\pi}_i = m_i^{\pi'}$ and $\alpha_i^{\vee}(\wt \pi)=\alpha_i^{\vee}(\wt \pi')$.\\
(2) If $f_i\pi \neq 0$, then $f_i\pi'\neq 0$ and $d_J(f_i\pi,\, f_i\pi')<2c_J\epsilon$.\\
(3) For $i\in\ I^{re}\cap J$, if $e_i\pi\neq 0$, then $e_i\pi'\neq 0$. If $e_i\pi,\, e_i\pi'\neq 0$ then
$d_J(e_i\pi,\, e_i\pi')<2c_J\epsilon$.
\end{lemma}
\begin{proof} Statement (1) is an immediate consequence of the definitions and integrality. By section \ref{crystalstructure} and (1) we obtain
$\varepsilon_i(\pi)=\varepsilon_i(\pi')$ and $\varphi_i(\pi)=\varphi_i(\pi')$ and thus the first part of (2)
and (3) follow by normality for $i\in\ I^{re}$. For $i\in\ I^{im}$ the first part of (2) follows by (1) and
section \ref{iff}. The second part of (2) follows exactly as in \cite[Lemma 6.4.25]{J}. A key point is to
show that the intervals $[f_+^i(\pi), f_-^i(\pi)]$ and $[f_+^i(\pi'), f_-^i(\pi')]$ have non-empty
intersection. A similar comment applies to the second part of (3).
\end{proof}

\noindent {\bf Remark.} Notice that we do not obtain that $e_i\pi\neq 0$ implies $e_i\pi'\neq 0$ as it does
for real indices since the ``only if'' of lemma \ref{fimgls} is violated (see section \ref{imposeeim}). This
leads to an extra difficulty ultimately resolved by lemma \ref{unique2}.

\subsection{Joining Generalized Lakshmibai-Seshadri paths}\label{joining}
Throughout this section fix $\lambda,\, \mu\in\ P^+$.

\subsubsection{}
Recall the definition of an $a$-chain \ref{achain}. We call the $a$-chain in \ref{achain} {\it weak}, if $a\beta_i^{\vee}(\nu_i)\in \mathbb N^+$, for all $i$, with $1\le i\le s$. Clearly, an $a$-chain is also a weak $a$-chain. Note that the condition of being a weak $1$-chain is empty. We denote by $\widehat{\mathbb P}_\lambda$ the set of paths given by (\ref{gls}) such that there exists a weak $a_i$-chain for $(\lambda_i, \lambda_{i+1})$, for all $i$, with $1\le i\le k$. Notice that since $\lambda_s\in T\lambda,\, \lambda_s\ge \lambda$. Clearly, $\Plambda\subset \widehat{\mathbb P}_{\lambda}$.

The set $\widehat{\mathbb P}_{\lambda}$ is stable under the $f_i,\, i\in I$. To show this one may first note that the only possible change is when $i \in I^{im}$. Recall Sections \ref{definitionofh}-\ref{definitionoff} the definition of the function $h_i^{\pi}$ and the action of $f_i$ on a path $\pi$. Since we still retain Bruhat
sequences in the definition of $\widehat{\mathbb P}_{\lambda}$, it
follows that the function $t\mapsto h_i^\pi(t)$ for $\pi \in
\widehat{\mathbb P}_\lambda$ is increasing for $i \in I^{im}$ as
in Lemma \ref{fimgls}. Thus even for these more general
paths we must still have $h^\pi_i(f^i_-(\pi))=1$ and it was this
we required to obtain Proposition \ref{actionfim}.

\subsubsection{}
Recall section \ref{fusion}. We will join under certain conditions paths in $\widehat{\mathbb{P}}_{\lambda}$ with paths in
$\widehat{\mathbb{P}}_{\mu}$.

Let $\tau\in\ T$ and suppose that $\tau\mu>\mu$. By definition we may write $\tau\mu = r_{\beta_1}\cdots
r_{\beta_s}\mu$ with $\beta_t^{\vee}(r_{\beta_{t+1}}\cdots r_{\beta_s}\mu)>0$, for all $t$, with $1\leq t\leq
s$. By Lemma \ref{twodominant}, one has that $\beta_t^{\vee}(r_{\beta_{t+1}}\cdots r_{\beta_s}\lambda)\geq 0$ for
all $t$ and so $r_{\beta_t}\cdots r_{\beta_s}\lambda\geq r_{\beta_{t+1}}\cdots r_{\beta_s}\lambda$. In the
expression $\tau\lambda=r_{\beta_1}\cdots r_{\beta_s}\lambda$, omit the $r_{\beta_t}$ if
$\beta_t^{\vee}(r_{\beta_{t+1}}\cdots r_{\beta_s}\lambda)=0$, that is if $r_{\beta_t}\in\
\operatorname{Stab}_T(r_{\beta_{t+1}}\cdots r_{\beta_s}\lambda)$, and denote by $\overline{\tau}$ the new
element in $T$. One has $\overline{\tau}\lambda=\tau\lambda$ and $\overline{\tau}\lambda\geq \lambda$. Notice
that if $\tau_1\tau_2\mu>\tau_2\mu>\mu$ then $\overline \tau_1 \overline
\tau_2\lambda=\overline{\tau_1\tau_2}\lambda\geq \overline \tau_2 \lambda\geq \lambda$.

\subsubsection{}\label{gluingpair}
{\bf Definition.} Fix two rational numbers $0< s\leq s'< 1$ and let
\begin{equation*}\pi = (\lambda_1, \dots, \lambda_k;0, a_1, \dots, a_{k-1}, 1)\in\ \widehat{\mathbb{P}}_{\lambda}, \,\pi'=(\mu_1, \dots, \mu_{\ell};0, b_1, \dots, b_{\ell-1}, 1)\in\
\widehat{\mathbb{P}}_{\mu},
\end{equation*}
be such that $a_{k-1}<s\leq s'<b_1$. Observe that by equation (\ref{otherexpression}), $\pi(t)$ is a
translate of $(t-a_{k-1})\lambda_k$ in $[a_{k-1}, s]$ and $\pi'(t)=t\mu_1$ in $[s', b_1]$.

We will allow $b_t=b_{t+1}$ if necessary, so that $\dist\, (\mu_t, \mu_{t+1})=1$ for all $t$, with $1\leq
t\leq \ell-1$; we may assume that $\mu_t=r_{\beta_t}\mu_{t+1}$ with $\mu_{\ell}=\mu$ and thus get a sequence of positive roots $\beta_1, \beta_2, \dots, \beta_{\ell-1}$. Note that one may be able to write the path $\pi'$ in different ways, which will give different sequences of positive roots. For any such sequence, set $\tau =
r_{\beta_1}r_{\beta_2}\cdots r_{\beta_{\ell-1}}$. Then $\mu_1=\tau\mu\geq \mu$. We say that the paths
$\pi,\,\pi'$ can be properly joined across $[s, s']$ if there exists a sequence $\beta_1, \beta_2, \dots, \beta_{\ell-1}$ as above for which the following two conditions hold :
\begin{enumerate}
\item $\lambda_k\geq \overline{\tau}\lambda$
and if $\lambda_k>\overline \tau \lambda$ there exists an $s$-chain for the pair $(\lambda_k,
\overline{\tau}\lambda)$.
\item For all $t$, with $1\leq t\leq \ell-1$, if $\beta_t\in\ \Delta_{im}$ one has that $s\beta_t^{\vee}(r_{\beta_{t+1}}\cdots
r_{\beta_{\ell-1}}\lambda)<1$.
\end{enumerate}
We call (1) and (2) the joining conditions. Note that it is enough to consider the second condition for the
roots $\beta_t$ appearing in $\overline \tau$, that is $r_{\beta_t}\notin \Stab_T(r_{\beta_{t+1}}\cdots
r_{\beta_{\ell-1}}\lambda)$. We may write
\begin{equation}\label{pitheta} \pitheta = (\lambda_1, \lambda_2,\dots, \lambda_k, 0, \mu_1, \mu_2,\dots, \mu_{\ell};0, a_1,\dots, a_{k-1}, s, s', b_1,\dots, b_{\ell-1}, 1),
\end{equation}
where we interpret the right hand side as a path using (\ref{gls}).
We denote by $\Pithetahat$ the set of paths $\pitheta$ were $\pi\in\
\widehat{\mathbb{P}}_{\lambda},\,\pi'\in\ \widehat{\mathbb{P}}_{\mu}$ can be properly joined across $[s,s']$ and are such that $\pitheta(1)\in P$. Of course if
$\lambda =\mu$ and $s=s'$, the set $\Pithetahat$ is equal to
$\widehat{\mathbb{P}}_{\lambda}$.

\noindent {\bf Remark.} Let $\mu_{\ell}\stackrel{\beta_{\ell-1}}{\leftarrow}\mu_{\ell-1}\cdots
\mu_2\stackrel{\beta_1}{\leftarrow}\mu_1$ and suppose that the second joining condition holds with
$\tau=r_{\beta_1}r_{\beta_2}\cdots r_{\beta_{\ell-1}}$ as specified above. Assume $i\in\ I^{re}$.\\

\noindent (1) If $r_i\mu_t>\mu_t$ for all $t$, with $1\leq n\leq t\leq m\leq \ell$ and $r_i\mu_{n-1}\leq
\mu_{n-1}$, then
$$\mu_{\ell}\stackrel{\beta_{\ell-1}}{\leftarrow}\mu_{\ell-1}\cdots \mu_m\stackrel{\alpha_i}{\leftarrow}
r_i\mu_m\stackrel{r_i\beta_{m-1}}{\leftarrow} \cdots
r_i\mu_{n+1}\stackrel{r_i\beta_n}{\leftarrow}r_i\mu_n=\mu_{n-1}\cdots\mu_2\stackrel{\beta_1}{\leftarrow}\mu_1.$$
As above, this specifies the element $\tilde \tau =r_{\beta_1}\cdots r_{\beta_{n-2}}r_{\beta_n'}\cdots
r_{\beta_{m-1}'}r_ir_{\beta_m}\cdots r_{\beta_{\ell-1}}$, where $\beta'_t=r_i\beta_t$, for all $t$, with
$n\leq t\leq m-1$, relative to which the second joining condition holds because no new scalar products
appear.\\

\noindent (2) If $r_i\mu_t<\mu_t$ for all $t$, with $1\leq n\leq t\leq m\leq \ell$ and $r_i\mu_{m+1}\geq
\mu_{m+1}$, then
$$\mu_{\ell}\stackrel{\beta_{\ell-1}}{\leftarrow}\mu_{\ell-1}\cdots
\mu_{m+1}=r_i\mu_m\stackrel{r_i\beta_{m-1}}{\leftarrow} \cdots
r_i\mu_{n+1}\stackrel{r_i\beta_n}{\leftarrow}r_i\mu_n
\stackrel{\alpha_i}{\leftarrow}\mu_n\cdots\mu_2\stackrel{\beta_1}{\leftarrow}\mu_1.$$ Similarly, relative to
$\tilde \tau =r_{\beta_1}\cdots r_{\beta_{n-1}}r_ir_{\beta_n'}\cdots r_{\beta_{m-1}'}r_{\beta_{m+1}}\cdots
r_{\beta_{\ell-1}}$, where $\beta'_t=r_i\beta_t$, for all $t$, with $n\leq t\leq m-1$, the second joining
condition holds.

\subsubsection{}
The subsets $\Pithetahat$ of $\mathbb{P}$ are more general than the
sets of Generalized Lakshmibai-Seshadri paths and they still have their nice properties as we show in the
following lemmata. Recall (see section \ref{weakchain}) what is meant by an integral path. We alter the
definition of a monotone path (section \ref{monotone}) by requiring
$h_i^{\pi}$ to be increasing and not necessarily strictly increasing in $[f_i^+(\pi), f_i^-(\pi)]$.\\

\begin{lemma}\label{integralmonotone} A path $\pi\in\ \Pithetahat$ is integral and monotone.
\end{lemma}
\begin{proof} Let $\pi\in\ \widehat{\mathbb{P}}_{\lambda},\, \pi'\in\ \widehat{\mathbb{P}}_{\mu},\, s,\, s'\in\ [0, 1]$ be as in section \ref{gluingpair} and assume that $\pitheta\in\
\Pithetahat$.

Set $h_i :=h_i^{\pitheta}$ and $m_i :=m_i^{\pitheta}$. Since the path $\pitheta$ is piecewise linear, a local
minimum of the function $h_i$ is attained at some $a_x,\, 0\leq x\leq k-1$ or $b_y,\, 1\leq y\leq \ell$ or at
$s, s'$. If a local minimum of $h_i$ is attained at $t\leq a_{k-1}$ or at $t\geq b_1$ then this number is an
integer by lemma \ref{integral}, since $\pi, \,\pi'$ are (weak) Generalized Lakshmibai-Seshadri paths and by the imposed condition that $\pitheta(1)\in P$.

It remains to examine the case where $\min \{h_i(t)|t\in\ [0, 1]\}=h_i(s)=h_i(s')$. This will mean that
$\alpha_i^{\vee}(\lambda_k)\leq 0$ and $\alpha_i^{\vee}(\mu_1)\geq 0$. If one of these numbers is zero,
then $h_i(s) = h_i(a_r)$ for some $r,\, 1\leq r\leq k-1$ or $h_i(s)=b_{r'}$ for some $r',\, 1\leq r'\leq
\ell$ and is an integer.

Assume then that $\alpha_i^{\vee}(\lambda_k)<0$ and $\alpha_i^{\vee}(\mu_1)>0$. Since we have
$\mu_1=\tau\mu$, lemma \ref{twodominant} gives that
$\alpha_i^{\vee}(\tau\lambda)=\alpha_i^{\vee}(\overline{\tau}\lambda)\geq 0$. By assumption, there exists an $s$-chain for
the pair $(\lambda_k, \overline{\tau}\lambda)$
$$\overline{\tau}\lambda:=\nu_t\stackrel{\beta_t}{\leftarrow}\nu_{t-1}\cdots \nu_1\stackrel{\beta_1}{\leftarrow}\nu_0=:\lambda_k,$$
and so by lemma \ref{simpleroot} one has that $\alpha_i=\beta_m$, for some $m$ with $1\leq m\leq t$. Then lemma
\ref{weakchain} gives $s\alpha_i^{\vee}(\lambda_k)=s\beta_m^{\vee}(\lambda_k)\in\ \mathbb{Z}$. Yet
$\pitheta(s)=\sum\limits_{j=1}^ka_j(\lambda_j-\lambda_{j+1})+s\lambda_k$ and since $\pi_1$ is a
GLS path, remark \ref{remarks} gives $a_j(\lambda_j-\lambda_{j+1})\in\ Q$ for all $j$ with $1\leq j\leq k$. Hence
$h_i(s)\in\ \mathbb{Z}$. We conclude that the path $\pitheta$ is integral.

Now if $f_i^-(\pi)<s$ or $f_i^+(\pi)>s'$ monotonicity follows by lemma \ref{monotone}. In the case where
$f_i^+(\pi)\leq s\leq s'\leq f_i^-(\pi)$ the path is monotone in the weaker sense (since $h_i(s)=h_i(s')$).
\end{proof}

\subsubsection{}
\begin{lemma}\label{Astable} Let $\lambda\in\
P^+$. The set $\Pithetahat$ is stable under the action of $f_i,\,
i\in\ I$.
\end{lemma}
\begin{proof} Let $\pitheta\in\ \Pithetahat$ and write it as in (\ref{pitheta}). We will show that if $f_i(\pitheta)\neq 0$ then $f_i(\pitheta)\in\ \Pithetahat$. First, note that since $f_i(\pitheta)(1)=\pitheta(1)-\alpha_i$ and $\pitheta(1)\in P$, we obtain that $f_i(\pitheta)(1)\in P$.

In this proof we set $f_+^i :=f_+^i(\pitheta),\,f_-^i:=f_-^i(\pitheta)$ and $h_i:=h_i^{\pitheta}$. If $i\in\
I^{re}$ and $f_-^i<s$ or $f_+^i>s'$, then the first joining condition follows trivially
and the second one trivially in the first case and by the remark in section \ref{gluingpair} in the second
case. Then $f_i(\pitheta)=\pi_1*\theta_s^{s'}*\pi_2'$ with $\pi_1\in \widehat{\mathbb{P}}_{\lambda}$ and $\pi_2'\in \widehat{\mathbb P}_{\mu}$ as in \cite[Proposition 5.6]{L}. On the other hand, if $i\in\ I^{im}$ and since $h_i$ is increasing in $[0, s]$ and $[s', 1 ]$ and is constant in $[s, s']$, hence increasing in $[0, 1]$, then $f_+^i=0$ or $s'$. If $f_+^i=0$ and $f_-^i<s$, then as in the case of real indices discussed above, the joining conditions follow trivially. Hence the only cases which need
to be checked are the following :

(1) Suppose that $f_-^i = s$; then
$$f_i(\pitheta)=\\
(\lambda_1,\dots, , \lambda_{t-1}, r_i\lambda_t, \dots r_i\lambda_k, 0, \mu_1,\dots, \mu_{\ell} ; 0, a_1, \dots, a_{k-1},
s, s', b_1, \dots, b_{\ell-1}, 1),
$$
with $t=1$, if $i\in\ I^{im}$.

We will show that $f_i(\pitheta)$ is the join of the paths
\begin{equation}\label{pi1}\pi_1 = (\lambda_1,\dots, , \lambda_{t-1}, r_i\lambda_t, \dots r_i\lambda_k; 0, a_1, \dots, a_{k-1}, 1)\in \widehat{\mathbb P}_{\lambda}
\end{equation}
(where $t=1$ if $i\in\ I^{im}$) and $\pi'\in \widehat{\mathbb P}_{\mu}$.

We first show that $\pi_1$ is in $\widehat{\mathbb P}_{\lambda}$.
Indeed, the existence of an $a_{t-1}$-chain for $(\lambda_{t-1},
r_i\lambda_t)$ and of $a_n$-chains for the pairs $(r_i\lambda_n,
r_i\lambda_{n+1})$ for all $n$ with $1\leq n\leq k-1$ follows as
in Propositions \ref{actionfre},
\ref{actionfim}. Note that $\pi_1$ is not
necessarily in $\mathbb{P}_{\lambda}$, since there might not exist
an $1$-chain for $(r_i\lambda_k, \lambda)$ for imaginary $i$.

It remains to show that $\pi_1$ and $\pi'$ can be properly joined across $[s, s']$.  Since $h_i(s)=m_i+1\in\ \mathbb{Z}$, by Lemma \ref{integralvalues} one has that
$s\alpha_i^{\vee}(\lambda_k)\in\ \mathbb{N}^+$ and as in Proposition \ref{actionfim} this number equals $1$ if $i\in I^{im}$, so there exists an $s$-chain for
$(r_i\lambda_k, \lambda_k)$. Combined with the given chain for $(\lambda_k, \overline \tau\lambda)$, we obtain an $s$-chain for $(r_i\lambda_k, \overline \tau\lambda)$. Since here $\tau$ is unchanged, the second joining condition immediately follows from the second condition on the starting path.\\

(2) Suppose now that $f^i_+<s$ and $f_-^i>s'$ and say $b_{m-1}<f_-^i\leq b_m$, with $1\leq m\leq \ell$. Then
\begin{align*}\label{pi} f_i(\pitheta) = (\lambda_1,\dots, , \lambda_{t-1}, r_i\lambda_t, \dots r_i\lambda_k, 0, r_i\mu_1, \dots, r_i\mu_m, \mu_m, \dots, \mu_{\ell}
;\\
0, a_1, \dots, a_{k-1}, s, s', b_1, \dots, b_{m-1}, f_-^i, b_m, \dots, b_{\ell-1}, 1),
\end{align*}
with $t=1$, if $i\in\ I^{im}$. One has that $\alpha_i^{\vee}(\lambda_k)> 0$ and $\alpha_i^{\vee}(\mu_1)> 0$.

We will show that $f_i(\pitheta)$ is the join of the paths $\pi_1$ given in (\ref{pi1}) (with $t=1$ if $i\in I^{im}$)
and
\begin{equation} \label{pi2}\pi_2' = (r_i\mu_1, \dots, r_i\mu_m, \mu_m, \dots, \mu_{\ell} ;
0, b_1, \dots, b_{m-1}, f_-^i, b_m,
\dots, b_{\ell-1}, 1)\in \widehat{\mathbb P}_{\mu}.
\end{equation}

Suppose first that $i\in\ I^{im}$. As in the previous case, we have that the path $\pi_1$ is in $\widehat{\mathbb P}_{\lambda}$. We will show that $\pi_2'\in \widehat{\mathbb P}_{\mu}$.

From the constraint imposed by the identity
$h_i(f_-^i)=\alpha_i^\vee(\pitheta(f_-^i))=1$, we
obtain as in say Proposition \ref{actionfim} that
$\alpha_i^\vee(\lambda_1)=\ldots
=\alpha_i^\vee(\lambda_k)$ and $\alpha_i^\vee(\mu_1)=\ldots=\alpha_i^\vee(\mu_m)\geq
\alpha_i^\vee(\mu_{m+1})$. Then this identity becomes
\begin{equation} \label{star} s\alpha_i^\vee(\lambda_k)+(f_-^i-s')\alpha_i^\vee(\mu_m)=1.
\end{equation}

It follows by Lemma \ref{middlechainim} that there exist $b_q$-chains for $(r_i\mu_q, r_i\mu_{q+1})$ for all $q$, with $1\le q\le m-1$. We need to show that there exists a weak $f_-^i$-chain for $(r_i\mu_m, \mu_m)$, i.e. that $f_-^i\alpha_i^\vee(\mu_m)$ is a positive
integer.

Since $h_i(f_-^i)=1$ by definition and $h_i(1)\in \mathbb N$, $\pitheta(1)$ being an integral weight, one has that $h_i(1)-h_i(f_-^i)\in \mathbb N$ that is
$\alpha_i^\vee(\pitheta(1))-\alpha_i^\vee(\pitheta(f_-^i)) \in \mathbb N$,
whilst it also equals
$$-f_-^i\alpha_i^\vee(\mu_m)+
\sum_{j=m}^{\ell-1}b_j\alpha_i^\vee(\mu_j-\mu_{j+1})+\mu_{\ell}.$$

Now $b_j(\mu_j-\mu_{j+1})\in -Q^+$ for all $j$, with $m\le j\le
\ell-1$, by the condition that $\pi'\in \widehat{\mathbb P}_\mu$
and equation (\ref{goodweight}) and so
$\sum\limits_{j=m}^{\ell-1}b_j\alpha_i^\vee(\mu_j-\mu_{j+1})\in
\mathbb N$. On the other hand, the weight $\mu_{\ell}$ is
integral. So necessarily $f_-^i\alpha_i^\vee(\mu_m)\in \mathbb N$
and $\pi_2'\in \widehat{\mathbb P}_{\mu}$ (recall that
$\alpha_i^{\vee}(\mu_m)>0$ by monotonicity).

We prove next that $\pi_1$ and $\pi_2'$ satisfy the joining conditions. Recall that there exists an $s$-chain for $(\lambda_k, \overline \tau\lambda)$. By equation (\ref{goodweight}) we obtain $s(\lambda_k-\overline \tau\lambda)\in\ -Q^+$ and so
$s\alpha_i^{\vee}(\lambda_k-\overline\tau \lambda)\in\ \mathbb{N}$, by Lemma \ref{improperties} (2). Yet
$s\alpha_i^{\vee}(\overline \tau\lambda)\geq 0$, by Lemma \ref{Tlambda}, whilst
$s\alpha_i^{\vee}(\lambda_k)<1$, since $f_-^i>s$, that is :

\begin{equation}\label{easy} 0\leq s\alpha_i^{\vee}(\overline \tau\lambda)\leq s\alpha_i^{\vee}(\lambda_k)<1.
\end{equation}

This forces $\alpha_i^{\vee}(\lambda_k)=\alpha_i^{\vee}(\overline{\tau}\lambda)$. By Lemma \ref{middlechainim} there exists an $s$-chain for $(r_i\lambda_k, r_i\overline{\tau}\lambda)$ and hence for
$(r_i\lambda_k, \overline{r_i\tau}\lambda)$.

Recall that $\mu_{\ell}\stackrel{\beta_{\ell-1}}{\leftarrow}\mu_{\ell-1}\cdots
\mu_2\stackrel{\beta_1}{\leftarrow}\mu_1$. Observe from Proposition \ref{actionfim} that $f_i (\pitheta)$
being specified as above means that $\alpha_i^{\vee}(\beta_t)=0$ for all $t$, with $1\leq t\leq m$ and
$\mu_{\ell}\stackrel{\beta_{\ell-1}}{\leftarrow}\mu_{\ell-1}\cdots
\mu_m\stackrel{\alpha_i}{\leftarrow}r_i\mu_m \stackrel{\beta_{m-1}}{\leftarrow}r_i\mu_{m-1}\cdots
r_i\mu_2\stackrel{\beta_1}{\leftarrow}r_i\mu_1$. Then $r_i\tau = r_{\beta_1}\cdots
r_{\beta_m}r_ir_{\beta_{m+1}}\cdots r_{\beta_{\ell-1}}$ and the second joining condition reduces to
$s\alpha_i^{\vee}(\overline \tau\lambda)<1$, verified in (\ref{easy}).

Suppose that $i\in\ I^{re}$. We show as in the first case that $\pi_1\in \widehat{\mathbb P}_{\lambda}$. On the other hand, $\pi_2'\in \widehat{\mathbb P}_{\mu}$; the existence of $b_q$-chains for $(r_i\mu_q, r_i\mu_{q+1})$ for all $q$ with $1\le q\le m-1$ follows by Corollary \ref{middlechain} and the fact that $f_-^i\alpha_i^{\vee}(\mu_m)$ is an integer follows by a similar computation (the only difference being that $h_i(f_-^i)$ is not necessarily $1$).

We finally prove that the joining conditions hold for $\pi_1,\, \pi_2'$. Since $\alpha_i^{\vee}(\tau\mu)>0$, we obtain by Lemma \ref{twodominant} that
$\alpha_i^{\vee}(\overline{\tau}\lambda)\geq 0$. However, the latter inequality is strict. Indeed, notice that since $\pi\in \widehat{\mathbb P}_\lambda$ and there exists an $s$-chain for $(\lambda_k, \overline\tau\lambda)$, the path $\hat\pi:=(\lambda_1, \dots, \lambda_k, \overline\tau\lambda; 0, a_1, \dots, a_{k-1}, s, 1)$ is also a path in $\widehat{\mathbb P}_\lambda$. If $\alpha_i^{\vee}(\overline\tau\lambda)=0$, and since $s\alpha_i^{\vee}(\lambda_k)\notin \mathbb Z$, one would get $\hat\pi(1)\notin P$, which is a contradiction. So we must have $\alpha_i^{\vee}(\overline\tau\lambda)>0$. Again $\alpha_i^{\vee}(\lambda_k)>0$ and so there exists an
$s$-chain for $(r_i\lambda_k, \overline{r_i\tau}\lambda)$ by Corollary \ref{middlechain} (1). The second
joining condition in this case follows by the remark of Section \ref{gluingpair}.

(3) Finally suppose that $f_+^i=s'<f_-^i$. Then
\begin{align*}f_i (\pitheta) = (\lambda_1,\dots, \lambda_k, 0, r_i\mu_1, \dots, r_i\mu_m, \mu_m, \dots, \mu_{\ell};\\
0, a_1,\dots, a_{k-1}, s, s', b_1,\dots, b_{m-1}, f_-^i, b_m,\dots, b_{\ell-1}, 1).
\end{align*}
We will show that $f_i(\pitheta)$ is the join of $\pi\in \widehat{\mathbb P}_{\lambda}$ and $\pi_2'\in \widehat{\mathbb P}_{\mu}$, where $\pi_2'$ is given by (\ref{pi2}). The fact that $\pi_2'\in \widehat{\mathbb P}_{\mu}$ follows exactly as in the previous case. Notice that equation (\ref{star}) becomes ($i\in I^{im}$)
\begin{equation} \label{twostar} (f_-^i-s')\alpha_i^{\vee}(\mu_m)=1.
\end{equation}
We will prove that the joining conditions hold for $\pi,\, \pi_2'$.

Suppose that $i\in\ I^{im}$. Then, we have that $\alpha_i^{\vee}(\lambda_k)=0$ and so
$\alpha_i^{\vee}(\overline{\tau}\lambda)=0$ and $\alpha_i^{\vee}(\mu_1)>0$. But then
$\overline{r_i\tau}\lambda=\overline{\tau}\lambda$ so there exists an $s$-chain for $(\lambda_k,
\overline{r_i\tau}\lambda)$. The second joining condition follows by the vanishing of
$\alpha_i^{\vee}(\overline \tau \lambda)$.

Suppose now that $i\in\ I^{re}$. We have that $\alpha_i^{\vee}(\lambda_k)\leq 0$ and
$\alpha_i^{\vee}(\mu_1)=\alpha_i^{\vee}(\tau \mu)>0$. Lemma \ref{twodominant} gives
$\alpha_i^{\vee}(\overline{\tau}\lambda)\geq 0$; if the latter inequality is strict, by Lemma \ref{leftchain} we obtain an $s$-chain for
$(\lambda_k, \overline{r_i\tau}\lambda)$, hence the first joining condition holds. If $\alpha_i^{\vee}(\overline{\tau}\lambda)=0$, then $\overline{r_i\tau\lambda}=\overline{\tau\lambda}$ and there is nothing to prove. Finally, the second one follows by
the remark of Section \ref{gluingpair}. Then $f_i (\pitheta) = \pi*\theta_s^{s'}*\pi_2'$, where $\pi_2'$ is
as in (\ref{pi2}). The assertion follows.

\end{proof}

\subsubsection{} Our aim now is to prove that $\mathcal{A}(\pi_{\lambda}\otimes \pi_{\mu})$ is a highest
weight crystal generated by $\pi_{\lambda}\otimes \pi_{\mu}$ over $\mathcal{F}$ (see lemma \ref{unique2}).
The following two lemmata are preliminary results for this purpose.\\

\begin{lemma}\label{twochains} Let $\nu : = \nu_s\stackrel{\beta_s}{\leftarrow} \nu_{s-1}\stackrel{\beta_{s-1}}{\leftarrow}\dots
\stackrel{\beta_2}{\leftarrow}\nu_1\stackrel{\beta_1}{\leftarrow}\nu_0=:\mu$ be an $a$-chain for $(\mu,
\nu)$, such that $\beta_l=w\alpha_i$, where $i\in\ I^{im}$ for some $l$, with $1\leq l\leq s$. Suppose
further that $a\alpha_i^{\vee}(\mu)=1-a_{ii}$. Then $\beta_l=\alpha_i$, $\mu=r_i\mu'$ and there exists an
$a$-chain for $(\mu', \nu)$.
\end{lemma}

\begin{proof}

Suppose that $\beta_l=w\alpha_i\neq \alpha_i$. Since $-\alpha_i$ is dominant, we obtain
$\beta_l=\alpha_i+\beta\in\ \alpha_i+\mathbb{N}\Pi_{re}$ and $\alpha_i^{\vee}(\beta)\leq -1$. By the
hypothesis,
\begin{equation}\label{first}a\alpha_i^{\vee}(\mu)=1-a_{ii}.
\end{equation}
On the other hand,\\
\begin{equation}\label{second}a\alpha_i^{\vee}(\mu)= a\alpha_i^{\vee}(\nu_l)- a\beta_l^{\vee}(\nu_l)\alpha_i^{\vee}(\beta_l)-\sum\limits_{q=1}^{l-1}a\beta_q^{\vee}(\nu_q)\alpha_i^{\vee}(\beta_q).
\end{equation}
Now $a\beta_l^{\vee}(\nu_l) = 1$ by \ref{defgls} (b), whereas
$\alpha_i^{\vee}(\beta_l)=\alpha_i^{\vee}(\alpha_i+\beta)\leq a_{ii}-1$. Then (\ref{first}) and
(\ref{second}) give that :

\begin{equation}\label{third} 0 \geq \underbrace{a\alpha_i^{\vee}(\nu_l)}_{\geq \,0} - \sum\limits_{q=1}^{l-1}\underbrace{a\beta_q^{\vee}(\nu_q)\alpha_i^{\vee}(\beta_q)}_{\leq \,0},
\end{equation}
which means that all the summands in (\ref{third}) are equal to zero so $\alpha_i^{\vee}(\nu_l)=0$ and
$\alpha_i^{\vee}(\beta_l)=a_{ii}-1$. This on one hand means that $\alpha_i^{\vee}(\beta_q)=0$ for all $q$,
with $1\leq q\leq l-1$ and so $r_i$ commutes with all $r_{\beta_q}$ with $1\leq q\leq l-1$. On the other hand
we can write $\beta_l = w\alpha_i=w_1r_k\alpha_i=w_1(\alpha_i-\alpha_k^{\vee}(\alpha_i)\alpha_k)$, with
$w_1\alpha_k=\alpha_k+\beta_1$ and moreover $w_1\alpha_i=\alpha_i+\beta_2$, with $\beta_1, \beta_2\in\
\mathbb{N}\Pi_{re}$. Yet $\alpha_i^{\vee}(\beta_l)=a_{ii}-1$, which forces $\alpha_i^{\vee}(\beta_1)=0$,
$\alpha_i^{\vee}(\beta_2)=0$ and $\alpha_i^{\vee}(\alpha_k)\alpha_k^{\vee}(\alpha_i)=1$. The second condition
forces $w_1\alpha_i=\alpha_i$ and $\alpha_i^{\vee}(\alpha_k)=\alpha_k^{\vee}(\alpha_i)=-1$. Then
$\gamma:=w_1\alpha_k$ is such that $\beta_l = w\alpha_i=r_{\gamma}\alpha_i$. Note that
$\gamma^{\vee}(\alpha_i)=\alpha_k^{\vee}(\alpha_i)=-1$ and
$\alpha_i^{\vee}(\gamma)=\alpha_i^{\vee}(\alpha_k)=-1$. Also $\beta_l=w_1r_k\alpha_i=\alpha_i+\gamma$. Then

\begin{equation}1 = a\beta_l^{\vee}(\nu_l) = a\alpha_i^{\vee}(r_{\gamma}\nu_l) =-a\gamma^{\vee}(\nu_l)\alpha_i^{\vee}(\gamma) = a\gamma^{\vee}(\nu_l),
\end{equation}
since $\alpha_i^{\vee}(\nu_l)=0$ and $\alpha_i^{\vee}(\gamma)=a_{ik}=-1$.  Again,

\begin{equation}a(r_ir_{\gamma}\nu_l-\nu_l)=-a\gamma^{\vee}(\nu_l)r_i\gamma =
-(\gamma +\alpha_i)=-\beta_l,
\end{equation}
whilst

$$a(r_{\beta_l}\nu_l-\nu_l)=-\beta_l.$$
Then $\nu_{l-1}=r_{\beta_l}\nu_l=r_ir_{\gamma}\nu_l>r_{\gamma}\nu_l>\nu_l$ and so $\dist\, (\nu_{l-1},
\nu_l)\geq 2$, which contradicts our hypothesis. Then necessarily $\beta_l=\alpha_i$ and $r_i$ commutes with
$r_{\beta_q}$ for all $q$, with $1\leq q<l$. It then follows that $\nu : =
\nu_s\stackrel{\beta_s}{\leftarrow}
\nu_{s-1}\stackrel{\beta_{s-1}}{\leftarrow}\cdots\nu_{l+1}\stackrel{\beta_{l+1}}{\leftarrow}\nu'_l\stackrel{\beta_{l-1}}{\leftarrow}\nu'_{l-1}\cdots
\stackrel{\beta_1}{\leftarrow}\nu'_1\stackrel{\alpha_i}{\leftarrow}\nu_0=:\mu$ is an $a$-chain, where
$r_i\nu'_q=\nu_q$ for all $q$, with $1\leq q\leq l$ and so there exists an $a$-chain for $(\nu_1', \nu)$.
\end{proof}

\subsubsection{}
\begin{lemma}\label{twodom}
Let $\lambda,\, \mu\in P^+$ and suppose that $\mu : = \mu_\ell\stackrel{\beta_\ell}{\leftarrow}\mu_{\ell-1}\stackrel{\beta_{\ell-1}}{\leftarrow}\cdots \mu_1\stackrel{\beta_1}{\leftarrow}\mu_0$. Set $\tau=r_{\beta_1}r_{\beta_2}\cdots r_{\beta_\ell}$. If $\alpha_i^\vee(\tau\lambda)<0$ for some $\alpha_i\in \Pi_{re}$, then also $\alpha_i^\vee(\tau\mu)<0$.
\end{lemma}
\begin{proof}
By lemma \ref{twodominant}, one has that $\alpha_i^\vee(\tau\mu)\le 0$; we need to show that this number is strictly negative. As in Section \ref{domred}, we may write $\tau$ in its dominant reduced expression $\tau=w\tau^*$, where $w\in W$ and $\tau^*\in T$ is such that for all $\nu\in P^+,\, \tau^*\nu\in P^+$.

Set $\mu^*:=\tau^*\mu$; if $\alpha_j^\vee(w\mu^*)<0$ for some $\alpha_j\in \Pi_{re}$, one has that $\ell(r_jw)<\ell(w)$ and so $w=r_jw'$ with $w'\in W$ and $\ell(w)=\ell(w')+1$.

We write $w$ as $w=r_{i_1}r_{i_2}\cdots r_{i_n}\tilde w$ where $\tilde w\mu^*=\mu^*$, $\ell(w)=\ell(\tilde w)+n$ and $\alpha_{i_t}^\vee(r_{i_t}\cdots r_{i_n}\mu^*)<0$, for all $t$, with $1\le t\le n$. We set $w':=r_{i_1}r_{i_2}\cdots r_{i_n}$. We also set $w_t'=r_{i_t}\cdots r_{i_n}$ and $\hat w_t=w'w_t'^{-1}$, for all $t$, with $1\le t\le n$.

Let $S(w^{-1})$ denote the set of positive roots which become negative if we apply $w^{-1}$; one has $S(w'^{-1})=\{\gamma_t:=\hat w_t\alpha_{i_t}\,|\,1\le t\le n\}\subset S(w^{-1})$. Then $w'^{-1}\gamma_t=w_t'^{-1}\alpha_{i_t}$ and so
$$\gamma_t^\vee(w\mu^*)=\gamma_t^\vee(w'\mu^*)=(w'^{-1}\gamma_t)^\vee(\mu^*)=(w_t'^{-1}\alpha_{i_t})^\vee(\mu^*)=\alpha_{i_t}^\vee(w_t'\mu^*)<0$$
and in particular non-zero. Since $\alpha_i^\vee(\tau\lambda)<0$, we have that $\alpha_i\in S(w^{-1})$ but by the above $\alpha_i\notin S(w'^{-1})$.

On the other hand, $\alpha_{i_1}^\vee(\tau\mu)<0$ and $\alpha_{i_1}^\vee(\mu)\ge 0$; by Lemma \ref{simpleroot} it follows that $\alpha_{i_1}$ is equal to $\beta_k$ for some $k$ with $1\le k\le \ell$.

Then one has that
\begin{equation}\alpha_{i_1}^\vee(\mu_{k})>0, \alpha_{i_1}^\vee(\mu_{k-1})<0.
\end{equation}

Let us assume $k$ is minimal having the above property (which is equivalent to $\beta_k=\alpha_{i_1}$).

Set $\beta_s':=r_{i_1}\beta_s$, for all $s=1, 2,\dots, k-1$; they are positive roots by the minimality of $k$. Then $(\beta_s')^\vee(r_{\beta_{s+1}'}\cdots r_{\beta_{k-1}'}\mu_{k})=\beta_s^\vee(\mu_{s})>0$ and $\tau\mu=r_{i_1}r_{\beta_1'}\cdots r_{\beta_{k-1}'}r_{\beta_{k+1}}\cdots r_{\beta_\ell}\mu=r_{i_1}r_{i_2}\cdots r_{i_n}\mu^*$ and so we may cancel $r_{i_1}$ and get $$r_{\beta_1'}\cdots r_{\beta_{k-1}'}r_{\beta_{k+1}}\cdots r_{\beta_\ell}\mu=r_{i_2}\cdots r_{i_n}\mu^*.$$

We continue in the same way and obtain an expression:
$$\mu^*=\tilde w\mu^*=r_{\tilde \beta_1}\cdots r_{\tilde \beta_{t}}\mu,$$
where $\tilde \beta_j\in \Delta^+\cap W\Pi$ and such that $\tilde \beta_j^\vee(r_{\tilde \beta_{j+1}}\cdots r_{\tilde \beta_{t}}\mu)>0$ and $\tilde w\in \Stab_W(\mu^*)$.

Thus we are reduced to the situation where $\tau \mu$ is dominant and $\tau=w\tau^*$ with $w \in \Stab_W(\mu),\tau^* \in T$. Since the stabilizer of a dominant weight is the product of simple reflections each of which stabilizes it, it is enough to obtain a contradiction with $w=r_i$.

Notice that in the above we may assume $\ell\geq 1$. By lemma \ref{chainfordom} and since $\tau\mu\in P^+$ we necessarily have that $\beta_1\in \Pi_{im}$. Then by Lemma \ref{dominance} $\alpha_i^{\vee}(\tau\lambda)<0$ implies that $r_i$ and $r_{\beta_1}$ commute, which means that $\alpha_i^\vee(\beta_1)=0$.

Then $r_i$ still stabilizes $\tau_2\mu:=r_{\beta_2}\cdots r_{\beta_\ell}\mu$ and $\alpha^\vee_i(\tau_2 \lambda)<0$.  Thus we can replace $\tau$ by $\tau_2$.  Again we can write $\tau_2=w_2\tau_2^*$, such that $\tau_2^*\lambda,\, \tau_2^*\mu$ are dominant and as before we may further cancel $r_{\beta_t}'s$, so that we are reduced to the situation where $\tau_2\mu$ is dominant and $w_2\in \Stab_W(\tau_2^*\mu)$. This procedure will have to stop, hence we get a contradiction.

%

\end{proof}

\subsubsection{}
\begin{lemma}\label{unique2} Let $\lambda,\, \mu\in\ P^+$. A path
in $\mathcal{A}(\pi_{\lambda}\otimes \pi_{\mu})$ is integral and monotone and the only path killed by the $e_i,\,
i\in\ I$ is $\pi_{\lambda}\otimes \pi_{\mu}$. In particular, $\mathcal{A}(\pi_{\lambda}\otimes \pi_{\mu}) =
\mathcal{F}(\pi_{\lambda}\otimes \pi_{\mu})$.
\end{lemma}
\begin{proof} We can write $\pi_{\lambda}\otimes \pi_{\mu} = \pi_{2\lambda}*\theta_{1/2}^{1/2}*\pi_{2\mu}\in\ \Pithetahatdouble$, since the joining conditions become trivial.

By lemma \ref{Astable} the set $\Pithetahatdouble$ is stable under
the action of $f_i,\, i\in\ I$, hence $\mathcal F(\pi_\lambda\otimes \pi_\mu)\subset \Pithetahatdouble$. Note that $\mathcal F(\pi_\lambda\otimes \pi_\mu)\subset \Plambda\otimes \mathbb P_\mu\subset \Plambdahat\otimes \Pmuhat$. Set $\textbf{P}:= \Plambda\otimes \mathbb P_\mu\cap \Pithetahatdouble$; by the above, $\mathcal F(\pi_\lambda\otimes \pi_\mu)\subset \textbf{P}$. Since both sets $\Pithetahatdouble$ and $\Plambda\otimes \mathbb P_{\mu}$ are $\mathcal F$-stable, $\textbf P$ is also $\mathcal F$-stable. We will show that as a subset of $\Plambda\otimes \mathbb P_\mu$, the set $\textbf{P}$ is
stable under the action of $e_i,\, i\in\ I$ and that the only path in
$\textbf{P}$ killed by all the $e_i,\, i\in\ I$ is the path
$\pi_{\lambda}\otimes \pi_{\mu}$. This will give $\textbf{P}= \mathcal{A}(\pi_{\lambda}\otimes \pi_{\mu}) = \mathcal{F}(\pi_{\lambda}\otimes \pi_{\mu})$. Finally the
integrality and monotonicity of the paths in $\mathcal{F}(\pi_{\lambda}\otimes \pi_{\mu})$ will follow by
lemma \ref{integralmonotone}.

We first show that $\textbf{P}$ is $e_i$ stable in the above sense for all $i\in\ I$. Let $\tilde \pi=\pi*\theta_{1/2}^{1/2}*\pi'$ be as in paragraph \ref{gluingpair}. We
will show that if $e_i\tilde \pi\neq 0$,  then $e_i\tilde \pi\in\
\Pithetahatdouble$. Since $\tilde \pi\in \textbf P$, then $\pi\in \mathbb P_{2\lambda}$ and $\pi'\in \mathbb P_{2\mu}$. On the other hand $\tilde \pi\in \Plambda\otimes \mathbb P_{\mu}$, hence $\tilde \pi=\pi_1\otimes \pi_2$ with $\pi_1\in \Plambda,\, \pi_2\in \mathbb P_\mu$

Since $\tilde \pi(1/2)\in\ P$, we can only have either $e_i(\pi_1\otimes \pi_2) = (e_i\pi_1)\otimes \pi_2$ or
$e_i(\pi_1\otimes \pi_2) = \pi_1\otimes (e_i\pi_2)$ the choice depending on the crystal tensor product rules
(section \ref{crystaltensorproduct}). It remains to check the joining condition for the new paths.

(1) Suppose first that $e_i(\pi_1\otimes \pi_2) = (e_i\pi_1)\otimes \pi_2\neq 0$.

If $e_+^i(\pi_1)<1$ (and so $e_+^i(\tilde \pi)<1/2$), $\lambda_k, \mu_1$ are unchanged, so there is nothing
to check.\\
Suppose then that $e_+^i(\pi_1)=1$, equivalently $e_+^i(\tilde \pi)=1/2$.

Since $\mu_1$ is unchanged we only need to show that there exists an $1/2$-chain for $(r_i2\lambda_k,
\overline \tau2\lambda)$. This is equivalent to the existence of an $1$-chain for $(r_i\lambda_k, \overline
\tau\lambda)$.\\
Let $i\in\ I^{re}$.

By the definition of $e_+^i(\tilde \pi)$, one has that $\alpha_i^{\vee}(\lambda_k)<0$ and
$\alpha_i^{\vee}(\tau \mu)\geq 0$. Then, by lemma \ref{twodom} we have that $\alpha_i^{\vee}(\overline \tau\lambda)=\alpha_i^{\vee}(\tau \lambda)\geq 0$ and so by lemma
\ref{leftchain} there exists an $1$-chain for $(r_i\lambda_k, \overline{\tau}\lambda)$, since there exists
one for $(\lambda_k, \overline \tau\lambda)$.
\\
Let now $i\in\ I^{im}$.

This means that there exists a path $\pi_1'\in\ \Plambda$ with $\pi_1'=(\lambda_1', \dots, \lambda_k';0,
a_1,\dots, a_{k-1}, 1)$ such that $f_i\pi_1'=\pi_1$, which in turn gives $\lambda_t=r_i\lambda_t'$ for all
$t$, with $1\leq t\leq k$.  Let $$\overline
\tau\lambda:=\nu_s\stackrel{\gamma_s}{\leftarrow}\nu_{s-1}\stackrel{\gamma_{s-1}}{\leftarrow}\cdots
\stackrel{\gamma_2}{\leftarrow}\nu_1\stackrel{\gamma_1}{\leftarrow}\nu_0=:\lambda_k,$$ with $\gamma_j\in\
W\Pi\cap \Delta^+$ for all $j$, with $1\leq j\leq s$, be an $1$-chain for $(\lambda_k, \overline
\tau\lambda)$. We need to show that there also exists an $1$-chain for $(\lambda_k', \overline \tau\lambda)$.
By the second joining condition, if $\overline \tau=r_{\beta_1}\cdots r_{\beta_n}$, then
$\beta_t^{\vee}(r_{\beta_{t+1}}\cdots r_{\beta_n}\lambda)<1$ for all $t$ such that $\beta_t\in\ \Delta_{im}$.
This forces $\overline \tau\in\ W$ which in turn implies that $\gamma_l=w\alpha_i$ for some $l$, with $1\leq
l\leq s$ and some $w\in\ W$. By assumption that $e_+^i(\pi_1)=1$, we obtain
$\alpha_i^{\vee}(\lambda_k)=1-a_{ii}$. Then the assertion follows by lemma \ref{twochains}.

(2) Suppose that $e_i(\pi_1\otimes \pi_2) = \pi_1\otimes (e_i\pi_2)\neq 0$.

If $e_-^i(\pi_1)>0$ (and so $e_-^i(\tilde \pi)>1/2$), then since the $\lambda_k, \mu_1$ are unchanged, the
first joining condition is trivial. The second one follows by the remark of section \ref{gluingpair}.
\\
Suppose then that $e_-^i(\pi_2)=0$ and so $e_-^i(\tilde \pi)=1/2$.\\
Suppose that $i\in\ I^{re}$.

We have that $\alpha_i^{\vee}(\lambda_k)\leq 0$ and $\alpha_i^{\vee}(\tau\mu)<0$, which implies that
$\alpha_i^{\vee}(\overline \tau\lambda)\leq 0$. If the latter is zero, the first joining condition trivially
follows. If not, then $\overline{r_i\tau}\lambda\stackrel{\alpha_i}{\leftarrow}\overline \tau\lambda$ and so
there exists an $1$-chain for $(\lambda_k, \overline{r_i\tau}\lambda)$. The second joining condition follows
by the remark of section \ref{gluingpair}.
\\
Suppose now that $i\in\ I^{im}$.

Then $\alpha_i^{\vee}(\lambda_k)=0$, and so since $\lambda_k\geq \overline\tau\lambda$, by lemma
\ref{positivityofhim} we obtain that $\alpha_i^{\vee}(\lambda_k)=\alpha_i^{\vee}(\overline \tau\lambda)=0$.
There exists a path $\pi_2'\in \mathbb P_\mu$, such that $\pi_2=f_i\pi_2'$. Then by Proposition \ref{actionfim} $\pi_2'$ will be of the form $$\pi_2'=(\mu_1'. \mu_2', \dots, \mu_t', \mu_{t+2}', \dots, \mu_\ell; 0, b_1, \dots, b_{t-1}, b_{t+1}, \dots, b_{\ell-1}, 1),$$
with $\mu_k=r_i\mu_k'$, for all $k$, with $1\le k\le t$, $f_i^-(\pi'_2)=b_t$ and $\mu_t'=\mu_{t+1}$. Again by Proposition \ref{actionfim} $r_ir_{\beta_k}=r_{\beta_k}r_i$ and $\tau=r_i\tau'$.
Since $0=\alpha_i^\vee(\overline \tau\lambda)=\alpha_i^\vee(\overline {r_i\tau'}\lambda)=(1-a_{ii})\alpha_i^\vee(\overline{\tau'}\lambda)$, it follows that the second joining condition holds and $\overline {\tau'}\lambda=\overline \tau\lambda$, hence there exists an $1$-chain for $(\lambda_k, \overline{\tau'}\lambda)$. We conclude that $\pi_1,
e_i\pi_2$ can be properly joined.

We finally show that $\pi_{\lambda}\otimes \pi_{\mu}$ is the only path in
$\textbf P$ killed by $e_i,\, i\in\ I$. For this we first
show that every $\pi\in\ \Plambda\otimes \mathbb{P}_{\mu}$, killed by the $e_i,\, i\in\ I$, takes the form
$\pi=\pi_{\lambda}\otimes \pi'$, with $\pi'\in\ \mathbb{P}_{\mu}$ and $\lambda+\pi'(1)\in\ P^+$.

Recall the tensor product crystal operations of section \ref{crystaltensorproduct}. Take $\pi=\pi_1\otimes
\pi_2\in\ \Plambda\otimes \mathbb{P}_{\mu}$ and assume that $e_i\pi=0$ for all $i\in\ I$. Let $i\in\ I^{re}$.
If $\varepsilon_i(\pi_2)>\varphi_i(\pi_1)$ one has that $e_i(\pi_1\otimes \pi_2)=\pi_1\otimes e_i\pi_2\neq 0$
by normality. So we must have
\begin{equation}\label{tp} \varphi_i(\pi_1)\geq \varepsilon_i(\pi_2),\quad \mbox{for
all} \quad i\in\ I^{re}.
\end{equation}
But then $e_i\pi=e_i\pi_1\otimes \pi_2$ and consequently, we must have $e_i\pi_1=0$, for all $i\in\ I^{re}$. Now take $i\in\ I^{im}$
and recall lemma \ref{fimgls}. One has that $\varphi_i(\pi_1)>-a_{ii}\Leftrightarrow \alpha_i^{\vee}(\wt
\pi_1)\geq 1-a_{ii}\Leftrightarrow e_i\pi=e_i\pi_1\otimes \pi_2$. On the other hand, if
$\alpha_i^{\vee}(\pi_1)<1-a_{ii}$, then $e_i\pi_1=0$ again by lemma \ref{fimgls}. In both cases $e_i\pi_1=0$.
We conclude that $e_i\pi=0$ for all $i\in\ I$, only if $e_i\pi_1=0$ for all $i\in\ I$ which forces
$\pi_1=\pi_{\lambda}$.
Notice also by (\ref{tp}) one has that
\begin{equation}\label{lambda}
\alpha_i^{\vee}(\lambda+\wt \pi_2)\geq 0,\quad \mbox{for all} \quad i\in\ I,
\end{equation}
that is $\lambda+\pi_2(1)\in\ P^+$. Set $J=\{i\in\ I\,|\, \alpha_i^{\vee}(\lambda)=0\}$ and assume that
$e_i(\pi_{\lambda}\otimes \pi)=0$ for all $i\in\ I$. Then $e_i\pi = 0$ for all $i\in\ J$. A path in
$\textbf P$ killed by all the $e_i,\, i\in\ I$ will then be
of the form $\tilde \pi = \pi_{2\lambda}*\theta_{1/2}^{1/2}*\pi$, where
$$\pi=(\mu_1, \dots, \mu_{\ell};0, b_1, \dots, b_{\ell-1}, 1)$$ is a GLS path of shape $\mu$ with $1/2<b_1$
and $\mu_1=\tau\mu$. Now the first joining condition forces $\tau\in\ \operatorname{Stab}_T(\lambda)$. If we
set $\mu_r=\tau_r\mu$ with $1\leq r\leq \ell$, then since $\mu_r>\mu_{r+1}$ we will have that $\tau_r\in\
\operatorname{Stab}_T(\lambda)$ for all $r$ with $1\leq r\leq \ell$. Take $i\in\ I\setminus J$. Then since
$\operatorname{Stab}_T(\lambda)=\langle r_i\,|\,i\in\ J\rangle$ by lemma \ref{stab}, we have that $\wt \pi\in\
\mu-\sum\limits_{j\in\ J}\mathbb{N}\alpha_j$ and so $e_i\pi=0$, for $i\in\ I\setminus J$ since the set of
weights of $\mathbb{P}_{\mu}$ lies in $\mu-Q^+$. Combined with our previous result, namely that $e_i\pi=0$
for all $i\in\ J$, this forces $\mu_1=\mu$ and the only path in
$\textbf P$ annihilated by all the $e_i$ is
$\pi_{2\lambda}*\theta_{1/2}^{1/2}*\pi_{2\mu} = \pi_{\lambda}\otimes \pi_{\mu}$.
\end{proof}

\subsection{The isomorphism theorem}\label{distofpix}
Recall the family $\{\pi^x,\, x\in\ [0, 1]\}$ of section \ref{pix}, which deforms the path
$\pi_{\lambda}\otimes \pi_{\mu}$ to $\pi_{\nu}$. We will show that $\mathcal{F}\pi^0\simeq \mathcal{F}\pi^1$
and then that $\pi_{\lambda}\otimes \pi_{\mu}$ and $\pi_{\nu}$ generate isomorphic crystals.

\subsubsection{}
By the construction of \ref{pix} and proposition
\ref{integralmonotone} it follows that $f\pi^x$ is integral and monotone for all $x\in\ [0,
1]$ and all $f\in\ \mathcal{F}$.\\

\begin{lemma}\label{pix2} Let $x,\, y\in\ [0, 1]$ and let $\pi^x,\, \pi^y$ be the paths defined in section \ref{pix}. Then $\mathcal{F}\pi^x\simeq \mathcal{F}\pi^y$.
\end{lemma}
\begin{proof} Let $J$ be a finite subset of $I$. By a direct computation, for all $i\in\ J$ we obtain:
\begin{eqnarray} d_J(\pi^x, \pi^y) = \max_{i\in\ J, t\in\ [0,
1]}|\alpha_i^{\vee}(\pi^x(t))-\alpha_i^{\vee}(\pi^y(t))| =\notag\\
=|x-y| \max_{i\in\ J, t\in\ [0, 1]}|\alpha_i^{\vee}((\pi_{\lambda}\otimes
\pi_{\mu})(t))-\alpha_i^{\vee}(\pi_{\lambda+\mu}(t))|=|x-y|d_J(\pi_{\lambda}\otimes \pi_{\mu}, \pi_{\nu}).
\end{eqnarray}
We reduce the distance of $x$ and $y$
so that $d_J(\pi^x, \pi^y)<(1/2c_J)^n$, which by lemma \ref{distance} implies that
$\mathcal{F}_J^n\pi^x\simeq \mathcal{F}_J^n\pi^y$. Since $n$ and $J$ are arbitrary the assertion follows.
\end{proof}

\subsubsection{} \label{iso}
Recall that $\lambda,\, \mu\in\ P^+$ and consider $\pi_{\lambda}\otimes \pi_{\mu}\in\ \Plambda\otimes
\mathbb{P}_{\mu},\, \pi_{\lambda+\mu}\in\ \mathbb{P}_{\lambda+\mu}$. The following obtains by combining lemmata \ref{unique}, \ref{unique2} and \ref{pix2}.\\

\begin{theorem} The crystals generated by $\pi_{\lambda}\otimes \pi_{\mu}$ and $\pi_{\lambda+\mu}$ are
isomorphic.
\end{theorem}
\\\\
{\bf Remark.} The crystals generated by $\pi_{\lambda}\otimes \pi_{\mu}$ and $\pi_{\lambda+\mu}$ viewed as
paths in $\mathbb{P}$ need not be isomorphic. For example, take $i\in\ I^{im}$, with $a_{ii}=-1$ and
$\lambda\in\ P^+$, with $\alpha_i^{\vee}(\lambda)=1$. Then by definition, $e_i\pi_{\lambda}=0$ in
$\mathbb{P}$ and so $e_i(\pi_{\lambda}\otimes \pi_{\lambda})=0$, by the tensor product rules. On the other
hand, $e_i\pi_{2\lambda}\neq 0$, again by lemma \ref{fimgls}. Of course $e_i\pi_{2\lambda}=0$ in
$\mathbb{P}_{2\lambda}$.


\section{Crystal Embedding Theorem}
We proved that the family of path crystals $\{\mathbb{P}_{\lambda}\,|\,\lambda\in\ P^+\}$ is closed. We will
now define the limit $\mathbb{P}_{\infty}$ of the family $\{\mathbb{P}_{\lambda}\,|\,\lambda\in\ P^+\}$ and
show that it is isomorphic to $B(\infty)$ (see theorem \ref{Binfty}).

\subsection{The limit $\mathbb{P}_{\infty}$}
\subsubsection{}
Let $\lambda, \mu\in\ P^+$ be two dominant weights and let $\pi \in\ \Plambda$. Denote by $\psi_{\lambda,
\lambda+\mu}$ the application $\psi_{\lambda, \lambda+\mu} : \Plambda\rightarrow \mathbb{P}_{\lambda}\otimes
\mathbb{P}_{\mu}$ which sends $\pi$ to $\pi\otimes \pi_{\mu}$. \\

\begin{lemma}\label{inductive} The application $\psi_{\lambda, \lambda+\mu}$ commutes with the $e_i, \,i\in\ I$, $\wt (\pi_{\lambda}\otimes \pi)=\wt\pi+\mu$ and if $f_i\pi\neq 0$, then $f_i\psi_{\lambda, \lambda+\mu}(\pi)=\psi_{\lambda,
\lambda+\mu}(f_i\pi)$. Thus $\psi_{\lambda, \lambda+\mu}$ is a crystal embedding up to translation of weight
by $\mu$.
\end{lemma}
\begin{proof} One has that $\varepsilon_i(\pi_{\mu})=0\leq \varphi_i(\pi)$, for all $i\in\ I$. If
$\varphi_i(\pi)=0$, then $f_i\pi=0$, for all $i\in\ I$, by section \ref{categoryB} (3). If $i\in\ I^{im}$ and
$\varphi_i(\pi)=0$, then $e_i\pi=0$ by lemma \ref{conditionfore}. Then apply the tensor product rules of
sections \ref{crystaltensorproduct} and \ref{tensorinB}.
\end{proof}

\subsubsection{}\label{limit}
As a set $\mathbb{P}_{\infty}$ is the inductive limit of the $\mathbb{P}_{\lambda}$ with respect to the above
embeddings. We will endow $\mathbb{P}_{\infty}$ with a crystal structure. Let  $\pi\in\ \mathbb{P}_{\infty}$,
then $\pi\in\ \Plambda$ for some $\lambda\in\ P^+$. Define $e_i\pi$ in $\mathbb{P}_{\infty}$ as $e_i\pi$ in
$\Plambda$. Define $f_i\pi$ again as in $\Plambda$ but we put $f_i\pi=0$ only if $f_i\psi_{\lambda,
\lambda+\mu}(\pi)=0$ for all $\mu\in\ P^+$. Finally, we define the weight of $\pi$ to be $-\mu$ if $\pi\in\
(\mathbb{P}_{\lambda})_{\lambda-\mu}$. This is clearly well defined. There is a unique element of weight
zero, since $\pi_{\lambda}$ is also unique in $(\Plambda)_{\lambda}$. We will denote this element by
$\pi_{\infty}$. It satisfies $e_i\pi_{\infty}=0$ for all $i\in\ I$. Notice that we may now forget about the
path crystal and consider any closed family of highest weight crystals $\{B(\lambda)|\lambda\in\ P^+\}$.

\subsection{The embedding theorem}
\subsubsection{} Recall the elementary crystals $B_i,\, i\in\ I$ defined in section \ref{elementary}.\\

\begin{theorem}\label{embedding} For all $i\in\ I$ there exists a unique strict embedding $\Psi_i : \mathbb{P}_{\infty}\longrightarrow \mathbb{P}_{\infty}\otimes B_i$, sending $\pi_{\infty}$ to $\pi_{\infty}\otimes b_i(0)$.
\end{theorem}
\begin{proof} Fix $i\in\ I$ and $f\in\ \mathcal{F}$. Call $f'$ a submonomial of $f$, if $f'$ obtains from $f$ by erasing
some of its factors. We say $f'$ is an $i$-submonomial of $f$ if it obtains by erasing some of the factors
$f_i$ in $f$. Let $\lambda\in\ P^+$ be such that $\alpha_i^{\vee}(\lambda)=0$ and $\alpha_j^{\vee}(\wt
f'\pi_{\lambda})>0$ for all $j\in\ I\setminus \{i\}$ and for all submonomials $f'$ of $f$. Let $\mu\in\ P^+$
be such that $\alpha_j^{\vee}(\mu)=0$, for all $j\in\ I\setminus \{i\}$. We will show that there exists an
integer $m\geq 0$ and an $i$-submonomial $f''$ of $f$ such that $f(\pi_{\lambda}\otimes
\pi_{\mu})=f''\pi_{\lambda}\otimes f_i^m\pi_{\mu}$. We argue by induction on the length of $f$.

For $f=\id$ the assertion is obvious. Let it be true for $f$ and set $f(\pi_{\lambda}\otimes
\pi_{\mu})=f''\pi_{\lambda}\otimes f_i^m\pi_{\mu}$. First notice that

$$f_i(f''\pi_{\lambda}\otimes f_i^{m}\pi_{\mu}) =\left\{
\begin{array}{ll}f_if''\pi_{\lambda}\otimes f_i^m\pi_{\mu}, &\varphi_i(f''\pi_{\lambda})>\varepsilon_i(f_i^m\pi_{\mu}),\\
f''\pi_{\lambda}\otimes f_i^{m+1}\pi_{\mu}, &\varphi_i(f''\pi_{\lambda})>\varepsilon_i(f_i^m\pi_{\mu}),
\end{array}
\right.
$$
which is of the required form.

Now for $j\in\ I\setminus \{i\}$ by assumption we have that $\varphi_j(f''\pi_{\lambda})\geq
\alpha_j^{\vee}(\wt f''\pi_{\lambda})>0$. On the other hand, $\varepsilon_j(f_i^m\pi_{\mu})=0$. Indeed, for
$j\in\ I^{im}$ this follows by definition. For $j\in\ I^{re}$, since $\mu-m\alpha_i+\alpha_j\notin \nu-Q^+$,
one has that $e_jf_i^m\pi_{\mu}=0$, hence by normality $\varepsilon_j(f_i^m\pi_{\mu})=0$. Then
$f_j(f''\pi_{\lambda}\otimes f_i^m\pi_{\mu}) = f_jf''\pi_{\lambda}\otimes f_i^m\pi_{\mu}$ which is also of
the required form.

By section \ref{limit} :
\[\varphi_i(f''\pi_{\lambda}) =
\varepsilon_i(f''\pi_{\lambda})+\alpha_i^{\vee}(\wt f''\pi_{\lambda}) =
\varepsilon_i(f''\pi_{\infty})+\alpha_i^{\vee}(\wt f''\pi_{\infty})+\alpha_i^{\vee}(\lambda)=
\varphi_i(f''\pi_{\infty}),\] which means that $\varphi_i(f''\pi_{\lambda})$ is independent of $\lambda$.

Finally one has $\varepsilon_i(b_i(-m))=\varepsilon_i(f_i^m\pi_{\mu})$ (and equal to $m$ for real indices and
$0$ for imaginary ones) and $\varepsilon_j(b_i(-m))=-\infty<\varphi_j(f''\pi_{\lambda})$, so that if
$f(\pi_{\lambda}\otimes \pi_{\mu})=f''\pi_{\lambda}\otimes f_i^m\pi_{\mu}$ then also $f(\pi_{\infty}\otimes
b_i(0))=f''\pi_{\infty}\otimes b_i(-m)$.
\end{proof}

\subsubsection{}
\begin{corollary} The crystal $\mathbb{P}_{\infty}$ is isomorphic as a crystal to $B(\infty)$.
\end{corollary}
\begin{proof} Notice that $\mathbb{P}_{\infty}$ has properties (1)-(4) of definition \ref{Binfty}.
Indeed the first three follow by construction and (4) follows by theorem \ref{embedding}. Then the assertion
follows by the uniqueness of $B(\infty)$.
\end{proof}

\section{The Character Formula}

\subsection{Weyl-Kac-Borcherds character formula}\label{wkb}

\subsubsection{}\label{not}
Assume the Borcherds-Cartan matrix $A$ to be symmetrizable. Recall $\rho\in\ \lieh^*$ of section
\ref{Wgeometry}. Let $\mathcal{P}(\Pi_{im})$ denote the set of all finite subsets $F$ of $\Pi_{im}$ such that
$\alpha_i^{\vee}(\alpha_j)=0$ for all $\alpha_i, \alpha_j\in\ F$. For all $\lambda\in\ P^+$ set
$$\mathcal{P}(\Pi_{im})^{\lambda}=\{F\in\ \mathcal{P}(\Pi_{im})\,|\,\alpha_i^{\vee}(\lambda)=0, \quad \mbox{for all}\quad \alpha_i\in\ F\}.$$
Given $F\in\ \mathcal{P}(\Pi_{im})$, let $|F|$ denote its cardinality and $s(F)$ the sum of its elements.
Then the character of the unique irreducible integrable highest weight module of $\lieg$ of highest weight
$\lambda\in\ P^+$ is given by the following formula known as the {\it Weyl-Kac-Borcherds character formula} :
\begin{equation}\label{character formula} \operatorname{char}\,V(\lambda) = \frac{\sum\limits_{w\in\ W}\sum\limits_{F\in\ \mathcal{P}(\Pi_{im})^{\lambda}}(-1)^{\ell(w)+|F|}e^{w(\lambda+\rho-s(F))}}{\sum\limits_{w\in\ W}\sum\limits_{F\in\
\mathcal{P}(\Pi_{im})}(-1)^{\ell(w)+|F|}e^{w(\rho-s(F))}}.
\end{equation}

\subsubsection{}
{\bf Remark.} It is not known if the above holds when $A$ fails to be symmetrizable. For $\Pi=\Pi_{re}$ of
finite cardinality, Kumar \cite{Ku} and Mathieu \cite{M} independently showed that the right hand side is the
correct character formula for the largest integrable quotient of the Verma module of highest weight
$\lambda$.

\subsubsection{}
Drop the assumption that the Borcherds-Cartan matrix $A$ is symmetrizable. Notice that the right hand side of
(\ref{character formula}) is still defined in this case. Define the character of $\mathbb{P}_{\lambda}$ by
$$\operatorname{char} \,\mathbb{P}_{\lambda} =
\sum\limits_{\pi\in\ \mathbb{P}_{\lambda}}e^{\pi(1)}.$$ Our main result is the following :\\

\begin{theorem}\label{character} The character of $\mathbb{P}_{\lambda}$ is given by the Weyl-Kac-Borcherds formula, that is to say the right hand side of (\ref{character formula}).
\end{theorem}
\\\\
The rest of the section is devoted to the proof of this theorem.

\subsection{The action of the Weyl group}\label{Weyl}
For all $i\in\ I^{re}$ define $\tilde r_i$ on $\pi\in\ \Plambda$ as follows :
$$\tilde r_i\pi = \left\{
\begin{array}{lll}f_i^{\alpha_i^{\vee}(\pi(1))}\pi, &\mbox{if}&\alpha_i^{\vee}(\pi(1))\geq 0,\\
e_i^{-\alpha_i^{\vee}(\pi(1))}\pi, &\mbox{if}&\alpha_i^{\vee}(\pi(1))\leq 0.
\end{array}
\right.$$ Then by \cite[Section 8]{L2} one has that $r_i\mapsto \tilde r_i$ extends to a representation
$W\rightarrow \operatorname{End}_{\mathbb{Z}}\Plambda$ and $w(\pi(1))=(w\pi)(1)$. Here we note that
$\mathbb{P}=\Pi_{\operatorname{int}}$ in the sense of \cite{L2} and the root operators $e_i,\,f_i,\, i\in\
I^{re}$ are defined as in \cite{L2}.

\subsection{The Kashiwara function}\label{Kashiwara}
Recall the crystal $B_J(\infty)$ of section \ref{BJcrystal} and that any element in $B_J(\infty)$ takes the
form
\begin{equation}\label{BJ} b=\dots\otimes b_{i_2}(-m_2)\otimes b_{i_1}(-m_1),
\end{equation}
with $m_k\in\ \mathbb{N}$ and $m_k=0$ for $k\gg 0$.
\subsubsection{}
Define the Kashiwara functions on $B_J(\infty)$ through
\begin{equation}r_i^k(b)=\varepsilon_i(b_{i_k}(-m_k))-\sum\limits_{j>k}\alpha_i^{\vee}(\wt
b_{i_j}(-m_j))=\varepsilon_i(b_{i_k}(-m_k))+\sum\limits_{j>k}m_ja_{i,i_j},
\end{equation}
noting that this sum is finite since $m_j=0$ for $j\gg 0$. Observe that $r_i^k(b)\in\ \{0, -\infty\}$ for
$k\gg 0$. Set $R_i(b)=\max_k\{r_i^k(b)\}$. From the definition of $J$ it follows that $R_i(b)\geq 0$ for all
$i\in\ I$ and all $b\in\ B_J(\infty)$, and $R_i(b)=0$ for all $i\in\ I^{im}$ and all $b\in\ B_J(\infty)$.
Note that if $R_i(b)=r_i^{k_0}(b)$ for some $k_0$, then $i_{k_0}=i$.

\subsubsection{}
The Kashiwara function determines at which place $e_i$ (resp. $f_i$) enters when computing $e_ib$ (resp.
$f_ib$). Let $\ell_i(b)$ (resp.
$s_i(b)$) be the largest (resp. smallest) value of $k$ such that $r_i^k(b)=R_i(b)$. Exactly as in the Kac-Moody case one has the following lemma :\\

\begin{lemma} \label{rfe}For all $b\in\ B_J(\infty)$ one has :
\begin{enumerate}
\item $\varepsilon_i(b)=0$ if and only if $R_i(b)=0$ for all $i\in\ I^{re}$ and $\varepsilon_i(b)=0$ for $i\in\
I^{im}$.
\item For all $i\in\ I$, $f_i$ enters at the $s_i(b)$th place.
\item For all $i\in\ I^{re}$, $e_i$ enters at the $\ell_i(b)$th place.
\item For all $i\in\ I^{im}$, $e_i$ enters at the $s_i(b)$th place.
\end{enumerate}
\end{lemma}

\noindent {\bf Remark.} It can happen that $\ell_i(b)=+\infty$ for $i\in\ I^{re}$, but then simply $e_ib=0$.
\subsubsection{}
\begin{lemma}\label{2elts} Let $b, b'\in\ B_J(\infty)$ be such that $f_ib=f_jb'$ for $i,\,j\in\ I^{im}$ and $i\neq j$. Then $f_i,\,
f_j$ commute and there exists $b''\in\ B_J(\infty)$ such that $b=f_jb''$.
\end{lemma}
\begin{proof} Write $b, b'$ as in (\ref{BJ}) with $m_k$ replaced by $m'_k$ for the latter. Suppose that $f_i$ enters $b$ at the $\ell$th place
and $f_j$ enters $b'$ at the $\ell'$th place. Since $i\neq j$ we have that $\ell\neq \ell'$. We can assume
that $\ell'<\ell$ interchanging $i, j$ if necessary. Note that $f_ib=f_jb'$ forces $m'_{\ell}=m_{\ell}+1>0$.

We have
$r_j^{\ell'}(b')=\varepsilon_j(b_{i_{\ell}}(-m'_{\ell}))+\sum\limits_{s>\ell'}m'_s\alpha_i^{\vee}(\alpha_{i_s})=R_j(b')\geq
0$. Since $j=i_{\ell'}$ one has that $\varepsilon_j(b_{i_{\ell}}(-m'_{\ell}))=0$. On the other hand $m'_s\geq
0$ and $\alpha_i^{\vee}(\alpha_{i_s})\leq 0$ for all $s>\ell'$, forcing $m_s'\alpha_i^{\vee}(\alpha_{i_s})=0$
for all $s>\ell'$. In particular, since $i_{\ell}=i$, $\ell>\ell',\, m_{\ell}'>0$ we obtain that
$\alpha_j^{\vee}(\alpha_i)=0$, that is $a_{ji}=0$ and so $a_{ij}=0$.

Take $c\in\ B_J(\infty)$. Since $\alpha_j^{\vee}(\alpha_i)=0$, one has $r_j^k(c)=r_j^k(f_ic)$, for all $c\in\
B_J(\infty)$, and all $k\in\ \mathbb{N}^+$. Then $s_i(f_jc)=s_i(c)$. Similarly $s_j(f_ic)=s_j(c)$ and so
$f_if_jc=f_jf_ic$, as required.

On the other hand, since $f_ib=f_jb'$ we have that $m_{\ell'}=m_{\ell'}'+1>0$. By lemma \ref{rfe}, $e_j$
enters $b$ at the $\ell'$th place and $e_jb\neq 0$. Set $e_jb=b''$, then $b=f_jb''$. Yet
$f_jb'=f_ib=f_if_jb''=f_jf_ib''$ and so $b'=f_ib''$.
\end{proof}

\subsubsection{}
\begin{corollary}\label{2paths} Let $\pi,\pi'$ be two paths in $\Plambda$ for $\lambda\in\ P^+$ such that $f_i\pi=f_j\pi'$ for $i,\,j\in\ I^{im}$ and $i\neq j$. Then $f_i,\,
f_j$ commute and there exists $\pi''\in\ \Plambda$ such that $\pi=f_j\pi''$.
\end{corollary}
\begin{proof} Embed $\Plambda$ in $B_J(\infty)$ : $$\Plambda\stackrel{\psi_1}{\hookrightarrow} B(\infty)\stackrel{\psi_2}{\hookrightarrow}
B_J(\infty).$$ Then if we assume $f_i\pi=f_j\pi'\neq 0$ in $\Plambda$ $f_i,\, f_j$ commute with $\psi_1$ and
$\psi_2$ so that $f_i\psi_2\psi_1(\pi)=f_j\psi_2\psi_1(\pi')$. The assertions follow by lemma \ref{2elts}. We
note here that $\psi_1$ does not in general commute with $f_i,\, f_j$ that is why we have to assume
$f_i\pi=f_j\pi'\neq 0$ (see lemma \ref{inductive}).
\end{proof}

\subsection{Proof of theorem \ref{character}}

\subsubsection{}
Recall section \ref{wkb}; we will show that the character of $\Plambda$ is given by the Weyl-Kac-Borcherds
formula. We need to show that :
\begin{equation}\label{char1}\sum\limits_{\pi\in\ \Plambda}\sum\limits_{w\in\ W}\sum\limits_{F\in\ \mathcal{P}(\Pi_{im})}(-1)^{\ell(w)+|F|}e^{w(\rho-s(F))+\pi(1)} = \sum\limits_{w\in\ W}\sum\limits_{F\in\ \mathcal{P}(\Pi_{im})^{\lambda}}(-1)^{\ell(w)+|F|}e^{w(\lambda+\rho-s(F))}.
\end{equation}

\subsubsection{}
For all $\mu\in\ P$ set $O(\mu)=\{(w, F, \pi)\in\ W\times \mathcal{P}(\Pi_{im})\times
\mathbb{P}_{\lambda}\,|\, w(\rho-s(F))+\pi(1)=\mu\}$. By section \ref{Weyl} we have an action of $W$ on
$O(\mu)$ by $w(w', F, \pi)=(ww', F, w\pi)$, where $wO(\mu)=O(w\mu)$. Moreover, since
$\ell(ww')=\ell(w)+\ell(w') \mod 2$, the sum
$$S(\mu):=\sum\limits_{(w, F, \pi)\in\ O(\mu)}(-1)^{\ell(w)+|F|}$$ satisfies $S(w\mu)=(-1)^{\ell(w)}S(\mu)$.
Now the left hand side of (\ref{char1}) becomes
\begin{equation}\label{reddom} \sum\limits_{\mu\in\ P}S(\mu)e^{\mu}=\sum\limits_{w\in\ W}\sum\limits_{\mu\in\ P^+}(-1)^{\ell(w)}S(\mu)e^{w\mu}.
\end{equation}
Then we can assume $\mu :=w(\rho-s(F))+\pi(1)$ to be dominant and in this case it remains to show that
$S(\mu)=0$, unless $O(\mu)=\{\id\}\times \mathcal{P}(\Pi_{im})\times \{\pi_{\lambda}\}$.

\subsubsection{}
Since $\mu$ is dominant and $t\mapsto \pi(t)$ is continuous, either
\begin{enumerate}
\item there exists some $t\in\ [0, 1]$ such that $w(\rho - s(F)) +\pi(t)$ is dominant but not regular or
\item $w(\rho - s(F)) +\pi(t)$ is regular and dominant for all $t\in\ [0, 1]$.
\end{enumerate}

Thus define
$$O_1(\mu) : = \{(w, F, \pi)\in\ O(\mu)\, | w(\rho - s(F))
+\pi(t) \, \mbox{is dominant but not regular for some}\, t\in\ [0, 1]\},$$ and
$$O_2(\mu) : = \{(w, F, \pi)\in\ O(\mu)\, | w(\rho - s(F))
+\pi(t) \, \mbox{is dominant and regular for all}\, t\in\ [0, 1]\}.$$

\subsubsection{}
In case (1) exactly as in \cite[Theorem 9.1]{L2} we obtain that $$\sum\limits_{(w, F, \pi)\in\
O_1(\mu)}(-1)^{\ell(w)+|F|}e^{\mu}=0.$$ In case (2), $w(\rho-s(F))+\pi(t)$ being dominant at $t=0$, implies
$w=\id$. Thus we define
$$\tilde O_2(\mu) := \{(F, \pi)\in\  \mathcal{P}(\Pi_{im})\times \Plambda\,|\, (\id, F, \pi)\in\ O_2(\mu)\}.$$

The formula we have to prove becomes :

\begin{equation}\label{f2}\sum\limits_{\mu\in\ P^+}\sum\limits_{(F, \pi)\in\ \tilde O_2(\mu)}(-1)^{|F|}e^{\rho-s(F)+\pi(1)}=\sum\limits_{F\subset \mathcal{P}(\Pi_{im})^{\lambda}}(-1)^{|F|}e^{\rho-s(F)+\lambda}.
\end{equation}

\subsubsection{}
For all $(F, \pi)\in\ \mathcal{P}(\Pi_{im})\times \Plambda$ set
\begin{equation}\label{s} S(F, \pi):=\{\alpha_i\in\ \Pi_{im}\setminus F \, |\, \alpha_i^{\vee}(s(F))=0 \quad \mbox{and}\quad e_i\pi\neq
0\}.
\end{equation}
Take $i, j\in\ I^{im}$ distinct. Notice that if $\alpha_i,\,\alpha_j\in\ S(F, \pi)$, then $a_{ij}=a_{ji}=0$.
In particular, $F\cup S(F, \pi)\in\ \mathcal{P}(\Pi_{im})$. Indeed, since $e_i\pi,\, e_j\pi\neq 0$, one has
$\pi = f_i\pi_1=f_j\pi_2$, for $\pi_1,\, \pi_2\in\ \Plambda$, and the assertion follows by lemma
\ref{2paths}. We call a pair $(F, \pi)\in\ \mathcal{P}(\Pi_{im})\times \Plambda$ minimal, if $S(F,
\pi)=\emptyset$.

For any subset $S = \{\alpha_{i_1}, \alpha_{i_2},\dots \alpha_{i_k}\}\subset \mathcal{P}(\Pi_{im})$, set
$f_S:=f_{i_1}f_{i_2}\cdots f_{i_k}$ and similarly $e_S:=e_{i_1}e_{i_2}\cdots e_{i_k}$. Notice that since the
$f_{i_j}$ (resp. $e_{i_j}$) mutually commute, the monomial $f_S$ (resp. $e_S$) does not depend on the order
of the indices. Suppose that $\pi\in\ \Plambda$ satisfies $e_i\pi\neq 0$, for all $\alpha_i\in\ S$. Then
$e_S\pi\neq 0$. Indeed this follows from lemmata \ref{rfe} and \ref{inductive}. Again if $f_i\pi\neq 0$ for
all $\alpha_i\in\ S$, then $f_S\pi\neq 0$. This follows from lemma \ref{star} (2) and section \ref{iff}.

\subsubsection{}
For all $(F, \pi)\in\ \mathcal{P}(\Pi_{im})\times \Plambda$, set $(F_0, \pi_0)=(F\cup S(F, \pi),
e_{S(F,\pi)}\pi)$. Clearly $(F_0, \pi_0)$ is minimal. For a minimal element $(F_0, \pi_0)$ define
$F_0^{\pi_0}:=\{S\subset F_0\,|\,\forall \, \alpha_i\in\ S, f_i\pi\neq 0\}$. Then set
$$\Omega(F_0, \pi_0) = \{(F_0\setminus S, f_S\pi_0)| S\in\ F_0^{\pi_0}\}.$$
The following is straightforward.\\

\begin{lemma} The only minimal element in $\Omega(F_0, \pi_0)$ is $(F_0, \pi_0)$. Moreover, if $(F, \pi)=(F_0\setminus S, f_S\pi_0)\in\ \Omega(F_0, \pi_0)$, for $S\in\ F_0^{\pi_0}$, then $S(F,
\pi)=S$.
\end{lemma}
\\

An immediate consequence of the above is that for any two minimal elements $(F_0, \pi_0)\neq (F_0', \pi_0')$,
one has $\Omega(F_0, \pi_0)\cap \Omega(F_0', \pi_0')=\emptyset$.
\\

\noindent {\bf Remark.} Note that for all $(F, \pi)\in\ \Omega(F_0, \pi_0)$, the weight $-s(F)+\pi(1)$ is
fixed (but it does not uniquely define $\Omega(F_0, \pi_0)$!).

\subsubsection{}\label{f1f2}
We show in section \ref{last} that if $\Omega(F_0, \pi_0)\cap \tilde O_2(\mu)\neq \emptyset$, then
$\Omega(F_0, \pi_0)\subset \tilde O_2(\mu)$. Then if we set $\Omega(\mu)=\{(F_0, \pi_0)\,|\,\Omega(F_0,
\pi_0)\subset \tilde O_2(\mu) \}$ we have

\begin{equation} \label{O2mu} \tilde O_2(\mu) = \coprod\limits_{(F_0, \pi_0)\in\ \Omega(\mu)}\Omega(F_0, \pi_0).
\end{equation}

Admit \ref{last}, so then (\ref{O2mu}) holds. We have :
\begin{equation}\label{f32}\sum\limits_{(F, \pi)\in\ \tilde O_2(\mu)}(-1)^{|F|}=\sum\limits_{(F_0, \pi_0)\in\ \Omega(\mu)}\left(\sum\limits_{(F, \pi)\in\ \Omega (F_0, \pi_0)}(-1)^{|F|}\right).
\end{equation}
We will compute the following sum :
\begin{equation}\label{onlyim} \Sigma:=\sum\limits_{(F, \pi)\in\ \Omega (F_0, \pi_0)}(-1)^{|F|}e^{\rho-s(F)+\pi(1)}.
\end{equation}

\subsubsection{}
\begin{lemma}\label{sum} If $|\Omega(F_0, \pi_0)|>1$, then the sum $\Sigma$ above is zero.
\end{lemma}
\begin{proof}
Write $F_0$ as $F_0=F_0'\sqcup F_0''$, where $F_0' :=\{\alpha_i\in\ F_0\,|\, \alpha_i^{\vee}(\pi(1))\neq 0\}$
(equivalently, by \ref{fimgls}, $F_0' :=\{\alpha_i\in\ F_0\,|\, f_i\pi\neq 0\}$). Our hypothesis that
$|\Omega(F_0, \pi_0)|>1$ implies that $F_0'\neq \emptyset$. Set $|F_0'|=n\geq 1$. Then the cardinality of
$\Omega(F_0, \pi_0)$ is equal to the number of subsets of $F_0'$. Moreover, the coefficient of $e^{\rho
-s(F)+\pi(1)}$ in $\Sigma$ is
$$(-1)^{|F_0''|}((-1)^n+(-1)^{n-1} \left(\begin{array}{cc}n\\1
\end{array}
\right) +(-1)^{n-2} \left(\begin{array}{cc}n\\2
\end{array}\right)
+\cdots +(-1) \left(\begin{array}{cc}n\\n-1
\end{array}\right)+1)=0.
$$
\end{proof}

\subsubsection{}
\begin{lemma} \label{sums} If $|\Omega(F_0, \pi_0)|=1$, then $\pi_0=\pi_{\lambda}$ and $F_0\in\
\mathcal{P}(\Pi_{im})^{\lambda}$.
\end{lemma}
\begin{proof}
Let $\Omega(F_0, \pi_0)= \{(F_0, \pi_0)\}$ be a singleton and recall section \ref{iff}. Then for all
$\alpha_j\in\ F_0$, one has that $\alpha_j^{\vee}(\pi_0(1))=0$. In particular, if $\pi_0=f_i\pi$ for some
$\pi\in\ \Plambda$ and $i\in\ I$, then $\alpha_i^{\vee}(s(F_0))=0$ and $\alpha_j^{\vee}(\pi(1))=0$, for all
$\alpha_j\in\ F_0$.

Assume that $\pi_0\neq \pi_{\lambda}$. Then we may write $\pi_0=f_i\pi$ as above. Suppose that $i\in\
I^{im}$, then $e_i \pi_0\neq 0$ and by the above remark $\alpha_i^{\vee}(s(F_0))=0$. By the minimality of
$(F_0, \pi_0)$, this implies that $\alpha_i\in\ F_0$. Yet $\pi_0=f_i\pi$ implies $f_i\pi_0\neq 0$ by lemma
\ref{computationofhfi} and so $(F_0\setminus \{i\}, f_i\pi_0)\in\ \Omega(F_0, \pi_0)$, which contradicts the
hypothesis.

Let now $i\in\ I^{re}$. Then $\alpha_i^{\vee}(\rho-s(F_0)+
\pi_0(t))=1+\alpha_i^{\vee}(\pi_0(t))=1+h_i^{\pi_0}(t)$. But $e_i \pi_0\neq 0$ which means (by definition)
that $h_i^{\pi_0}(t)$ takes integral values $\leq -1$. Hence $\rho-s(F_0)+ \pi_0(t)$ is not regular for all
$t\in\ [0, 1]$, again a contradiction.

We obtain that $\pi_0=\pi_{\lambda}$ and $F_0\in\ \mathcal{P}(\Pi_{im})^{\lambda}$.
\end{proof}
By lemmata \ref{sum}, \ref{sums} the only remaining terms in the left hand side of (\ref{f2}) is the right
hand side. Thus to complete the proof of theorem \ref{character} it remains to prove (\ref{O2mu}). As we
noted in section \ref{f1f2}, this follows by the lemma below.

\subsubsection{}\label{last}
\begin{lemma} Let $\Omega(F_0, \pi_0)\cap \tilde O_2(\mu)\neq \emptyset$. Then $\Omega(F_0, \pi_0)\subset \tilde O_2(\mu)$.
\end{lemma}
\begin{proof}
Fix $(F, \pi)\in\ \Omega(F_0, \pi_0)\cap \tilde O_2(\mu)$.

Assume that $(F\cup \{\alpha_i\}, e_i\pi)\in\ \Omega(F_0, \pi_0)$. This means that $\alpha_i\notin F$,
$\alpha_i^{\vee}(s(F))=0$ and $e_i\pi\neq 0$. We show that $(F\cup \{\alpha_i\}, e_i\pi)\in\ \tilde
O_2(\mu)$.

By definition of $e_i\pi$, there exists a piecewise linear function $c(t)$ with $0\leq c(t)\leq 1$ for all
$t\in\ [0, 1]$ such that $e_i\pi(t)=\pi(t)+c(t)\alpha_i$. Then since $\rho-s(F)+\pi(t)$ is regular and
dominant for all $t\in\ [0, 1]$ and $\alpha_i$ is anti-dominant, we obtain that $\rho-s(F\cup
\{\alpha_i\})+e_i\pi(1) = \rho -s(F)-\alpha_i+\pi(t)+c(t)\alpha_i=\rho-s(F)+\pi(t)+(c(t)-1)\alpha_i$ is also
regular and dominant for all $t\in\ [0, 1]$, as required.

Now suppose that $(F\setminus \{\alpha_i\}, f_i\pi)\in\ \Omega(F_0, \pi_0)$. It follows that $\alpha_i\in\ F$
and $f_i\pi\neq 0$. We show that $(F\setminus \{\alpha_i\}, f_i\pi)\in\ \tilde O_2(\mu)$.

Set $F'=F\setminus \{\alpha_i\}$, then $F=F'\cup \{\alpha_i\}$. By assumption, $M(t) :=
\rho-s(F)+\pi(t)=\rho-s(F')-\alpha_i+\pi(t)$ is regular and dominant for all $t\in\ [0, 1]$. We need to show
that $M'(t):=\rho-s(F)+\alpha_i+(f_i\pi)(t)=\rho-s(F')+(f_i\pi)(t)$ is regular and dominant for all $t\in\
[0,1]$. Now for $t\in\ [f_-^i(\pi), 1]$ one has $f_i\pi(t)=\pi(t)-\alpha_i$ and so $M'(t)=M(t)$, hence
$M'(t)$ is regular and dominant for all $t\in\ [f_-^i(\pi), 1]$.

Suppose that for some $t\in\ [0, f_-^i(\pi)]$, $M'(t) = \rho-s(F')+(f_i\pi)(t)$ is not regular. This means
that there exists $j\in\ I^{re}$ such that $\alpha_j^{\vee}(M'(t))=0$, for some $t\in\ [0, f_-^i(\pi)[$. In
this region, $(f_i\pi)(t)=r_i\pi(t)$, hence
\begin{equation}h_j(t):=\alpha_j^{\vee}(M'(t)) =
h_j^{f_i\pi}(t)+\alpha_j^{\vee}(\rho-s(F'))=\alpha_j^{\vee}(\pi(t))-\alpha_i^{\vee}(\pi(t))a_{ji}+\alpha_j^{\vee}(\rho-s(F'))=
0,
\end{equation}
for some $t\in\ [0, f_-^i(\pi)]$. On the other hand $h_j(0) = \alpha_j^{\vee}(\rho-s(F'))>0$ and
$h_j(f_-^i(\pi))=\alpha_j^{\vee}(M'(f_-^i(\pi)))=\alpha_j^{\vee}(M(f_-^i(\pi)))>0$, hence the function $h_j$
attains a local minimum at some $t_0\in\ ]0, f_-^i(\pi)[$ and consequently  $h_j^{f_i\pi}$ attains a local
minimum at $t_0$.

Let $\pi = (\lambda_1, \lambda_2, \dots, \lambda_s;0, a_1, a_2\dots, a_s=1)$ and recall proposition
\ref{actionfim}, choosing $p$ as defined there. One has $a_{p-1}<f_-^i(\pi)\leq a_p$ and
$$f_i\pi=(r_i\lambda_1, r_i\lambda_2, \dots, r_i\lambda_p, \lambda_p, \dots, \lambda_s;0, a_1, \dots,
a_{p-1}, f_-^i(\pi), a_p, \dots, a_s=1).$$ By lemma \ref{integral}, we must have $t_0=a_k$ for some $k\leq
p-1$ and so either $\alpha_j^{\vee}(\lambda_k)\leq 0$ and $\alpha_j^{\vee}(\lambda_{k+1})>0$, or
$\alpha_j^{\vee}(\lambda_k)< 0$ and $\alpha_j^{\vee}(\lambda_{k+1})\geq 0$, depending on whether the minimum
at $t_0$ is right or left. Then by lemma \ref{simpleroot}, if
$\lambda_{k+1}\stackrel{\beta_t}{\leftarrow}\cdots \stackrel{\beta_1}{\leftarrow}\lambda_k$, we obtain
$\beta_{\ell}=\alpha_j$ for some $\ell$, with $1\leq \ell\leq t$. On the other hand, by proposition
\ref{actionfim}, $\alpha_i^{\vee}(\lambda_k)=\alpha_i^{\vee}(\lambda_{k+1})$, for all $k$, with $1\leq k\leq
p-1$ and so $\alpha_i^{\vee}(\beta_s)=0$ for all $s$, with $1\leq s\leq t$. In particular $a_{ij}=0$ and so
$h_j(t)=\alpha_j^{\vee}(M(t))$ which is strictly positive by assumption. This contradiction proves that
$M'(t)$ is regular for all $t\in\ [0, 1]$.
\end{proof}

\end{document}